\documentclass[10pt]{article}

\usepackage{hyperref}

\usepackage{enumitem}
\usepackage{enumitem} 
\usepackage{array} 
\usepackage{amsfonts}
\usepackage{amsmath}
\usepackage{amssymb}
\usepackage{graphicx}
\usepackage{color}
\usepackage{bm}
\usepackage[utf8]{inputenc}
\usepackage[T1]{fontenc}
\usepackage{nameref}
\usepackage{esint}
\usepackage{epstopdf}

\newtheorem{thm}{Theorem}[section]
\newtheorem{cor}[thm]{Corollary}
\newtheorem{lem}[thm]{Lemma}
\newtheorem{prop}[thm]{Proposition}
\newtheorem{defn}[thm]{Definition}
\newtheorem{rem}[thm]{Remark}
\numberwithin{equation}{section}

\newcommand{\D}{\mathcal D}
\newcommand{\dx}{\,{\rm d}x}
\newcommand{\dy}{\,{\rm d}y}
\newcommand{\dz}{\,{\rm d}z}

\newcommand{\ds}{\,{\rm d}s}
\newcommand{\dt}{\,{\rm d}t}
\newcommand{\dtau}{\,{\rm d}\tau}

\def\LL{\mathrm{L}} 
\def\supp{\mathrm{supp}} 
\def\sr{r^{\ast}}

\def\one{\mbox{1\hspace{-4.25pt}\fontsize{12}{14.4}\selectfont\textrm{1}}}

\newcommand{\ka}{\overline{\kappa}}
\newcommand{\ke}{k_{\varepsilon}}
\newcommand{\kb}{\underline{\kappa}}
\newcommand{\ta}{\overline{t}}
\renewcommand{\k}{\kappa}

\newcommand{\RR}{\mathbb{R}}
\newcommand{\NN}{\mathbb{N}}

\newcommand{\M}{\overline{V}}
\newcommand{\m}{\underline{V}}

\def\ee{\mathrm{e}} 
\DeclareMathOperator*{\osc}{osc}
\def\qed{\,\unskip\kern 6pt \penalty 500
\raise -2pt\hbox{\vrule \vbox to8pt{\hrule width 6pt
\vfill\hrule}\vrule}\par}
\definecolor{darkblue}{rgb}{0.05, .05, .65}
\definecolor{darkgreen}{rgb}{0.1, .65, .1}
\definecolor{darkred}{rgb}{0.8,0,0}

\begin{document}

\title{\textbf{Quantitative a Priori Estimates\\ for Fast Diffusion Equations with \\
Caffarelli-Kohn-Nirenberg weights.\\
Harnack inequalities and Hölder continuity. }\\[7mm]}

\author{\Large Matteo Bonforte$^{\,a}$
~and~ Nikita Simonov$^{\,b}$\\[3mm]
Departamento de Matem\'{a}ticas,\\
Universidad Aut\'{o}noma de Madrid,\\
Campus de Cantoblanco, 28049 Madrid, Spain.\\
}
\date{} 

\maketitle

\begin{abstract}
We study a priori estimates for a class of non-negative local weak solution to the weighted fast diffusion equation $u_t = |x|^{\gamma} \nabla\cdot (|x|^{-\beta} \nabla u^m)$, with $0 < m <1$ posed on cylinders of $(0,T)\times\RR^N$. The  weights $|x|^{\gamma}$ and $|x|^{-\beta}$, with $\gamma < N$ and $\gamma -2 < \beta \leq \gamma(N-2)/N$ can be both degenerate and singular and  need not belong to the class $\mathcal{A}_2$, a typical assumption  for this kind of problems. This range of parameters is optimal for the validity of a class of Caffarelli-Kohn-Nirenberg inequalities, which play the role of the standard Sobolev inequalities in this more complicated weighted setting.

The weights that we consider are not translation invariant and this causes a number of extra difficulties and a variety of scenarios: for instance, the scaling properties of the equation change when considering the problem around the origin or far from it. We therefore prove quantitative - with computable constants - upper and lower estimates for local weak solutions, focussing our attention where a change of geometry appears. Such estimates fairly combine into forms of Harnack inequalities of forward, backward and elliptic type.  As a consequence, we  obtain H\"older continuity of the solutions, with a quantitative (even if non-optimal) exponent. Our  results  apply to a quite large variety of solutions and problems. The proof of the positivity estimates requires a new method and represents the main technical novelty of this paper.

Our techniques are flexible and can be adapted to more general settings, for instance to a wider class of weights or to similar problems posed on Riemannian manifolds, possibly with unbounded curvature. In the linear case, $m=1$, we also prove quantitative estimates, recovering known results in some cases and extending such results to a wider class of weights.
\end{abstract}

\vspace{1cm}

\noindent {\sc Keywords. }Fast diffusion with weights; Parabolic Regularity; Positivity; Smoothing effects; Harnack inequalities; H\"older Continuity; Caffarelli-Kohn-Nirenberg inequalities.

\noindent{\sc Mathematics Subject Classification}. 35B45, 35B65, 35K55, 35K67, 35K65.

\vfill
\begin{itemize}
\item[(a)] E-mail:\texttt{~matteo.bonforte@uam.es }\;\;
Web-page:\texttt{~http://verso.mat.uam.es/\~\,matteo.bonforte/}
\item[(b)] E-mail:\texttt{~nikita.simonov@uam.es }\;\;
\end{itemize}

\newpage
\small

\tableofcontents

\normalsize

\newpage

\section{Introduction}
We investigate quantitative a priori estimates and regularity properties of nonnegative solutions to nonlinear singular diffusion equations with weights, possibly degenerate or singular, whose prototype is given by the following Weighted Fast Diffusion Equation
\begin{equation}\label{WFDE}\tag{WFDE}
u_t=|x|^\gamma\nabla\cdot\left(|x|^{-\beta}\nabla u^m\right)
\end{equation}
posed on a domain of $(0,+\infty)\times\RR^N$\,, with $N\ge 3$ and $m\in (0,1)$. We will always consider the following range of parameters, see also Figure \ref{fig.3} \,:
\begin{equation}\label{paramt.range}
\gamma < N\qquad\mbox{and}\qquad\gamma -2 < \beta \leq \frac{N-2}{N}  \gamma \,.\vspace{-2mm}
\end{equation}
This range of parameters is optimal for the validity of a family of the so-called Caffarelli-Kohn-Nirenberg inequalities \cite{Caffarelli-Kohn-Nirenberg-84}, then as follows: let $r^*:= 2(N-\gamma)/(N-(2+\beta))$,
\begin{equation} \label{w_sobolev0}
\left( \int_{\RR^N}{|f|^{\sr} |x|^{-\gamma} \dx} \right)^{1/\sr} \leq  \overline{S}_{\gamma, \beta} \left(\int_{\RR^N }{ |\nabla f|^2  |x|^{-\beta}\dx}  \right)^{1/2}\,,
\end{equation}
see Subsection \ref{sect.funct.ineq} for more details; these inequalities are deeply connected with the above WFDE, in its evolutionary or stationary version, see for instance \cite{BDMN2016b,BDMN2016a,Catrina-Wang-01,1406,DEL2015,DELM2015,DETT}; some further connection will be discussed and explored in this paper.

A priori estimates are the cornerstone of the theory of nonlinear partial differential equations: the main purpose of this paper is to prove precise quantitative local upper and lower bounds which combine into different forms of Harnack inequalities; as a consequence we also prove interior H\"{o}lder continuity for solutions to this class of equations with a (small) quantified exponent: the optimal H\"older exponent is not known. Indeed, in the case of the Cauchy problem some explicit (Barenblatt-type) solutions  are known to be only H\"older continuous at $x=0$, as we shall discuss later, see also \cite{BDMN2016b,BDMN2016a}.

The weights that we consider are not translation invariant and this causes a number of extra difficulties and a variety of scenarios that we explain in Subsection \ref{result.organization}. Roughly speaking, the scaling properties of the equation change from $R^{2+\beta-\gamma}$ to (a multiple of) $R^2$, when we are considering the problem around the origin or far from it, respectively. We focus on the cases in which the change of geometry plays a role: in the other cases, the results are essentially the same as the classical ones, cf. \cite{BV-ADV,DiBenGianVesBook}.

Our quantitative interior estimates are formulated for nonnegative local strong solutions, defined  in Subsection \ref{result.organization}. A number of interesting problems fall into our setting, for instance, the Cauchy problem on $\RR^N$, problems on Euclidean domains with different boundary condition (Dirichlet, Neumann, Robin, etc.), as well as the so-called ``large solutions'' (which tend to $+\infty$ at the boundary of the domain). Moreover, our estimates can be extended to a wider class of solutions, through lengthy but standard approximations.  We prove analogous results also in the linear case $m=1$, as we shall discuss below.

The above nonlinear equation has been introduced in the 80s by Kamin and Rosenau \cite{KR-CPAM81, KR-CPAM82, KR-JMPH}, to model heat propagation -or more generally singular/degenerate diffusion- in inhomogeneous media; the parabolic problem has been studied by many authors since then, mostly in the case $m\ge 1$ and with only one weight \cite{AP,BS2000,  DGP-SIAM,Eidus,KE-PAMS, GMPo13, GMP15, IS14, KKT, KRV-DCDS,KK-MA, K-SNS,  NR13, RV-CPAA09,RV-CPAA08,RV-JEMS06,RV-NHM,RV-JMPA,Surnachev-JDE, Surnachev-TMM, T-CPAM68, WNC-NonAnal}.

In the non-weighted case $\gamma=\beta=0$, the WFDE becomes the standard Fast Diffusion Equation (FDE) which has been intensively studied in the recent years by many authors: it is hopeless to give here a complete bibliography, hence we refer to the monographs \cite{JLVSmoothing, JLVPorousMedioum} and \cite{DaskaBook,DiBenGianVesBook} for a complete account, as well as for the physical relevance of the model. We just remark that our results hold also in the non-weighted case, and we recover the previous results with a different proof.

More recently, \ref{WFDE} has been investigated for its deep relation with the so-called Caffarelli-Kohn-Nirenberg inequalities, \cite{Caffarelli-Kohn-Nirenberg-84}; in particular, the intriguing  issue of symmetry/symmetry breaking, has attracted the attention of many prominent researchers,  \cite{BDMN2016b, BDMN2016a, DMN2017, Catrina-Wang-01, 1406,DEL2015, DELM2015, DETT,  Felli-Schneider-03}. The study of such problem partly relies on the study of the Cauchy problem for the WFDE on $\RR^N$ for which the regularity estimates of this paper are fundamental and were not present in the literature: sometimes an extra hypothesis had to be added to fix this issue. This happens for instance in \cite{BDMN2016b, BDMN2016a}, where the sharp asymptotic behaviour of solutions to the Cauchy problem for WFDE is studied: the regularity estimates proven here are indeed essential to ensure the validity of those results in full generality.

Lately, new connections between weighted parabolic equations and nonlinear diffusions on Riemannian manifolds were explored in \cite{BGV-JEE,BGV-ARMA,DEL2015, GM14,JLV-Hyp}. This intriguing connection motivates the present work, which makes a preliminary step towards understanding the behaviour of singular nonlinear diffusion on manifolds possibly with unbounded curvature; it has to be noticed that in this latter case the weights are locally regular and degenerate only at infinity, see for instance \cite{GMV-manifolds}.

Since the pioneering paper of Fabes, Kenig and Serapioni \cite{FKS}, weighted (degenerate or singular) elliptic and parabolic equations have been investigated in the linear case $m=1$, \cite{AB-JDE,AB-CM,XC-XRS, CF-AA84,  CS-AA87, CS-RSMUP85, CS-AMPA84, CS-CPDE84,I-NM,IM-MZ,IM-SNS,MadernaSalsa,Pinchover}; in many cases, the weights are assumed to belong to the ``natural'' (for second order differential operators) Muckenhoupt class $\mathcal{A}_2$: the reader may notice that the weight $|x|^{-\gamma}$ does not belong to $\mathcal{A}_2$ when $\gamma \in \left(-\infty, - N\right)$, a case that we consider here.  Our contribution in this direction are quantitative Harnack inequalities and H\"older continuity for weak solutions to linear equations  with measurable coefficients. Our results agree with the known results  \cite{GW2,GW, Moser,MoserCpam71} and extend those in some range of parameters.

\smallskip

\noindent\textbf{Ideas of the main results and organization of the paper. }
The behaviour of solutions to the non-weighted FDE presents strong differences  between  two ranges: \textit{good fast diffusion range} $m_c < m <1$ and \textit{very fast diffusion range} $0 < m \le m_c$: the \textit{critical exponent}  being $m_c=(N-2)/N$, see \cite{BV-ADV, JLVSmoothing}. We show here that  also solutions to the \ref{WFDE} behave quite differently in the two ranges; the \textit{critical exponent} $m_c$ now depends on the weights through $\gamma$ and  $\beta$, in the range given by \eqref{paramt.range}\vspace{-2mm}
\begin{equation}\label{critical.exponent}
m_c=\frac{N-\left(2+ \beta \right)}{N-\gamma}\in (0,1).\vspace{-2mm}
\end{equation}
Our first main result consists in quantitative upper bounds, see Theorems \ref{upper.easy} and \ref{local_upper_bounds} proven in Part I, which take the form of \textit{local smoothing effects}, that in a simplified form read\vspace{-1mm}
\begin{equation*}
\sup_{y\in B_{R}(0) }{u\left(t,y\right)} \leq \frac{\ka_1}{t^{(N-\gamma)\vartheta_p}}
\left[\int_{B_{2R}(0)}{|u_0(y)|^{p} \ |y|^{-\gamma} \ \dy} \right]^{(2+\beta-\gamma)\vartheta_p} + \ka_2 \left[\frac{t}{R^{2+\beta-\gamma}} \right]^{\frac{1}{1-m}},\vspace{-1mm}
 \end{equation*}
where the exponent $\vartheta_p=[(2+\beta-\gamma)(p-p_c)]^{-1}$ is sharp (see below) and the constants $\ka_1,\ka_2>0$,  depending only on $N, \gamma$ and $\beta$,  have an almost explicit expression. In the so-called \textit{good fast diffusion range}, $m_c<m<1$, solutions corresponding to $u_0 \in \LL^1_{\rm loc}( |x|^{-\gamma}\dx) $, turn out to be locally bounded. In the very fast diffusion range, $0<m\le m_c$, a counterexample given in Remark \ref{rem.upper} shows that this is not necessarily true. Indeed, the smoothing effect holds only for data in  $\LL^p_{\rm loc}( |x|^{-\gamma}\dx)$  with $p > p_c$, the so-called \textit{critical integrability exponent}, defined as\vspace{-1mm}
\begin{equation}\label{critical.integrability.exponent}
p_c= \frac{\left(1-m \right)\left(N-\gamma \right)}{2+\beta-\gamma}.\vspace{-1mm}
\end{equation}
Notice that $\vartheta_p>0$ whenever $p>p_c$ and that $p_c>1$ only when $m\in (0,m_c)$. We refer to the monograph \cite{JLVSmoothing} for a more detailed exposition of the relevance of such exponents in the non-weighted case $\beta=\gamma=0$ both for the smoothing estimates and for extinction phenomena.

The second main result is a precise quantitative lower bound for positive solutions, and it shows a remarkable property of WFDE, called ``instantaneous positivity'': as it happens in the case without weights, see \cite{BV-ADV-plap,BV-ADV,DiBenGianVesBook,DGV,DiBenedettourbVesp}, non zero data immediately produce strictly positive solutions.  A simplified version of our result reads:\vspace{-2mm}
\begin{equation}\label{lower.easy}
 \inf_{x \in B_{2R}(0)} u(t,x) \ge \kb \left[\frac{t}{R^{2+\beta-\gamma}}\right]^{\frac{1}{1-m}}\qquad\mbox{for any $t \in[0, t_*]$,}\vspace{-1mm}
\end{equation}
where $ t_\ast=t_*(u_0)\sim \|u_0\|_{\LL^1_\gamma(B_{R}(0))}^{1-m}$, is precisely defined in \eqref{choice_of_t}. We call $t_*$ the \textit{``minimal life time''} of the solution $u$, following \cite{BV-ADV}, since it represents the amplitude of the time interval in which any nonnegative local solution stays positive. Roughly speaking, if the solution is nonnegative  in  a small ball, then it becomes instantaneously positive  in  a bigger ball (expansion of positivity) and for some more time, precisely quantified by the minimal life time $t_*$; as a consequence, it becomes also H\"older continuous. The above lower bound is somehow optimal, indeed solutions to  \textit{fast diffusion}-type  equations may extinguish after a finite time $T=T(u_0)$: $t_*$ is a (sharp) lower bound for $T$; see Section \ref{sec:PART.II} for more details.

Our estimates are  quantitative and we show an (almost) explicit expression of $\kb>0$, which depends on the parameters $N, m, \gamma, \beta$, and possibly on $u_0$ or  other geometric quantities.  Note that in the good fast diffusion range $m\in (m_c,1)$, $\kb$ does not depend on the initial data. This does not happen in the linear case $m=1$, where the lower bound  depends on the initial data, and also, it is in contrast with the degenerate case $m>1$, where the finite speed of propagation forces to wait some time in order to have strict positivity, see \cite{JLVPorousMedioum}.  In the very fast diffusion range, $\kb>0$ also depends on $u_0$ through  $H_p\sim \|u_0\|_{\LL^p_\gamma(B_{R})}/\|u_0\|_{\LL^1_\gamma(B_{R})}$, see \eqref{def.HP} for a precise definition.

The proof of \eqref{lower.easy} is complex and contains the main new technical novelties of this paper. Due to the presence of the weights, the approach developed in \cite{BV-ADV} for the model equation ($\beta=\gamma=0$) that relies on moving plane methods (Alexandrov reflection principles) can not be applied.   Moreover, parabolic De Giorgi-type methods, typically used for equations with coefficients, see \cite{DiBenGianVesBook}, can be also applied to the case with weights in appropriate Muckenhoupt classes, see \cite{Surnachev-JDE,Surnachev-TMM}; however,  to our best knowledge, these techniques do not provide quantitative results, indeed the constants in the estimates are not always computable. We therefore develop in Part II a new strategy that allows us to keep the constants explicit.

Upper an lower bounds fairly combine in the form of parabolic Harnack inequalities,  our third main result proven in Part III.   In the non-weighted case already  it has been a longstanding problem to understand which form the Harnack inequality may take (if any) in the very fast diffusion range; the first answer has been given in  \cite{ BV-ADV-plap,BV-ADV} and then generalized to other contexts, see the monograph \cite{DiBenGianVesBook}.   A simplified version of our Harnack inequalities reads\vspace{-1mm}
\begin{equation}\label{harnack.easy}
\sup_{x \in B_R(0)} u(t, x) \leq \ka_3  \inf_{x \in B_R(0)}u(t, x) \qquad\mbox{for any}\qquad \frac{t_*}{2} < t < t_*.\vspace{-1mm},
\end{equation}
 where the constant $\ka_3 >0$ depends on $N, \gamma, \beta$ and possibly on the initial data $u_0$ in some ranges.
The minimal life time $t_*$  quantifies the size of the so-called intrinsic cylinders, which roughly speaking represent the right domain where regularity estimates can be stated in a clean form. The size of such intrinsic cylinders depends on the solution itself, due to the singular/degenerate character of the nonlinearity $u^m$, see \cite{DBbook, DiBenGianVesBook}.

In the very fast diffusion range, $\ka_3$ may depend on $u_0$ through some weighted $\LL^p$ norms, and this dependence cannot be avoided in view of explicit counter-examples, see \cite{BV-ADV, DGV} for the non-weighted case. On the other hand, in the good fast diffusion range, $\ka_3$ does not depend on $u_0$ anymore. In all cases we provide an explicit expression for $\ka_3$, see Remark \ref{rem.harnack.intro}.
It is remarkable that in \eqref{harnack.easy} we may take the infimum  and the supremum the same time (elliptic-type Harnack inequality), or even at a previous time (backward-type Harnack inequality): this feature is typical of fast diffusion or of singular evolutions \cite{BV-ADV-plap,BV-ADV,DiBenGianVesBook,DGV}, and is compatible with the fact that solutions can extinguish in finite time; this happens to be false (for general local weak solutions) even for the linear equation $m=1$, in which case forward Harnack inequalities typically hold.  See Theorem \ref{Harnack.2} below for a more general statement and remarks.

An important consequence of Harnack-type estimates is H\"older continuity, that we also establish in Part III. A simplified version of our estimates states that there exists $\alpha\in (0,1) $ and $\ka_\alpha >0$ such that, if $ 0\le u\le M_0$  on $ (t_0, t_0+ t_*]\times B_{4R_0}(0)$ we have, letting $\sigma=2+\beta-\gamma$\vspace{-1mm}
\[
|u(t,x)-u(\tau,y)|
\le\ka_\alpha R_0^{-\sigma} M_0 \bigl(|x-y|^{\sigma} + M_0^{m-1}|t-\tau| \bigr)^\frac{\alpha}{\sigma}\,, \vspace{-1mm}
\]
for all $t_0+\frac{5}{8}t_*\le t,\tau\le t_0+\frac{7}{8}t_*$  and  all $x, y \in  B_{R_0/4}(0)$. See Theorem \ref{thm.holder.nonlin} for a precise statement.  The H\"older exponent $\alpha$ depends on the constant $\ka_3$ of \eqref{harnack.easy}\,, and it will be chosen  uniformly  in the good fast diffusion range, where $\ka_3$ does not depend on $u$ (nor on $u_0$). On the other hand, in the very fast diffusion case, $\alpha$ may depend  on $u_0$ through some weighted $\LL^p$ norm: this is somehow natural, since  solutions  corresponding to data in $\LL^1_\gamma$ may be unbounded, as already discussed above.  We provide a (non optimal, but explicit) expression of the exponent $\alpha$ in Part III, and we show that it can vary depending on the cylinder: this may seem strange at a first sight, but indeed it is perfectly reasonable in view of the example given in Remark \ref{rem.holder}.
We can appreciate here the effect of the weights on the regularity of solutions.

The proof of H\"older continuity in the nonlinear case  depends on the regularity results for linear  equations  with weights. We prove in Part III both Harnack inequalities and H\"older continuity for weak solutions to linear equations  with measurable coefficients, whose prototype is given by $ u_t=|x|^\gamma\nabla\cdot\left(|x|^{-\beta}a(t,x)\nabla u\right)$, with $0<\lambda_0\le a(t,x)\le \lambda_1$; in our results we keep track  of  the dependence on $\lambda_0,\lambda_1$ in all constants, as Moser did in the non-weighted case, \cite{MoserCpam71}. We refer to Section \ref{sec.linear.regularity} for more details.

Finally, the Appendix contains the proof of the energy estimates of Part I, the proof of some functional inequalities that we use, together with a number of technical results. We have posponed those long and technical proofs there in order not to break the flow of the paper, and focus more on the main ideas.

We shall now present the main results and the different scenarios, together with the notation and definition of solutions that we are going to use.

\subsection{Precise statement of the main results in the different scenarios}\label{result.organization}
In order to state our main results, we need to introduce first some notations and definitions. We will write $a \asymp b$ whenever there exists constants $c_0, c_1 > 0$ such that $c_0 a \leq b \leq c_1 a$; we let $a \vee b = \mathrm{max}\{a,b\}$ and $a \wedge b = \mathrm{min}\{a,b\}$.

\noindent\textbf{Functional spaces. }Let $ p \geq 1 $ and $ \Omega \subseteq \RR^N $  be an open connected set with smooth boundary (at least $C^2$). For any $\gamma \in \RR$, $ \mu_\gamma $ will denote the measure $\mu_\gamma (\Omega) := \int_{\Omega}|x|^{-\gamma} dx$ and $\|f\|_{\LL^p_\gamma\left( \Omega \right)}:=\left( \int_{\Omega}|f|^p \, |x|^{-\gamma} \dx\right)^{\frac{1}{p}}$. We will denote by $\LL^p_{\gamma}\left( \Omega \right)$  the weighted $\LL^p$-space with respect to $\mu_\gamma$; it is known that $\LL^p_{\gamma}\left( \Omega \right)$ is a Banach space, see \cite{KO}. In what follows we will systematically deal with doubly weighted  Sobolev spaces, in which the norms of the function and of its derivatives are taken with respect to different measures. In the present weighted setting the usual definition of Sobolev spaces may not yield to a complete space, see \cite{KO}, therefore, we will follow the ideas of Fabes, Kenig and Serapioni \cite{FKS}, see also \cite{HKM}, and to avoid technical difficulties, we shall always assume that the parameters $\gamma$ and $\beta$ satisfy assumption \eqref{paramt.range}.
We define $H^1_{\gamma, \beta}\left( \Omega \right)$ to be the closure of $C^{\infty}\left( \bar{\Omega} \right)$ with the topology given by the norm $\|\phi\|^2_{H^1_{\gamma, \beta}\left( \Omega \right)}= \|\phi\|^2_{2, \gamma} + \|\nabla \phi\|^2_{2, \beta}$, and $\D_{\gamma, \beta}\left( \Omega \right)$  to be the closure of $C^{\infty}_{c}\left( \Omega \right)$ under the norm $\|\phi\|_{\D_{\gamma, \beta}\left( \Omega \right)}= \|\nabla \phi\|_{2, \beta}$. This procedure lead to the definition of a complete space in which functions have a unique weak gradient, obtained by approximation. Without any further assumption on the weights, the limit of such approximation may fall out of $\LL^1_{\rm loc}\left( \Omega \right)$, see \cite[section 2 and 3]{FKS} and  \cite{HKM}. As a consequence, solutions to \ref{WFDE} need to be considered in a suitable weak sense, as follows.

 \begin{defn}[Weak and strong solutions]\label{def.local.weak.sol}
Let $Q=(T_0, T] \times \Omega\subseteq (0,\infty)\times \RR^N$. A function $u: Q \rightarrow \RR $ is a \textsl{local weak solution }to equation (WFDE) in $Q$ if
\begin{equation}
u \in C_{\rm loc}((T_0, T); \LL^{2}_{\gamma, \rm loc}(\Omega)) \qquad \mbox{and} \qquad u^m \in \LL^2_{\rm loc}((T_0, T); H^{1}_{\gamma, \beta, \rm loc}(\Omega)),
\end{equation}
and the following identity holds true,
\begin{equation}\label{weak.formulation}\begin{split}
 \int_{\Omega}\left[u(t_2, x) \phi\left(t_2, x \right) \right. &-\left. u(t_1, x) \phi\left(t_1, x \right) \right] |x|^{-\gamma} \dx \\
 &=   \int_{t_1}^{t_2}\int_{\Omega}{u \phi_t \ |x|^{-\gamma} \dx } \dt - \int_{t_1}^{t_2}\int_{\Omega}{ \nabla{u^m}\cdot\nabla{\phi} \ |x|^{-\beta} \dx } \dt ,
\end{split}\end{equation}
for every open subset $[t_1, t_2]\times K \subset Q$ and for any test function $\phi$ such that
\begin{equation*}\label{test.space}
\phi \in W^{1, 2}_{\rm loc}((T_0, T); \LL^{2}_{\gamma}(K)) \cap \LL^2_{\rm loc}((T_0, T); \D_{\gamma, \beta}(K)).
\end{equation*}
A \textsl{local strong solution }to equation \ref{WFDE} is a local weak solution such that $u_t\in \LL^1_{\rm loc}((0,T)\,;\,\LL^1_{\gamma,\rm loc}(\Omega))$. \\
A local weak (or strong) sub (resp. super) solution satisfies  \eqref{weak.formulation} with $\le$ (resp. $\ge$) for any nonnegative test function in the same class.
\end{defn}

\noindent\textbf{About the class of solutions. }As already mentioned in the Introduction, most of our results will be proven for local strong solutions: lengthy (but nowadays standard) approximation procedures  allows one  to extend our results to a wider class of solutions, the so-called limit solutions \cite{JLVPorousMedioum},  sometimes also called SOLA, Solutions Obtained by Limit of Approximations \cite{BoGa, MiKuu}. In particular, our result apply to weak solutions in the sense of the above definition. However, such approximations are often long and technical in the framework of local solutions, but easier when dealing with global problems, like Cauchy, Dirichlet, Neumann, Robin, or even for large problems (Dirichlet problem whose solutions go to $+\infty$ on the lateral boundary). We shall say that most of the weaker concepts of solutions are included in the so-called class of limit solutions, i.e. limit of strong solutions, for which our estimates apply by a simple limiting process.

\medskip

\begin{figure}[h]
  \centering
  \includegraphics[width=19cm, height=8cm]{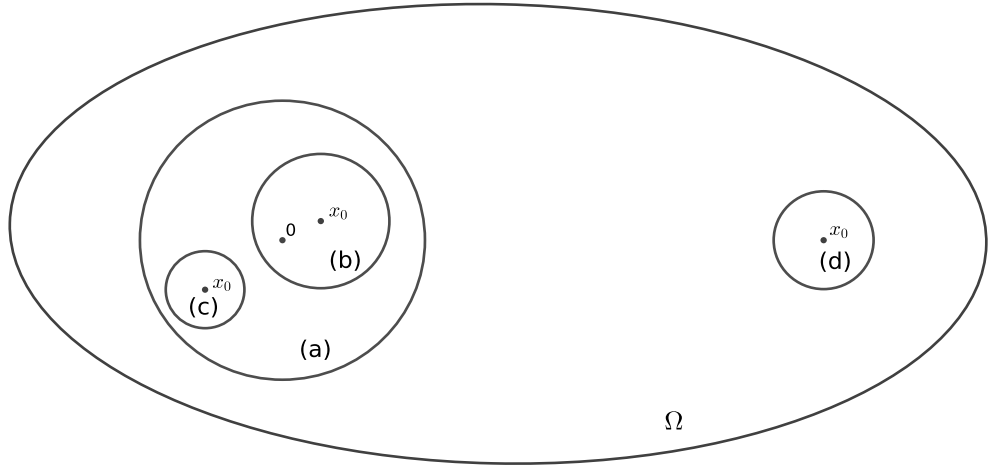}
  \caption{A representation of possible scenarios (a), (b), (c) and (d). }\label{fig.1}
\end{figure}

\noindent\textbf{Weights and different scenarios. }Our results concern quantitative a priori upper estimates of local type. We will consider a fixed cylinder of reference $Q=(T_0, T] \times \Omega\subseteq (0,\infty)\times \RR^N$, as in Definition \ref{def.local.weak.sol}, and the estimates will take place on a smaller cylinder, typically sufficiently far from the boundary of $\Omega$. Due to the lack of translation invariance of the weights, we need to find the right quantity that takes into account for the change of geometry, following \cite{CS-CPDE84,NPS-JDE2015,Surnachev-JDE,Surnachev-TMM}  we define for any $x_0 \in \RR^N$ and any $R>0$:
\begin{equation}\label{rho.pseudo}
\rho_{x_0}^{\gamma,\beta}(R):=\left( \int_{B_R(x_0)} |x|^{(\beta - \gamma) \frac{N}{2}} \dx \right)^{\frac{2}{N}}\,,\,\,\,\,\,\,\mbox{and}\,\,\,\,\,\,\ka_{16}^{-1}\,\rho_{x_0}^{\gamma,\beta}(R) \leq R^2  \frac{\mu_\gamma(B_R(x_0))}{\mu_\beta(B_R(x_0))} \leq  \ka_{16}\,\rho_{x_0}^{\gamma,\beta}(R)\,,
\end{equation}
where the latter equivalence is proven in Lemma \ref{technical.lemma.measures}.
Roughly speaking, the scaling properties of the equation change from $R^{2+\beta-\gamma}$ to (a multiple of) $R^2$, when we are considering problems around the origin or far from it, respectively.
There are at least four possible scenarios, see figure \eqref{fig.1}:
\begin{itemize}[leftmargin=.7cm]\itemsep2pt \parskip0pt \parsep0pt
\item[(a)] When $x_0=0$ and $ R_0 > 0 $ , we have $\rho_{x_0}^{\gamma,\beta}(R_0)\sim R_0^{2+\beta-\gamma}$.
\item[(b)] When $x_0\neq 0$ and $0 \in B_{R_0}(x_0)$, we have $\rho_{x_0}^{\gamma,\beta}(R_0)\sim R_0^{2+\beta-\gamma}$.
\item[(c)] When $x_0\neq 0$, $0 \not\in B_{R_0}(x_0)$ and $R_0 > |x_0|/2$, namely $x_0$ is relatively far from the origin but the singularity is still felt by the equation,  and  in this case we have $\rho_{x_0}^{\gamma,\beta}(R_0)\sim R_0^{2+\beta-\gamma}$.
\item[(d)]  When $x_0 \neq 0$,  $0 \not\in B_{R_0}(x_0)$ and $ 0<R_0\leq |x_0|/2$. This is the case  where $x_0$ is relatively far from the origin and does not heavily affect the geometry of the parabolic cylinders.  In this case we are essentially dealing with a nonlinear singular parabolic equation (governed by the nonlinearity $u^m$) and where the diffusion is driven by a uniformly elliptic operator; more specifically the standard (non-intrinsic) parabolic cylinders, depend on the ellipticity constants which in turn are proportional to $|x_0|^{\beta-\gamma}$, more precisely $\rho_{x_0}^{\gamma,\beta}(R_0)\sim R_0^2 \,|x_0|^{\beta-\gamma}$; note that all these latter quantities are bounded and bounded away from zero in this case;  see for instance \cite{MoserCpam71} for the linear case and \cite{BV-ADV,DiBenGianVesBook} for the nonlinear case.
\end{itemize}
We will focus only on the cases (a), (b) and (c) in which we have novel results, and where the geometry of the weights really plays a role; as already mentioned, the case (d) follows from nowadays standard results. For the sake of simplicity, from now on, we will always make one of the following assumptions on the ball $B_{R_0}(x_0)$ where our local estimates will take place:
\begin{itemize}[leftmargin=.7cm]\itemsep2pt \parskip0pt \parsep0pt
\item[(1)] Let $x_0=0$ and any $ R_0>0 $, or
\item[(2)] Let $x_0\neq 0$ and $|x_0|/32 \leq R_0 \leq |x_0|/16$, or
\item[(3)] Let $x_0 \neq 0$ and $(5/2)|x_0| \leq R_0 \leq 4 |x_0|$.
\end{itemize}
 We notice that under the assumption (1), (2) or (3) a simple calculation shows that (see proof of Lemma \ref{technical.lemma.measures})
\begin{equation}\label{r.sigma}
\ka_{17}^{-1} R^{2+\beta-\gamma}\leq R^2  \frac{\mu_\gamma(B_R(x_0))}{\mu_\beta(B_R(x_0))} \leq \ka_{17} R^{2+\beta-\gamma}\,,
\end{equation}
where the constant $\ka_{17}>0$ depends only on $N, \gamma, \beta$.  Before stating our main results we first need to introduce the following parameters whose role has been already explained above (recall that $p_c$ is defined in \eqref{critical.integrability.exponent})
\begin{equation}\label{definition.sigma.thethap}
\sigma= 2+\beta-\gamma \qquad\mbox{and}\qquad \vartheta_p=\frac{1}{\sigma p - \left( N-\gamma \right)\left( 1-m \right)}=\frac{1}{\sigma(p-p_c)}.
\end{equation}

\begin{thm}[Local Upper Bounds]\label{upper.easy}
Let $ u $ be a nonnegative local strong  (sub)solution  to \ref{WFDE} on the cylinder $\Omega \times (0, T]$. Let moreover $ p \geq 1 $ if $ m \in (m_c, 1) $ and $ p > p_c $ if $ m \in (0, m_c] $. Let $B_{2R_0}(x_0) \subset \Omega$ and assume that $B_{R_0}(x_0)$ satisfies either (1),(2) or (3). Then there exist $ \ka_1,\ka_2>0 $ such that  for any $ t \in (0, T]$ we have
 \begin{equation}\label{upper.bounds.rho}
\sup_{y\in B_{R_0}(x_0) }{u\left(t,y\right)} \leq \frac{\ka_1}{t^{(N-\gamma)\vartheta_p}} \left[\int_{B_{2R_0}(x_0)}{|u_0(y)|^{p} \ |y|^{-\gamma} \ \dy} \right]^{\sigma  \vartheta_p} + \ka_2 \left[\frac{\mu_\beta(B_{R_0}(x_0))}{\mu_\gamma(B_{R_0}(x_0))}\frac{t}{R_0^{2}} \right]^{\frac{1}{1-m}},
 \end{equation}
 where $\vartheta_p$ and $\sigma$ are defined as in \eqref{definition.sigma.thethap}. The constants $\ka_1, \ka_2$ depend   only on $N, \gamma$ and $\beta$.
\end{thm}
Part I contains the proof of Theorem \ref{local_upper_bounds}, which implies the above theorem as a particular case.
\begin{rem}\label{rem.upper} \rm \begin{enumerate}[leftmargin=*]\itemsep2pt \parskip0pt \parsep0pt
\item[(i)] The above smoothing effect may fail when $m<m_c$, if we choose exponents $p < p_c$. Indeed, there is an explicit counterexample to the above $\LL^p_\gamma \rightarrow \LL^{\infty}$ smoothing effect for solutions with initial data in $\LL^p_{\gamma, \rm loc}$ with $p < p_c$, given by the following function:
\begin{equation*}
  U(t,x)=c\,(T-t)^{\frac{1}{1-m}} |x|^{-\frac{\sigma}{1-m}},
\end{equation*}
where $c=c(m, N, \gamma, \beta)$ is chosen in such a way that $U$ becomes a local solution to \ref{WFDE} in the cylinder $(0,T) \times \RR^N$. In the non-weighted case $\beta = \gamma = 0$, the above counterexample was shown in \cite{BV-ADV, JLVSmoothing}.

\item[(ii)]It is worth noticing that when $\gamma=2+\beta$ or when $\sigma=0$, i.e. outside our range of parameters \eqref{paramt.range}, also in linear case the smoothing effect fails: this has been proved in \cite{CS-CPDE84} by means of counterexamples.

\item[(iii)]The above upper bound \eqref{upper.bounds.rho} contains two terms, the first takes into account the influence of the initial data, while the second takes into account the ``worst case scenario''\,, that happens when the local weak solution comes from the so-called large solutions, namely solutions to the Dirichlet problem which go to $+\infty$  at  the lateral boundary, see for instance \cite{BV-ADV-plap,BV-ADV}.
\end{enumerate}
\end{rem}

Our second result concerns quantitative positivity estimates. In order to state our main positivity results, we need to introduce first the following intrinsic quantities, for $p\ge 1$:
\begin{align}\label{def.HP}
    H_p\left( f, x_0, R \right)&:= \frac{\mu_\gamma(B_R(x_0))^{\sigma \vartheta_p}}{\mu_\gamma(B_R(0))^{\sigma\vartheta_p}}  \left[ \frac{\mu_\gamma(B_R(x_0))\left(\int_{B_R(x_0)}{f^p |x|^{-\gamma}  \dx}\right)^\frac{1}{p}}{\mu_\gamma(B_R(x_0))^{\frac{1}{p}}\int_{B_R(x_0)}{f |x|^{-\gamma} \dx}} \right]^{p \sigma \vartheta_p},\\
    \widetilde{H}_p&:=\widetilde{H}_p(f, x_0, R):= 1+\left(\frac{|x_0|}{R}\vee 1\right)^{\beta-\gamma} H_p(f, x_0, R)^{1-m} \ge 1\,,\nonumber
\end{align}
where $\vartheta_p$ and $\sigma$ are defined as in \eqref{definition.sigma.thethap}; we notice that in the cases (1), (2) and (3)  the quantities $|x_0|/R_0$ and $\mu_\gamma(B_R(x_0))^{\sigma \vartheta_p}\, \mu_\gamma(B_R(0))^{-\sigma\vartheta_p}$, become independent of $x_0, R$.
The above quantity is an adaptation to the weighted case of a similar one introduced in \cite{BV-ADV} and it plays an essential role in the positivity estimates: in particular, $\widetilde{H}_p$ encodes the geometric information of the weights which is relevant in the estimates.
An important aspect of these quantities, that will play an important role in our main results, is that both $H_p$ and $\widetilde{H}_p$ are scaling invariant, with respect to the natural scaling of the equation, see for instance formula \eqref{HP.scaling.inv}. Finally, we would like to emphasize that  in the good fast diffusion range, i.e. when we can choose $p=1$, $H_p$ (hence $\widetilde{H}_p$) does not depend on $f$:
\begin{equation*}
H_{1}(f,x_0,R)=\frac{\mu_\gamma(B_R(x_0))^{\sigma \vartheta_1}}{ \mu_\gamma(B_R(0))^{\sigma\vartheta_1}}.
\end{equation*}

We are now in the position to state our main positivity result.
\begin{thm}[Local Lower Bounds]\label{LOCAL.LOWER.BOUNDS}
Let $ u $ be a nonnegative local strong (super)solution to \ref{WFDE} on $(0,T)\times\Omega$ and let $0\le u_0\in \LL^p_{\gamma, \rm loc}(\Omega)$ with $ p \geq 1 $ if  $ m \in (m_c, 1) $ and $ p > p_c $ if $ m \in (0, m_c] $.  Let $B_{4R}(x_0)\subseteq\Omega$ and assume that $B_{R}(x_0)$ satisfies either (1), (2) or (3).  Define the minimal life time $t_*$ as
    \begin{equation} \label{choice_of_t}
     t_{\ast}= t_{\ast}(u_0, x_0,R)= \kappa_*\,R^\sigma \frac{ \|u_0\|_{\LL^1_{\gamma}(B_{R}(x_0))}^{1-m}}{\mu_\gamma(B_{R}(x_0))^{1-m}} \,.
    \end{equation}
Then, there exists $\kb=\kb\big(H_p(u_0,x_0,R), R, N,m,\gamma,\beta\big)>0$ such that
\begin{equation}\label{LOCAL.LOWER.INEQUALITY}
    \inf_{x \in B_{2R}(x_0)} u(t,x) \ge \kb \left[\frac{\mu_\beta(B_R(x_0))}{\mu_\gamma(B_R(x_0))}\,\frac{t}{R^{2}}\right]^{\frac{1}{1-m}}\qquad\mbox{for any $t \in[0, t_*\wedge T]$.}
\end{equation}
 The constant $\kappa_*>0$ depends on $N,m,\gamma,\beta$ and it is given in Corollary \ref{weak.positivity.MDP.L1}; $\kb$ has an (almost)  explicit expression given  in \eqref{const.lower.cor.L1} and depends on $H_p(u_0, x_0, R)$ only when $m\in (0,m_c]$.
\end{thm}
\begin{rem}\label{rem.low.bdds.intro}\rm\begin{enumerate}[leftmargin=*]\itemsep2pt \parskip0pt \parsep0pt
\item[(i)] Roughly speaking, the above lower bound \eqref{LOCAL.LOWER.INEQUALITY} shows that any bounded nonnegative solution
becomes instantaneously (strictly) positive on a whole time interval $(t_0, t_*(t_0)]$\,. This result will allow us to give \textit{an estimate on the size of the intrinsic cylinders}, which are the natural domain for Harnack and H\"older continuity estimates, see also Part III. We construct intrinsic cylinders inside arbitrary ones:  in the literature this is often an assumption, cf. \cite{DiBenGianVesBook}.

\item[(ii)] All the constants are computable: from the expression \eqref{const.lower.cor.L1} of $\kb$, we deduce that when $\widetilde{H}_p$, defined in \eqref{def.HP}, is large enough, then there is a constant $c_1>0$ depending on $N,m,p,\beta,\gamma$ such that
\begin{equation*}\label{Hp.tilde.behaviour}
\kb\asymp  \widetilde{H}_p^{-\frac{c_1\widetilde{H}_p^{1/2}}{m(1-m)}}.
\end{equation*}

\item[(iii)] As already mentioned above, when $m_c < m < 1$, the constant  $H_p$ does not depend on $u_0$ anymore, hence formula \eqref{LOCAL.LOWER.INEQUALITY} provides an absolute lower bound, i.e. independent of $u$ and $u_0$\,:
\[
    \inf_{x \in B_{2R}(x_0)} u(t,x) \ge \kb' \left[\frac{\mu_\beta(B_R(x_0))}{\mu_\gamma(B_R(x_0))}\,\frac{t}{R^{2}}\right]^{\frac{1}{1-m}}\qquad\mbox{for any $t \in[0, t_*\wedge T]$, }
\]
where $\kb'$ only depends on $N, m, \gamma, \beta$. However the presence of $u_0$ is still felt through $t_*\sim \|u_0\|_{\LL^1_\gamma}^{1-m}$.
\item[(iv)]   Part II contains a detailed proof of Theorem \ref{LOCAL.LOWER.BOUNDS}, which is the major technical contribution of this paper.
    Our proof applies also to the non-weighted case $\gamma=\beta=0$, and we recover the result of \cite{BV-ADV} with a different proof.
\end{enumerate}
\end{rem}
Our quantitative upper and lower bounds fairly combine into Harnack-type inequalities.
\begin{thm}[Harnack Inequalities]\label{Harnack.2}
Let $ u $ be a nonnegative local strong solution to \ref{WFDE} on $(0,T)\times\Omega$ and let $0\le u_0\in \LL^p_{\gamma, \rm loc}(\Omega)$ with $ p \geq 1 $ if $ m \in (m_c, 1) $ and $ p > p_c $ if $ m \in (0, m_c] $.  Let $B_{8R}(x_0)\subseteq\Omega$, $t_0\in [0,T)$ and assume that $B_{2R}(x_0)$ satisfies either (1), (2) or (3). Define
\[
t_*=t_*(u(t_0), x_0, 2R)=\kappa_* (2R)^{\sigma} \mu_\gamma(B_{2R}(x_0))^{m-1} \|u(t_0)\|_{\LL^1_{B_{2R}(x_0)}}^{1-m}.
\]
Then, for any $\varepsilon\in (0, 1)$, there exists $\ka_3>0$ such that for any $t,t\pm\theta\in [t_0+\varepsilon t_*(t_0), t_0+t_*(t_0)]\cap (0,T)$
\begin{equation}\label{harnack.easy.thm}
\sup_{x \in B_R(x_0)} u(t, x) \leq \ka_3   \inf_{x \in B_R(x_0)}u(t \pm \theta, x).
\end{equation}
The constants $\kappa_*,\ka_3>0$ always depend on $N,m,\gamma,\beta$ and are given in \eqref{t.star.MDP.L1} and  \eqref{const.harnack} respectively; $\ka_3$ may also depend on $R, x_0$ and $\varepsilon$, and, when $0<m\le m_c$, it depends on $ H_p(u_0, x_0, 2R)$ defined in \eqref{def.HP}.
\end{thm}
\begin{rem}\label{rem.harnack.intro}\rm\begin{enumerate}[leftmargin=*]\itemsep2pt \parskip1pt \parsep0pt
\item[(i)]  In \eqref{harnack.easy.thm} we may take the infimum  and the supremum at the same time (elliptic-type Harnack inequality, $\theta=0$); we can even take the infimum at a previous time (backward-type Harnack inequality, $\theta<0$): both inequalities are in contrast with the classical parabolic Harnack inequality valid for the linear heat equation ($m=1$), which  needs to be forward in time (infimum at a later time, $\theta>0$), \cite{Hadamard, Moser,  Pini}.  Indeed, elliptic and backward Harnack inequalities are typical features of fast diffusion equations, as already observed  in \cite{BV-ADV, DiBenGianVesBook}. They are false in general for the Heat Equation ($m=1$) and for the Porous Medium Equation ($m>1$), when dealing with local solutions (i.e. regardless of the boundary conditions), or in the case of solutions to the Cauchy problem posed on $\RR^N$. However for solutions to the homogeneous Dirichlet problem, also when $m\ge 1 $, elliptic and backward inequalities have been proven, see \cite{BFR,BFV-Parabolic,BV-PPR1,FGS86,SY}.
\item[(ii)]  Our result is quantitative, in the sense that all the constants are computable: from the expression \eqref{const.harnack} of $\ka_3$, when $\widetilde{H}_p$ is large enough, we deduce\vspace{-4mm}
\begin{equation*}\label{k3.behaviour.intro}
\ka_3
\asymp \varepsilon^{-\frac{\sigma p \vartheta_p}{1-m}} H_p(u_0, x_0, 2R) \, \widetilde{H}_p^{\frac{c\,\widetilde{H}_p^{1/2}}{m(1-m)}}\,, \vspace{-2mm}
\end{equation*}
where $\widetilde{H}_p$ is given in \eqref{def.HP}, and $c>0$ only depends on $N,m,p,\beta,\gamma$.
\item[(iii)]Notice that in the good fast diffusion range, we can choose $p=1$ and obtain a more classical form of the above Harnack inequality i.e. with the constant independent of $u_0$:
\[
\sup_{x \in B_R(x_0)} u(t, x) \leq \ka'_3 \inf_{x \in B_R(x_0)}u(t \pm \theta, x),
\]
with $t_*$ as in the formula above \eqref{harnack.easy.thm} and where $\ka'_3$ depend only on $N,m,p, \gamma, \beta $ and possibly $R, x_0, \varepsilon$. As already explained in the Introduction, in the very fast diffusion range we can not eliminate the dependence on $u_0$ in the above Harnack inequalities \eqref{harnack.easy.thm}; indeed, explicit counterexamples as in the non-weighted case can be constructed, see \cite{BV-ADV, DGV}.
\end{enumerate}
\end{rem}
Our last main result concerns continuity estimates, and the proof relies on some results for linear weighted equations that we describe later.

\begin{thm}[Interior H\"older Continuity]\label{thm.holder.nonlin}Let $u$ be a nonnegative local weak solution to the \eqref{WFDE} on $Q:=[0,T)\times\Omega$. Let $t_0\in [0,T)$ and $B_{16R_0}(x_0) \subset \Omega$ and assume that $B_{4R_0}(x_0)$ satisfies either (1), (2) or (3).  If $u \le M_0<\infty$ on $\left(t_0, T\wedge(t_0+ t_*) \right]\times B_{4R_0}(x_0)$ with $t_*=t_*(u(t_0), x_0,4R_0)$ as in \eqref{choice_of_t}, letting \[
D_0:=1\wedge \ka_{19}^{-2}\left(T \wedge t_*/8\right)^{1/\sigma}  \wedge \ka_{19}^{-2}\left(\rho^{\gamma,\beta}_{ x_0}\right)^{-1}\left(T \wedge t_*/8\right),
\]
with a suitably small $\ka_{19}>0$ as in \eqref{inverse.technical.inequality.measures}, then there exist $\alpha\in (0,1) $ and $\ka_\alpha'>0$\,, such that
\begin{equation}\label{Prop.HoCont.ineq}
|u(t,x)-u(\tau,y)|
\le\frac{\ka_\alpha'\, M_0}{D_0^\alpha \left(R_0 \wedge R_0^{\frac{\sigma}{2\vee\sigma}}\right)^\alpha} \left(|x-y| + M_0^{\frac{m-1}{2\vee\sigma}}|t-\tau|^{\frac{1}{2\vee\sigma}} \right)^\alpha\,,
\end{equation}
for all $t,\tau \in \left[t_0+\frac{5}{8}t_*, t_0+\frac{7}{8}t_*\right]\cap(0,T)$ and for all $x, y \in B_{R_0}(x_0)$.
The constants $\alpha\in (0,1) $ and $\ka_\alpha' >0$ depend on $N,\gamma,\beta, H_p(u(t_0), x_0, 4R_0)$ and possibly on $R_0,x_0$.
\end{thm}

\begin{rem}\label{rem.holder}\rm\begin{enumerate}[leftmargin=*]\itemsep2pt \parskip1pt \parsep0pt
\item[(i)]  We provide a (non optimal, but explicit) expression of the exponent $\alpha$ in Part III, and we show that it can vary depending on the base point $x_0$ and on the radius $R_0$\,: this may seem strange at a first sight, but indeed it is perfectly reasonable in view of the following example. Consider the Cauchy problem on the whole space: the fundamental (or Barenblatt) solution, has a selfsimilar form $B(t,x)=t^aF(|x|t^{-b})$ where $F(|x|)=A(D+ |x|^{\sigma})^{1/(m-1)}$, see \cite{BDMN2016b,BDMN2016a}. Clearly this explicit solution is merely H\"older continuous at zero when $\sigma\in (0,1]$, is $C^{1,\alpha}$ when $\sigma\in (1,2]$ and so on,  but such a solution is always $C^{\infty}(\RR^N\setminus\{0\})$\,.
We can appreciate here the effect of the weights on the regularity of solutions. Again, it is worth noticing that when $\sigma=0$, H\"older continuity fails, as well as the upper bounds, see Remark \ref{rem.upper}\,(ii).

\item[(ii)] We have decided to state the Theorem in this simplified form, to focus on the main result. Indeed it is quite easy to show that it holds for $t,\tau \in \left[t_0+\frac{5}{6}\varepsilon t_*, t_0+\varepsilon t_*\right]\cap(0,T)$, for any $\varepsilon\in (0,7/8)$, just the price of having a dependence on  in the constant $\ka_\alpha$, as it happens for the Harnack inequality of Theorem \ref{Harnack.2}.

\item[(iii)] The above theorem is stated in a general form emphasizing the fact that bounded solutions are $C^\alpha$ on a smaller intrinsic cylinder, whose size depends both on $t_*$ (i.e. on the $\LL^1_\gamma$ norm) and on the $\LL^\infty$ bound $M_0$. A closer inspection of the proof reveals that by slightly changing $t_*$ to $t_*= t_*(u(t_0), x_0, 8R_0)$, and using the upper bounds of Theorem \ref{upper.easy} or \ref{local_upper_bounds}, we can choose $M_0= c H_p(u(t_0), x_0, 8R_0) (t_*/R_0^\sigma)^{1/1-m}$,  by means of the same computation \eqref{harnack.estimate3bbb} as in the proof of Theorem \ref{Harnack.2}.

    Also, the exponent $\alpha$ depends on $H_p$ and $t_*$ in a quantitative way,
   \[
   \alpha\sim {\rm exp}\left(- \frac{c_6}{t_*}H_p^{\frac{c_7(1-m)}{m}H_p^{(1-m)/2} }\right)
   \]
   where $c_i>0$ only depend on $N,m,p,\beta,\gamma$. See the end of the proof of Theorem \ref{thm.holder.nonlin}, for more details.
   Therefore, the size of the intrinsic cylinders in the good fast diffusion range can be chosen to depend only on $t_*$ (i.e. on the  $\LL^1_{\gamma, \rm loc}$  norm of $u(t_0)$), while in the very fast diffusion range it has to depend also on $H_p$. The same happens for the exponent $\alpha$. This reveals a typical feature of the fast diffusion equation, for which there are strong differences between the two regimes. We last notice that such quantities are stable in the limit $m\to 1^-$, in which case we recover the linear results. On the other hand, by the above formulae, it is also clear that $\alpha\to 0^+$ when either $m\to 0^+$ or $H_p\to +\infty$.  This is compatible with the fact that solutions to the Dirichlet problem with $m=0$ (i.e. the logarithmic diffusion) extinguish immediately, cf. \cite{DPDS-LOG, JLVSmoothing}, and with the fact that if $H_p=+\infty$ the solution maybe unbounded.
\item[(iv)] We have already mentioned above that the optimal H\"older exponent $\alpha$ is not known, since $\alpha$ in general has to depend on $R_0,x_0$ and $u_0$. However, in some particular cases, $\alpha$ can be chosen uniformly in the whole range $m\in (0,1)$. This happens for instance in the Cauchy problem on the whole space, when we deal with the class of solutions trapped between two Barenblatt, in which case $H_p$ can be shown to be a suitable constant, see for instance \cite{BDMN2016b,BDMN2016a}. Notice that in the latter case, solutions are $C^\alpha$ around the origin, and classical elsewhere.
\item[(v)]\textit{About a uniform H\"older exponent. }A closer inspection of the proof reveals that indeed it is possible to choose a uniform H\"older exponent: in the good fast diffusion range, we can let $p=1$ and choose $M_0\asymp \|u(t_0)\|_{\LL^1_\gamma(B_{4R_0})}$ to obtain an (explicit) exponent $\alpha$ which only depends only on $N,\gamma,\beta$.
\end{enumerate}
\end{rem}

\noindent\textbf{Harnack inequalities and H\"older continuity in the linear case. }The study of quantitative regularity estimates for linear parabolic equations with measurable coefficients  began  with Moser \cite{MoserCpam71}. We  show in Section \ref{sec.linear.regularity} analogous quantitative Harnack inequalities and H\"older continuity for weak solutions to linear equations  with degenerate/unbounded coefficients, whose prototype is given by $ u_t=|x|^\gamma\nabla\cdot\left(|x|^{-\beta}a(t,x)\nabla u\right)$, with $0<\lambda_0\le a(t,x)\le \lambda_1$; in our estimates, we keep  track of the dependence on $\lambda_0,\lambda_1$. We do not claim originality for these results, indeed in many ranges of parameters they were already known \cite{CS-AA87,CS-RSMUP85,CS-AMPA84,CS-CPDE84, GW,GW2}; however we did not find in the literature the quantitative result that we needed, hence we sketch the proof in Section \ref{sec.linear.regularity}. The motivation for this analysis  comes from the application to nonlinear equations: our proof of the H\"older continuity for solutions to WFDE, heavily depends on the linear estimates.

\newpage

\subsection{Possible generalizations. }\label{ssec.111}
\textit{Other range of parameters. }Formally our results extend to a wider zone of parameters, namely zone (II) in Figure \ref{fig.3}, which amounts to require:  $\gamma>N$ and $\gamma-2>\beta\ge\frac{N-2}{N}\gamma$, note that in this case the weight $|x|^{-\gamma}$ is not integrable at $x=0$, also $|x|^{-\beta}$ is allowed to be not integrable and $\sigma=2+\beta-\gamma<0$. Allowing this range of parameters would require more technical results about the weighted functional spaces and inequalities involved in our proofs: we have decided to not treat this case here, since a rigorous proof would require a significant amount of technical results that would increase the length of the paper.
For the sake of simplicity we assume that $N\ge 3$, but our method works with straightforward modifications also when $N=1,2$.

\noindent\textit{More general equations. }We can allow for more general weights, equations and nonlinearities.
For instance, all the results of the present paper easily generalize to nonnegative solutions of
\[
w_\gamma^{-1}(x)u_t=\sum_{i,j=1}^N\partial_i \left( A_{i,j}(x)\partial_j u^m+B_i(x)u^m\right),
\]
with
\[
w_\gamma(x)\asymp |x|^{\gamma}\,,\qquad 0<\lambda_0|x|^{-\beta} |\xi|^2\le \sum_{i,j=1}^N A_{i,j}(x)\xi_i \xi_j\le \lambda_1 |x|^{-\beta}|\xi|^2,\qquad\mbox{and}\qquad |B_i(x)|\le \lambda_1 |x|^{-\frac{\beta+\gamma}{2}}
\]
for some constants $0<\lambda_0\le \lambda_1$; we can also translate the singularity at another point $x_0\neq 0$. Note that the upper bounds extend also to signed solutions, recalling that in this case we have to work with odd nonlinearities $u^m:=|u|^{m-1}u$.

A close inspection of the proofs reveals that all our results can be adapted, with some extra work, to nonlinearities $F(u)$ with $F\in C^1(\RR\setminus\{0\})$ with $F/F'\in {\rm Lip}(\RR)$ such that there exist $0<m_0 \le m_1$ such that
\[
\frac{1}{m_1}\le \left(\frac{F}{F'}\right)'\le \frac{1}{m_0}.
\]
It is often convenient to take $m_0,m_1\in (0,1)$, but  we can also allow $m_1\ge 1$; the above assumption guarantees that $t\mapsto t^{\frac{m_0}{m_1(1-m_0)}}u(t,\cdot)$ is monotone non-increasing, see for instance \cite{CP-JFA82}; as a consequence, all the proofs of the present paper can be repeated with minor modifications. The rough idea is that our results extend to nonlinearities that behave at zero like a concave power, $F(u)\asymp |u|^{m_0-1}u$ for $|u|\sim 0$, and for large $u$ behave like another (not necessarily concave) power $F(u)\asymp |u|^{m_1-1}u$ for $|u|\gg 1$. We stress that on one hand the qualitative results are still true (boundedness, positivity and continuity) also for more general nonlinearities. On the other hand, although qualitatively the same, the quantitative results shall have a quite different form, namely the exponents in the estimates and the dependence by the data in the constants and in the estimates may change in function of $m_0,m_1$; one advantage of the present method is that all the quantities can be controlled in a quantitative way.

The space-time estimates of Theorem \ref{spacetime.smoothing.theorem} and Proposition  \ref{lower.spacetime.smoothing} in its space-time form \eqref{lower.first.inequality.9},  can be extended to even more general nonlinear operators of the form
\[
w_\gamma^{-1}(x)u_t=\nabla\cdot A(t,x,u,\nabla u) + B(t,x,\nabla u, u),
\]
with $w_\gamma$ as above and $|A(t,x,u,\nabla u)|\le \lambda_1 |x|^{-\beta} \big|\nabla |u|^m\big|$ and $A(t,x,u,\nabla u)\cdot \nabla u\ge \lambda_0 |x|^{-\beta} \big|\nabla |u|^m\big|^2 $. Again, the assumptions on the power-like nonlinearity can be weakened, as above, and we can allow a concave $F$ with $F(u)\asymp |u|^{m-1}u$ for $u\sim 0$, and regular outside zero. As for the lower order term $B$, the typical assumption would be $|B(t,x,u,\nabla u)|\le \lambda_1|x|^{-\beta}\big|\nabla |u|^m\big|+\lambda_1^2w_\gamma^{-1}(x) |u|^m$\,, but they can be weakened.
On one hand, it is possible to obtain upper bounds in a refined form like Theorem \ref{upper.easy} or \ref{local_upper_bounds} also in this generality.
On the other hand, precise lower bounds like in Theorem \ref{LOCAL.LOWER.BOUNDS} or \ref{positivity.MDP.thm.general},  are not easily extended in this degree of generality: the major technical difficulty is represented by the absence of time monotonicity for solution to a homogeneous Dirichlet problem, namely that $t\mapsto t^{\frac{1}{1-m}}u(t,\cdot)$ is monotone non-increasing.

Finally, our methods can be adapted to hold also on Riemannian manifolds; we can possibly  allow for  manifolds with unbounded curvature, as already mentioned, and this partially motivates the present paper.

\subsection{Weighted Functional Inequalities}\label{sect.funct.ineq}
\begin{figure}[!]
  \centering
  \includegraphics[width=16.50cm, height=7cm]{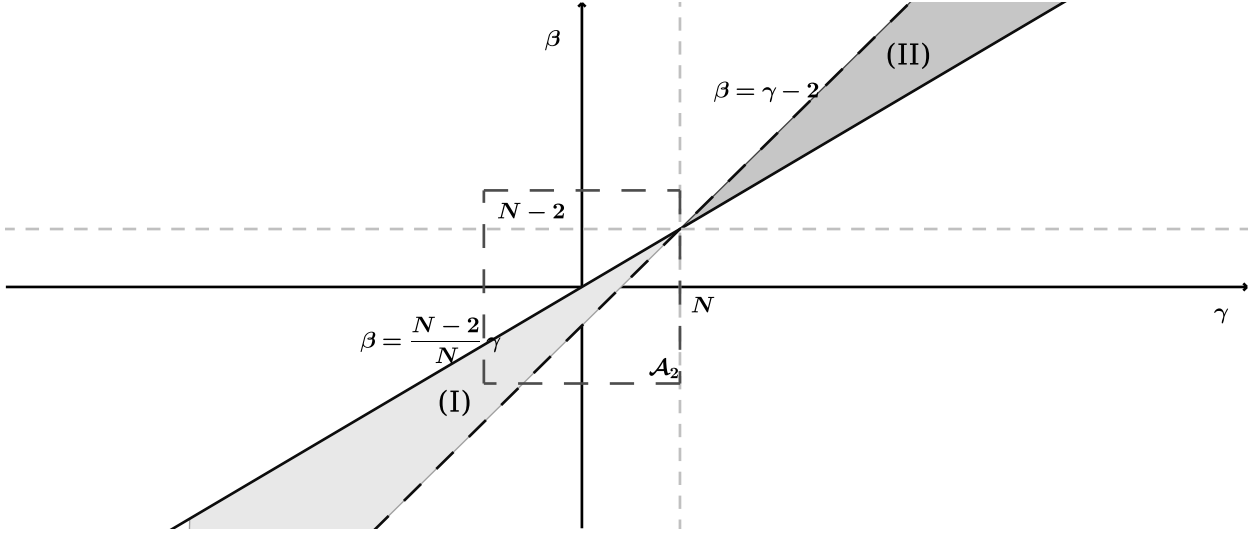}
  \caption{In light grey the region (I) of  parameters $\beta$ and $\gamma$ as in  \eqref{paramt.range}, where all of our results hold.  Note that a big part of the region falls outside the $\mathcal{A}_2$ region. The region (II) in dark grey correspond to a range of parameters that we do not treat, but where our results formally apply, as discussed in Subsection \ref{ssec.111}}\label{fig.3}
\end{figure}

In order to study regularity properties of the solution to \ref{WFDE} a key point in our approach is represented by weighted functional inequalities, that we briefly recall here. For any $\gamma < N$, consider the measure $\mu_\gamma(B):=\int_B|x|^{-\gamma}\dx$, which is known to be doubling, i.e.
\begin{equation}\label{Doubling.Constant}
\mu_\gamma(B_{2R}) \leq D_\gamma \mu_\gamma(B_R),
\end{equation}
where $B_R$ is a ball contained in $\RR^N$ and the constant $D_\gamma$ depends only on the dimension $N$ and the parameter $\gamma$,
see \cite[Chp. 15]{HKM}. On the whole $\RR^N$, there is a celebrated family of interpolation inequalities, the so-called Caffarelli-Kohn-Nirenberg inequalities \cite{Caffarelli-Kohn-Nirenberg-84}, that we state hereafter in a special case, namely as in  \eqref{w_sobolev0}. Let $ \gamma, \beta $ as in \eqref{paramt.range}, then there exists a constant $\overline{S}_{\gamma, \beta}>0$ such that for any $ f\in C^{\infty}_c \left(\RR^N\right) $ the following inequality  holds
\begin{equation} \label{w_sobolev}\tag{CKNI}
 \|f\|_{\LL^{\sr}_{\gamma}\left(\RR^N \right)}\le\overline{S}_{\gamma, \beta} \|\nabla f\|_{\LL^2_{\beta}\left( \RR^N \right)}\qquad\mbox{where}\qquad\sr= 2\frac{N-\gamma}{N-\left(2+\beta\right)}\,,
\end{equation}
where the weighted $\LL^p$ norms are defined in subsection \ref{result.organization}.
This family of inequalities contains both the classical Sobolev inequality ($ \gamma=\beta = 0 $) and the Hardy inequality ($ \beta=\gamma - 2 $), cf. \cite{KOBook}. In our range of parameters \eqref{paramt.range} (see Figure \ref{fig.3}) we always have $\sr\in[2, 2N/(N-2)]$.
\begin{prop}\label{prop.Sobolev.Balls}
Let $N\geq 3 $, $ \gamma, \beta $ be as in \eqref{paramt.range}, $r^*$ be as in \eqref{w_sobolev} and $\sigma$ as in \eqref{definition.sigma.thethap}. Let $x_0 \in \RR^N $ and $R>0$. Then, there exists a constant $S_{\gamma,\beta}>0$ depending only on $N,\beta, \gamma$, such that for any $ f\in H^1_{\gamma,\beta}(B_R(x_0)) $ the following inequality holds true
\begin{equation}\label{WSIcomplete}\tag{CKNI2}
    \|f\|_{\LL^{\sr}_{\gamma}\left( B_R(x_0) \right)} \leq S_{\gamma,\beta} \left( \|\nabla f\|_{\LL^2_{\beta}\left( B_R(x_0) \right)} +  \mu_\gamma(B_R(x_0))^{\frac{-\sigma}{2(N-\gamma)}}  \| f\|_{\LL^2_{\gamma}\left( B_R(x_0)\right)} \right)\,.
\end{equation}
\end{prop}
The above \eqref{WSIcomplete} will play the role of the classical Sobolev inequality in the proof of both upper and lower bounds, hence we have have given a short proof of the above inequality in \nameref{sec:appendix.C}.

Another essential tool in our proofs will be the following weighted Poincar\'{e} inequality.
\newpage
\begin{prop}[Poincar\'e Inequality]\label{zero.mean.Poincare}
   Let be $N \geq 3$ and $\gamma, \beta < N $ as in \eqref{paramt.range}, $x_0 \in \RR^N$ and $R>0$. Then there exists a constant $P_{\gamma, \beta}>0$ such that for any $\phi \in H^1_{\gamma, \beta}(B_{R}(x_0))$ we have
    \begin{equation}\label{Weighted.Poincare.ineq}
            \left( \frac{1}{\mu_\gamma(B_R(x_0))} \int_{B_R(x_0)} |\phi-\overline{\phi_{\gamma}}|^2 |x|^{-\gamma}\dx  \right)^{\frac{1}{2}} \leq P_{\gamma, \beta} \, \, R \left(\frac{1}{\mu_\beta(B_R(x_0))} \int_{B_R(x_0)} |\nabla \phi|^2 |x|^{-\beta}\dx  \right)^{\frac{1}{2}},
    \end{equation}
    where $ \overline{\phi_{\gamma}} = \mu_{\gamma}\left( B_R(x_0)\right)^{-1} \int_{B_R(x_0)}{\phi |x|^{-\gamma}  \dx } $; the constant $P_{\gamma, \beta}>0$ depends only on $N, \gamma, \beta$.
\end{prop}
The above Poincar\'e inequality is a direct consequence of a more general one, called \textit{Sobolev-Poincar\'{e} inequality} see \cite[Chp. 15]{HKM}, and also \cite{FKS,FGW,HK,Giusti} and  references  therein for the known results. For the sake of completeness, we have decided to give a proof of the above inequality in
\nameref{sec:appendix.C}.

\medskip

In the proof of positivity estimates we will also use BMO  -Bounded Mean Oscillation-  weighted spaces and a weighted John-Nirenberg inequality; we recall here the definition and  inequalities that we will use, for convenience of the reader.
\begin{defn}
  A function $f \in \LL^1_{\gamma, \rm loc}(\Omega)$ is said to be in $\textit{BMO}_{\gamma}(\Omega)$ if
  \begin{equation*}\label{BMO}
      \|f\|_{BMO_{\gamma}(\Omega)}:=  \sup_{B\Subset \Omega}\frac{1}{\mu_\gamma(B)}\int_{B}|f-\overline{f}_B| |x|^{-\gamma} \dx < +\infty,
  \end{equation*}
  where $B$ are balls compactly contained in $\Omega$ and $\overline{f}_{B}= \mu_\gamma(B)^{-1}\int_{B}f |x|^{-\gamma} \dx$.
\end{defn}
The following version of the John-Nirenberg Lemma can be found for instance in  \cite[Thm 18.3]{HKM}.
\begin{lem}[Weighted John-Nirenberg inequality, \cite{HKM}]\label{JOHN.NIRENBERG}
  Let $f \in BMO_{\gamma}(\Omega)$. Then, for any ball $B$ compactly contained in $\Omega$ the following inequality holds
  \begin{equation}\label{JonhNirenberg}
        \frac{1}{\mu_\gamma(B)}\int_{B}e^{\mathsf{s}|f-\overline{f}_{B}|}|x|^{-\gamma} \dx \leq \ka_{5} \qquad\, \text{for any} \qquad\ 0 < \mathsf{s} \leq \frac{1}{\ka_{6}  \|f\|_{BMO_{\gamma}(\Omega)}}\,,
  \end{equation}
  where $\ka_{5}$ and ${\ka_{6}}$ are positive constants depending only on $N, \gamma$.
\end{lem}
From the above inequality  there  follows Corollary \ref{Reverse.Holder.Corollay} that will play a crucial role in Part II.
\begin{cor}\label{Reverse.Holder.Corollay} Let $u: \Omega \rightarrow \RR$ be a positive measurable function such that $\log (u) \in BMO_{\gamma}(\Omega)$. Then  the following inequality holds
\begin{equation}\label{rev.hold.cor.intro}
\|u\|_{\LL^{\mathsf{s}}_{\gamma}(B)} \leq \ka_{7}^{\frac{2}{\mathsf{s}}} \mu_\gamma(B)^{\frac{2}{\mathsf{s}}} \|u\|_{\LL^{-\mathsf{s}}_{\gamma}(B)} \qquad \mbox{for any} \qquad 0 < \mathsf{s} < \frac{1}{\ka_{6} \|\log(u)\|_{BMO_{\gamma}(\Omega)}},
\end{equation}
where $B$ is any ball compactly contained in $\Omega$, $\ka_6>0$ is as in Lemma \ref{JOHN.NIRENBERG}  and $\ka_7>0$ is a constant depending  only  on $N, \gamma$\,.
\end{cor}
\noindent\textbf{Proof.} The proof of this result follows by the weighted John-Nirenberg inequality \eqref{JonhNirenberg}; it is a straightforward adaptation of the proof of \cite[Theorem 4]{Moser}, see also \cite[Proposition 4.4]{BGV}.\qed

\section{Part I. Local upper bounds and energy estimates}

The main result of this part is a precise and quantitative local upper bound, which ensures that sufficiently locally integrable solutions   are indeed bounded at a later time, as precisely stated   below; we state here the most general form of upper bounds, which includes the cases considered in Theorem \ref{upper.easy}.
 \begin{thm}\label{local_upper_bounds}
Let $ u $ be a nonnegative  local strong solution to \ref{WFDE} on the cylinder $\Omega \times (0, T]$. Let moreover $ p \geq 1 $ if $ m \in (m_c, 1) $ and $ p > p_c $ if $ m \in (0,  m_c] $. Let $B_{2R_0}(x_0) \subset \Omega$, and assume that $B_{R_0}(x_0)$ satisfies either (a), (b) or (c).  Then there exist $\ka_8, \ka_9 > 0$  such that  for any  $0< t-t_0 < T$
\begin{equation}\label{upper.provisional.2}
 \sup_{y\in B_{R_0}(x_0) }{u\left(t,y\right)} \leq \frac{\ka_8}{(t-t_0)^{(N-\gamma)\vartheta_p}} \left[\int_{B_{2R_0}(x_0)}{|u(t_0,y)|^{p} \ |y|^{-\gamma} \ \dy} \right]^{\sigma  \vartheta_p} + \ka_9 \left[\frac{t-t_0}{R_0^{\sigma}} \right]^{\frac{1}{1-m}},
\end{equation}
where $\vartheta_p$ and $\sigma$ are as in \eqref{definition.sigma.thethap}. The constants $\ka_8, \ka_9$  depend only on $N, \gamma, \beta$ and on the quotient $|x_0|/R_0$;   they both have an  explicit expression given in \eqref{constants.upper.bounds}.
 \end{thm}

We have already observed the main features of the above upper bounds, see Remark \ref{rem.upper}.  Note that Theorem \ref{upper.easy} is a particular case of the above, in the sense that the dependence on $|x_0|/R_0$ can be eliminated in the constants $\ka_8$ and $\ka_9$, simply by choosing the cases (1), (2) and (3).  We shall now proceed with the proof of the Theorem,  in  Subsection \ref{Proof.Thm2.1}. Before that, we need a number of preliminary results, some of them having their own interest.

\subsection{Local space-time energy estimates}
We collect in this Subsection all the energy inequalities that we will use in the rest of the paper, the proof is quite  technical;  but follows by nowadays standard ideas; the hidden difficulty  lies on  the careful approximations needed to deal with the singular/degenerate nature of the weights. We postpone the proof to \nameref{sec:appendix.A}, in order to not  to break  the flow of the proof.
In order to state the energy inequalities in all the possible scenarios, we introduce an auxiliary function: to avoid unnecessary complications, we consider balls $B_{R_1}(x_0) \subset B_{R_0}(x_0)$ such that $0 \not \in \overline{B_{R_0}(x_0)\setminus B_{R_1}(x_0)}$. Let $0<R_1<R_0$ and $ \sigma= 2+\beta-\gamma \in (0, \infty)$ and define
\begin{equation}\label{function_h}
h_{\sigma}\left(R_0, R_1, x_0\right):=
\left\{\begin{array}{lll}
 \left(\dfrac{R_0+|x_0|}{R_0-R_1} \right)^{2-\sigma}\,, & \mbox{if }0<\sigma<2\,,\\
 1\vee \left(\dfrac{R_0-R_1}{R_1-|x_0|} \right)^{\sigma -2} \,, & \mbox{if }\sigma\ge 2\mbox{ and }0 \in B_{R_1}(x_0)\,,\\
 1\vee \left(\dfrac{R_0-R_1}{|x_0|-R_0} \right)^{\sigma -2} \,, & \mbox{if }\sigma\ge 2\mbox{ and }0 \not \in \overline{B_{R_0}(x_0)}\,.\\
 \end{array}\right.
\end{equation}
The function $ h_{\sigma}$ takes into account the change of geometry and cover more general cases than  the ones defined in (1), (2) and (3). We moreover observe that whenever $B_R$ satisfies one of the hypothesis  (1),  (2) or (3)   we have
\begin{equation*}
 h_{\sigma}\left( 4R, R, x_0 \right) \asymp  h_{\sigma}\left(2R, R, x_0\right) \asymp C_{\gamma, \beta} \qquad \mbox{and} \qquad  \ka_{17}^{-1}R^{\sigma} \leq R^2 \frac{\mu_\gamma(B_R(x_0))}{\mu_\beta(B_R(x_0))}\leq \ka_{17} R^{\sigma}\,,
\end{equation*}
where the constant $C_{\gamma, \beta}>0$ depends only on the constants $\gamma, \beta$.
This is consistent with the weighted estimates proved in the linear case, see \cite{CS-RSMUP85, GW, GW2} and also Remark \ref{rem.GW} below.
\begin{lem}[Energy estimates] \label{theorem_energy_inequality}
Let $x_0 \in \RR^N$,  $0 < R_1 < R$ such that $0 \not \in \overline{B_R(x_0)\setminus B_{R_1}(x_0)}$ and let  $0<m < 1 $, $0\le T_0 < T_1 < T$ .
 \begin{itemize}[leftmargin=*]
 \item Let $u$ be a non-negative local strong subsolution to \ref{WFDE} in $(T_0,T)\times B_{R}(x_0)$. Let  $p > 1$  and assume $u\in \LL^p_{\gamma}((T_0,T)\times B_{R}(x_0))$. Then there exists $c_1>0$ depending only on $m,p,N $, such that

\begin{align}\label{sup.energy.inequality.upper}
   \sup_{\tau \in \left[T_1, T \right]}&{\left\{\int_{B_{R_1}\left( x_0 \right)}{u^p(\tau, x) \ |x|^{-\gamma} \dx}\right\}} + \int_{T_1}^T\int_{B_{R_1}(x_0)}{\left|\nabla u^{\frac{p+m-1}{2}} \right|^2 \ |x|^{-\beta} \dx  \dt} \nonumber \\
  & \leq c_1 \left[\frac{h_{\sigma}\left( R, R_1, x_0 \right)}{\left(R - R_1 \right)^{\sigma}}  + \frac{1}{T_1-T_0} \right] \int_{T_0}^T\int_{B_R(x_0)}{\left( u^{p+m-1} + u^{p} \right) \ |x|^{-\gamma} \dx \dt}.
\end{align}

\item Let $\delta>0$ and $u\ge \delta$ be a local strong supersolution to \ref{WFDE} in $(T_0,T)\times B_{R}(x_0)$.\\
- For all $0 < p < 1-m$ there exists $c_2>0$ depending only on $m,p,N $, such that
 \begin{align}\label{energy.inequality.last}
    &  \int_{B_{R_1}(x_0)} u(T_0,x)^{p}  |x|^{-\gamma} \dx+ \int_{T_0}^{T_1}\int_{B_{R_1}(x_0)} \left|\nabla u^{\frac{p+m-1}{2}} \right|^2  |x|^{-\beta} \dx  \dt   \nonumber \\
& \leq c_2 \left[\frac{h_\sigma(R, R_1, x_0)}{(R-R_1)^\sigma}+ \frac{1}{T-T_1}\right] \int_{T_0}^T\int_{B_R(x_0)}(u^{p+m-1}+ u^{p})  |x|^{-\gamma} \dx \dt;
 \end{align}

- For all $p>0$ there exists $c_3>0$ depending only on $m,p,N $, such that
     \begin{align}\label{sup.energy.inequality.lower}
 \sup_{\tau \in \left[T_1, T \right]}&\left\{\int_{B_{R_1}(x_0)}{u(\tau, x)^{-p} \ |x|^{-\gamma} \dx}\right\} + \int_{T_1}^T\int_{B_{R_1}(x_0)}{\left|\nabla u^{\frac{-p+m-1}{2}} \right|^2 \ |x|^{-\beta} \dx  \dt} \nonumber \\
  & \leq c_3  \left[\frac{h_\sigma(R, R_1, x_0)}{\left(R-R_1 \right)^\sigma} + \frac{1}{T_1-T_0} \right]
  \int_{T_0}^T\int_{B_R(x_0)}{\left( u^{-p+m-1} + u^{-p} \right) \ |x|^{-\gamma} \dx \dt}.
 \end{align}
\end{itemize}
The constants $c_1,c_2,c_3$ have an explicit expression given in the proofs.
\end{lem}
The next Lemma plays the role of the celebrated Caccioppoli inequality and corresponds to the above estimates in the borderline case $p=1-m$. It will be a key ingredient in the proof of the positivity estimates of \ref{sec:PART.II}.
\begin{lem}[Caccioppoli estimates]\label{Lem.Caccioppoli}
Let $m \in (0, 1)$, $\delta>0$ and $u\ge \delta$ be a local strong supersolution to \ref{WFDE} in $(T_0,T)\times B_{R}(x_0)$. Let $ \psi \in C_c^{\infty}(\Omega) $ with $\supp(\psi)\subseteq B_{R}(x_0)\subset \Omega$ and $T_0 \leq \tau < t \leq T$. We have
\begin{equation}\label{caccioppoli.inequality}\begin{split}
   & \int_{B_R(x_0)}{u(\tau, x)^{1-m}  \psi^2 \ |x|^{-\gamma} \dx} + \frac{m^2}{2}\left(1-m \right) \int_{\tau}^{t}{\int_{B_R(x_0)}{ \psi^2 \ |\nabla \log{u}|^2 \ |x|^{-\beta} \dx \dt}}   \\
& \leq  2\left(1-m \right)\int_{\tau}^{t}{\int_{B_R(x_0)}{  |\nabla \psi|^2 \ |x|^{-\beta} \dx \dt}} + \int_{B_R(x_0)}{u(t, x )^{1-m} \ \psi^2  \ |x|^{-\gamma} \dx}\,.
\end{split}\end{equation}
\end{lem}
We finally observe that letting $m\to 1^-$ in \eqref{caccioppoli.inequality}, and recalling that $u^{1-m}/(1-m)\to \log u$ in such limit, we recover the classical Caccioppoli estimate, valid for $m=1$ and $\beta=\gamma=0$, cf. \cite{AS-ARMA67,Moser,MoserCpam71,T-CPAM68}.

\subsection{Behaviour of local $\LL^p_\gamma$ norms.}
\begin{prop} \label{Herrero_Pierre_Thm}
Let $m \in (0, 1)$ and $u\left( t,x \right): \left[T_0, T_1 \right] \times \Omega   \rightarrow  \RR  $ be a nonnegative local strong solution to the \ref{WFDE}. Let $x_0 \in \Omega $ and $R > 0$  be  such that $B_{2R}(x_0)\subset \Omega$. Then, for any $ 0 \le t, \tau \in \left[T_0, T_1\right]$ we have
\begin{equation}\label{Herrero_Pierre_inequality}
\left[\int_{B_{R}(x_0)}{u\left(t, x \right) |x|^{-\gamma} \dx}\right]^{1-m} \leq \left[\int_{B_{2R}(x_0)}{u\left(\tau,  x \right)  |x|^{-\gamma} \dx}\right]^{1-m} + \ka_{10} \frac{\mu_\gamma(B_R(x_0))^{1-m}}{\rho^{\gamma, \beta}_{x_0}(R)} |t-\tau|,
\end{equation}
where the constant $\ka_{10}$ depends only on $ N, m, \gamma$ and $\beta$.  Moreover, under the assumption  that $B_{R_0}(x_0)$ satisfies either \textrm{(1)}, (2) or (3), \eqref{Herrero_Pierre_inequality} becomes
\begin{equation}\label{Herrero.Pierre.123}
\left[\int_{B_{R}(x_0)}{u\left(t, x \right) |x|^{-\gamma} \dx}\right]^{1-m} \leq \left[\int_{B_{2R}(x_0)}{u\left(\tau,  x \right)  |x|^{-\gamma} \dx}\right]^{1-m} + \ka_{10}' \frac{\mu_\gamma(B_R(x_0))^{1-m}}{R^\sigma} |t-\tau|,
\end{equation}
where $\ka_{10}'=\ka_{10}\,\ka_{17}\,\ka_{16}^{-1}$; $\ka_{16}\,,\ka_{17}$ depend only on $N, \gamma, \beta$ and are as in \eqref{rho.pseudo}, \eqref{r.sigma} respectively.
\end{prop}
\noindent\textbf{Remark. }The above Lemma quantifies the displacement of local mass backward and forward in time, and has been first proved in \cite{HerreroPierre} in the non-weighted case and for solutions to the Cauchy problem on $\RR^N$: inequality \eqref{Herrero_Pierre_inequality} implies conservation of mass (letting $R\to\infty$, when $m>m_c$). It has also been used to prove estimates from below for the extinction time in different contexts: for the Cauchy problem when $m>m_c$, for any $m\in (0,1)$ for the Dirichlet problem or on Manifolds with negative curvature, see \cite{BGV-JEE, BV-ADV, JLVSmoothing}. Also in this weighted case it allows  one  to prove the same results: we use it here (also) to prove lower bounds for the extinction time for the Minimal Dirichlet problem for all $m\in (0,1)$, see \eqref{Bounds.FET}.

\noindent\textbf{Proof.} Let $ \phi $ be a cut-off function supported in $ B_{2R}(x_0) $ and let $\phi = 1 $ in $ B_{R}(x_0) $. In what follow we will write $B_{R}$ instead of $B_{R}(x_0)$ when no confusion arises. We adopt the notation $\mathcal{L}_{\gamma, \beta}f= |x|^\gamma\nabla\cdot\left(|x|^{-\beta}\nabla f\right)$, see also  \eqref{operator_on_cutoff_functions}. Let us compute
\begin{equation}
\begin{aligned}\label{derivation_local_norm}
& \left| \frac{d}{dt} \int_{B_{2 R}}{u( t,x ) \phi\left( x \right) |x|^{-\gamma} \dx} \right|  =\left| \int_{B_{2 R}}{\mathcal{L}_{\gamma, \beta}\left( u^m \right) \phi\left( x \right) |x|^{-\gamma} \dx} \right| \\
& = \left| \int_{B_{2 R}}{u^m \mathcal{L}_{\gamma, \beta}\left( \phi\left( x \right) \right)  |x|^{-\gamma} \dx} \right|
 \leq  \int_{B_{2R}}{u^m \left|\mathcal{L}_{\gamma, \beta}\left(\phi\left( x \right)\right)\right|  |x|^{-\gamma} \dx}.
\end{aligned}
\end{equation}
The  H\"older's inequality with conjugate exponents $\frac{1}{m}$ and $\frac{1}{1-m}$  gives
\begin{equation*}
\begin{aligned}
 \int_{B_{2R}}{u^m \left|\mathcal{L}_{\gamma, \beta}\left(\phi \right)\right|  |x|^{-\gamma} \dx} & \leq \left[ \int_{B_{2 R}}{u \phi \left( x\right)  |x|^{-\gamma} \dx} \right]^{m} \left[ \int_{B_{2R}}{\phi ^{\frac{-m}{1-m}} \left|\mathcal{L}_{\gamma, \beta}\left(\phi \right)\right|^{\frac{1}{1-m}}  |x|^{-\gamma} \dx} \right]^{1-m} \\
 & := C\left( \phi \right) \left[ \int_{B_{2R}}{u \phi \left( x\right)  |x|^{-\gamma} \dx} \right]^{m}
\end{aligned}
\end{equation*}
Notice that joining the above estimate and (\ref{derivation_local_norm}) we get the closed differential inequality
\begin{equation*}
  \left| \frac{d}{dt} \int_{B_{2 R}}{u\left( t,x \right) \phi\left( x \right) |x|^{-\gamma} \dx} \right| \leq C(\phi) \left[ \int_{B_{2 R}}{u\left(t,x \right) \phi\left( x \right) |x|^{-\gamma} \dx} \right]^m.
\end{equation*}
An integration in time shows that for all $t,\tau\ge 0$ we have
\begin{equation*}
\left( \int_{B_{2 R}}{u\left( t,x \right) \phi\left( x \right) |x|^{-\gamma} \dx} \right)^{1-m} \leq  \left( \int_{B_{2 R}}{u\left( \tau,x \right) \phi\left( x \right) |x|^{-\gamma} \dx} \right)^{1-m}  + \left( 1-m \right)C\left(  \phi \right) \left| t-\tau \right|.
\end{equation*}
Since $ \phi $ is supported in $ B_{2R} $ and equal to $ 1$ in $  B_{R}$, this implies \eqref{Herrero_Pierre_inequality}.
 The above proof is formal when considering weak or very weak solutions, in which case, it is quite lengthy (although standard) to make it rigorous: we start by considering the integrated version of inequality \eqref{derivation_local_norm}, which follows by Definition \ref{def.local.weak.sol} of weak solution plus an integration by parts (that can be justified through approximation); we then conclude by a Grownwall-type argument.\\
The proof is concluded once we show that $ C\left( \phi \right) $ is bounded: choosing $\phi$ as in Lemma \ref{test.function.HP} we get
\begin{equation*}
 \phi^{\frac{-m}{1-m}} \left| \mathcal{L}_{\gamma, \beta}\left( \phi \right) \right|^{\frac{1}{1-m}} \le \ka_{10}  \left(\rho^{\gamma, \beta}_{x_0}(R)\right)^{-\frac{1}{1-m}},
\end{equation*}
where the constant $\ka_{10}>0$ does not depends on $x_0$ but only on $N,m, \gamma$ and $\beta$. Finally, we get
\[
\left( 1-m \right)C\left(  \phi \right)=(1-m) \left[ \int_{B_{2R}}{\phi\left( x \right)^{\frac{-m}{1-m}} \left|\mathcal{L}_{\gamma, \beta}\left(\phi\left( x \right)\right)\right|^{\frac{1}{1-m}}  |x|^{-\gamma} \dx} \right]^{1-m}  \le \frac{\ka_{10}}{\rho^{\gamma, \beta}_{x_0}(R)} \mu_\gamma(B_R)^{1-m}.
\]
 Using \eqref{rho.pseudo} and \eqref{r.sigma} one easily gets \eqref{Herrero.Pierre.123}. The proof is now concluded.\qed
When $p>1$ similar estimates still hold in the following slightly weaker form.
\begin{prop} \label{lemma_time}
Let $m\in (0,1)$ and $ u  $ be a nonnegative local strong solution to \ref{WFDE} on $(0,T]\times \Omega$.  Let $p>1$ and $u(\tau)\in \LL^p_{\gamma, \rm loc}(\Omega)$ for some $\tau\ge 0$. Let $x_0 \in \Omega$ and  $0 < R_1 < R_0$ be such that  $0 \not \in \overline{B_{R_0}(x_0)\setminus B_{R_1}(x_0)}$. Then for any $t\ge \tau$ we have
 \begin{equation}\label{lemma.time.ineq}
  \left[\int_{B_{R_1}\left( x_0 \right)}{u\left(t, x \right)^p \ |x|^{-\gamma} \dx} \right]^{\frac{1-m}{p}} \leq  \left[\int_{B_{R_0}\left( x_0 \right)}{u\left(\tau, x \right)^p \ |x|^{-\gamma} \dx} \right]^{\frac{1-m}{p}} + K_{R_0, R, p, \sigma, x_0}\left(t-\tau \right)
 \end{equation}
where the constant $ K_{R_0, R, p, \sigma, x_0}$ is given by
\begin{equation*}
 K_{R_0, R, p, \sigma, x_0}=c_p \frac{ h_{\sigma}\left( R_0, R_1, x_0 \right) }{\left(R_0-R_1\right)^{\sigma}}\left[ \mu_{\gamma}\left(B_{R_0}\left( x_0 \right) \setminus B_{R_1}\left( x_0 \right) \right)\right]^{\frac{1-m}{p}},
\end{equation*}
where $c_p\sim m(1-m)/(p-1)$ depends only on $p, m, N$.
\end{prop}
\noindent\textbf{Remark. }The above estimates prove stability of local $\LL^p_\gamma$ norms. Analogous estimates have been proven in \cite{BV-ADV,DaskaBook,DBbook, DiBenGianVesBook} in the non-weighted case and in \cite{BGV-JEE} on manifolds of non-positive curvature.

\noindent\textbf{Proof. }
The energy inequality \eqref{energy.inequality.1} can be  written as follows,
for any  $ 0\le \psi \in C_{c}^{\infty}{\left(\Omega \right)}$
\begin{equation*}\begin{split}
 \frac{d}{dt} \int_{\Omega}{u^p \psi^2  |x|^{-\gamma} \dx}
\leq  \ \frac{2mp}{p-1}\int_{\Omega}{u^{p+m-1}|\nabla \psi|^2 \ |x|^{-\beta} \dx}.
\end{split}\end{equation*}
 Using H\"older inequality with exponents $p/(1-m)$ and $p/(p+m-1)$ we get
\begin{align*}
 \int_{\Omega}u^{p+m-1}&|\nabla \psi|^2 \ |x|^{-\beta} \dx =  \int_{\Omega}{|\nabla \psi|^2  |x|^{-\beta}  |x|^{\frac{\gamma(p+m-1)}{p}} \psi^{\frac{-2(p+m-1)}{p}} \ u^{p+m-1} |x|^{\frac{-\gamma(p+m-1)}{p}} \psi^{\frac{2(p+m-1)}{p}}  \dx} \\
 \leq & \left[\int_{\Omega}{u^p \psi^2 \ |x|^{-\gamma} \dx} \right]^{1-\frac{1-m}{p}} \left[\int_{\Omega}{|\nabla \psi|^{\frac{2p}{1-m}} |x|^{-\beta \frac{p}{1-m}} |x|^{\gamma \frac{p+m-1}{1-m}} \psi^{- 2\frac{p+m-1}{1-m}} \dx} \right]^{\frac{1-m}{p}}
\end{align*}
Combining the above estimates we obtain the following differential inequality
\begin{equation*} \label{differential_inequality}
  \frac{d}{dt} \int_{\Omega}{u^p \psi^2  |x|^{-\gamma} \dx} \leq C_{\psi} \left[\int_{\Omega}{u^p \psi^2 \ |x|^{-\gamma} \dx} \right]^{1-\frac{1-m}{p}},
\end{equation*}
which, integrated over $\left(\tau, t \right)$ gives us
\begin{equation*}
 \left[\int_{\Omega}{u\left(t,x \right)^p \psi^2 \ |x|^{-\gamma} \dx} \right]^{\frac{1-m}{p}} \leq \left[\int_{\Omega}{u\left( \tau, x \right)^p \psi^2 \ |x|^{-\gamma} \dx} \right]^{\frac{1-m}{p}} + \frac{\left(1-m \right)}{p}C_{\psi}\left(t-\tau \right).
\end{equation*}
The above proof is formal when considering weak or very weak solutions: a rigorous derivation of the above inequality can be done in that case by using directly the energy inequality \eqref{energy.inequality.1} and a Grownwall type argument.
To conclude the proof, we need to show that $ C_{\psi} $ is finite. To this end, we choose $\psi=\phi^b$ with $0\le \phi\le 1$  so that $\supp(\psi)\subseteq B_{R_0}$, $\supp(|\nabla\psi|)\subseteq B_{R_0}\left( x_0 \right)\setminus B_{R_1}\left( x_0 \right):=A_{R_0,R_1}$; we take $b>p/(1-m)$ so that
\[
|\nabla\psi(x)|^{\frac{2p}{1-m}}\,\psi(x)^{-\frac{2(p+m-1)}{1-m}}
\le C_1 \phi^{2b-\frac{2p}{1-m}}|\nabla\phi(x)|^{\frac{2p}{1-m}}
\le C_2 (R_0-R_1)^{-\frac{2p}{1-m}}\,.
\]
As a consequence, we obtain
\begin{align*}
& \frac{1-m}{p}C_{\psi} =   \frac{2m(1-m)}{p-1} \left[\int_{A_{R_0,R_1}}\frac{|\nabla \psi|^{\frac{2p}{1-m}}}{\psi^{\frac{2(p+m-1)}{1-m}}}   |x|^{\gamma \frac{p+m-1}{1-m}-\beta \frac{p}{1-m}}  \dx \right]^{\frac{1-m}{p}}  \\
 &\leq   \left[\frac{C_2^{\frac{p}{1-m}}\, c_{0,p}^{\frac{p}{1-m}}}{(R_0-R_1)^{\frac{2p}{1-m}}}\int_{A_{R_0,R_1}} |x|^{\left(\gamma-\beta \right)\frac{p}{\left(1-m \right)}}    \frac{\dx}{|x|^{\gamma}} \right]^{\frac{1-m}{p}} \nonumber
\leq    c_p\frac{h_{\sigma}(R_0, R_1, x_0)}{\left(R_0-R_1\right)^{\sigma}}   \mu_{\gamma}\left(A_{R_0,R_1} \right) ^{\frac{1-m}{p}}
:= K_{R_0, R, p, \sigma, x_0}\,,
\end{align*}
where $c_p, c_{0,p}\sim m(1-m)/(p-1)$, and in the last step we have used inequalities \eqref{1111}, \eqref{2222}, \eqref{3333}, depending on the different cases.\qed

\subsection{Space-time smoothing effects for linear and nonlinear equations}
We prove here a weighted space-time $\LL^p_\gamma-\LL^\infty$ smoothing effects, through a Moser-type iteration. This result represents the core of the proof of our main upper estimates, Theorem \ref{local_upper_bounds}.  Here we will cover more general cases than the ones defined in (1), (2) and (3).

\begin{thm}\label{spacetime.smoothing.theorem}
    Let $ u\in \LL^p_{\gamma, \rm loc}(\left(0, T \right)\times B_R(x_0))$ be a nonnegative local strong (sub)solution to \ref{WFDE}, let $ p \geq 1 $ if $ m \in (m_c,1) $ and $ p > p_c $ if $ m \in(0, m_c]$. Let $x_0 \in \RR^N$, $ 0 < R_1 < R_0 < R $ be  such that $0 \not \in \overline{B_{R_0}(x_0)\setminus B_{R_1}(x_0)}$ and let $0 \leq T_0 < T_1 < T$. Then there exists a constant $\ka_{11}>0$ depending only on $\gamma, \beta, N, m,p$ such that the following inequality holds\vspace{-1mm}
    \begin{equation}\label{spacetime.smoothing} \begin{split}
        \sup\limits_{(\tau,y)\in (T_1, T] \times B_{R_1}(x_0)} \!\!\!\!  u(\tau,y) &\leq  \ka_{11}\!
        \left[\frac{h_{\sigma}\left( R_0, R_{1}, x_0 \right)}{\left(R_0 - R_{1}\right)^{\sigma}}\! +\! \frac{1}{T_{1} - T_0} \right]^{(N-\gamma+\sigma)\vartheta_p}
\left[ \int_{T_0}^T\int_{B_{R_0}(x_0)} \!\!\!\!\!\left(u^{p}+1\right) \frac{\dx\dt}{|x|^{\gamma}}  \right]^{\sigma\vartheta_p}
    \end{split}\end{equation}
where $\sigma$ and $\vartheta_p$ are defined in \eqref{definition.sigma.thethap},
$h_\sigma$ is defined in \eqref{function_h} and $\ka_{11}$ is given in \eqref{kappa11}.\vspace{-1mm}
\end{thm}
\begin{rem}\label{rem.GW}\rm This result is similar to Theorem 2.4 of \cite{BV-ADV}, but in this weighted case an important geometric factor appears: $h_\sigma$ defined in \eqref{function_h}. Since the weights are not translation invariant, the factor $h_\sigma$ will change in a strong way the behaviour of the local estimates in the different situations (a), (b), (c) and (d), described in Section \ref{result.organization}. The thechnical hypothesis $0 \not \in \overline{B_{R_0}(x_0)\setminus B_{R_1}(x_0)} $ assumed in Theorem \ref{spacetime.smoothing.theorem} gurantees that the quantity $h_\sigma\left( R_0, R_{1}, x_0 \right)$ is finite. In the non-weighted case, similar estimates  are proven in \cite{DaskaBook,DiBenGianVesBook,DGV}, with different operators and nonlinearities.
\end{rem}

\begin{rem}[The linear case with coefficients]\rm A close inspection of the proofs (both of the above Theorem and of the energy inequalities) reveals that indeed the above result still holds  in the limit $m\to 1^-$, and  even  for more general equations, see Proposition \ref{spacetime.smoothing.theorem.linear}; roughly speaking, we can consider solutions $v$ to the linear equation $v_t=\mathcal{L}_{\gamma, \beta}\, v$, where the prototype operator has the form $\mathcal{L}_{\gamma, \beta}\, v=|x|^{\gamma}\nabla\cdot \big(a(t,x)|x|^{-\beta}\nabla v\big)$ with $0<\lambda_0\le a(t,x)\le \lambda_1<\infty$. We refer to Subsection  \ref{sec.linear.regularity} for more details.
\end{rem}

\medskip

We now proceed with the proof of Theorem \ref{spacetime.smoothing.theorem}, which relies on a variant of the celebrated Moser iteration, adapting the proof of Theorem 2.4 of \cite{BV-ADV} to the weighted setting under consideration, for this reason we will be rather sketchy in the proofs. As already mentioned in Subsection \ref{result.organization}, the role of weighted Sobolev inequalities will be played here by the Caffarelli-Kohn-Nirenberg inequalities \eqref{WSIcomplete}, in the following form:\vspace{-1mm}

\begin{lem}[Iterative version of CKNI Inequality]
Let $r^*:=2(N-\gamma)/(N-2-\beta)$ with $ \gamma ,\beta $ as in \eqref{paramt.range}. Then for any ball $B_R(x_0)$  and for any $ f \in \LL^2(T_0, T_1 ;H^1_{\gamma, \beta} B_R( x_0))$ and  for any $a\in [1,r^*/2]$ the following inequality holds
 \begin{equation}
 \begin{split}\label{iterative_sobolev_inequality}
  \int\limits_{T_0}^{T_1}\int\limits_{B_R(x_0)} f^{2 a}\frac{\dy \dt}{|y|^{\gamma}}&\leq 2 \ S_{\gamma, \beta}^2 \left[   \int_{T_0}^{T_1}\int_{B_R(x_0)}f^{2} \frac{\dy \dt}{|y|^{\gamma}} +  \mu_\gamma(B_R(x_0))^\frac{\sigma}{N-\gamma}  \int_{T_0}^{T_1}\int_{B_R(x_0)}| \nabla f |^2  \frac{\dy \dt}{|y|^{\beta}}  \right]  \\
 & \hskip 5mm \times \sup_{t \in \left[ T_0, T_1 \right] }{\left(\mu_\gamma(B_R(x_0))^{-1}  \int_{B_R(x_0)}f^{2\left(a - 1 \right)q} \left(y, t \right) \ \frac{\dy }{|y|^{\gamma}} \right)}^{\frac{1}{q}}.
 \end{split}
 \end{equation}
where $q= r^*/(r^* - 2)=(N-\gamma)/\sigma$, $\sigma$ is  given  in \eqref{definition.sigma.thethap} and the constant $S_{\gamma, \beta}>0$ is the one appearing in \ref{WSIcomplete}.
\end{lem}
\noindent {\sl Proof.~}We will write $B_R$ instead of $B_R(x_0)$ and prove the result on  $\left(T_0, T_1 \right) \times B_R$.
 Since $r^*/2$ and $q$ are conjugate H\"older exponents, using H\"older and \eqref{WSIcomplete} inequalities we get
\begin{equation*}\begin{split}
 &\int_{B_R(x_0)}   {f^{2a}\frac{\dx}{|x|^{\gamma}}} =  \int_{B_R(x_0)}{f^{2} f^{2 \left(a - 1 \right)} \frac{\dx}{|x|^{\gamma}}} \leq \left[ \int_{B_R(x_0)}f^{\sr} \frac{\dx}{|x|^{\gamma}} \right]^{\frac{2}{\sr}} \left[ \int_{B_R(x_0)}{f^{2\left(a -1  \right)q} \frac{\dx}{|x|^{\gamma}}} \right]^{\frac{1}{q}} \\
 & \leq 2 S_{\gamma, \beta}^2 \left[ \mu_\gamma(B_R(x_0))^\frac{-\sigma}{N-\gamma}\int_{B_R(x_0)}{f^{2} \frac{\dx}{|x|^{\gamma}}} + \int_{B_R(x_0)}{| \nabla f|^{2} \ \frac{\dx}{|x|^{\beta}}} \right]  \sup_{t \in \left(T_0, T_1 \right)}{\left[ \int_{B_R(x_0)}{f^{2\left(a -1  \right)q} (x, t) \frac{\dx}{|x|^{\gamma}}} \right]^{\frac{1}{q}}}.
\end{split}\end{equation*}
Inequality \eqref{iterative_sobolev_inequality} follows by integrating in time.\qed

\noindent\textbf{Proof of Theorem \ref{spacetime.smoothing.theorem}. }  Throughout this proof $u$ will be a local strong solution to \ref{WFDE}  defined on $Q=(T_0, T] \times  B_R(x_0)$. We shall define $v(t, x)= u(t, x) \vee 1$, so that $v$ is a local strong \textit{subsolution} to the same equation. Notice that $v$ satisfies $u \leq v \leq u+1$ almost everywhere in $Q$. Let us fix $x_0\in \RR^N$\,, and simply denote $B_R=B_R(x_0)$ ($B_{R_1}=B_{R_1}(x_0)$ resp.) when there is no ambiguity. In what follows we will make some a priori estimates of the solution: the quantity $ h_{\sigma}$ will be involved and we can assume it to be bounded, as it will be clear later. We recall that $\sigma=2+\beta-\gamma>0$, $\sr=2(N-\gamma)/(N-2-\beta)$, and we set $q=(N-\gamma)/\sigma$.   We will split the proof into several steps. We first deal with the case $p>1$, the case $p=1$, which only affects the good fast diffusion range $m_c<m<1$, requires some extra work, and will be discussed at the last Step.

\noindent$\bullet~$\textsc{Step 1. }\textit{Preparation of the iteration step.}
Let $ 0 < R_1 < R_0 < R $ and $ 0 \leq T_0 < T_1 < T $ and define $ Q_0 := \left(T_0,  T\right] \times B_{R_0}\left( x_0 \right)  $ and $ Q_1 := \left(T_1,  T\right] \times B_{R_1}\left( x_0 \right) $. Recall that we are assuming $p>1$.  We are going to prove the following inequality:
\begin{align}\label{upper_iteration_first_step_4}
  \iint_{Q_1}v^{a \left(p+m-1 \right)} \frac{\dt\dx}{ |x|^{\gamma}}  \leq \left(2S_{\gamma, \beta}\right)^2 (2 c_1)^{1+\frac{1}{q}}\left[\frac{h_{\sigma}\left( R_0, R_1, x_0 \right)}{\left(R_0 - R_1 \right)^{\sigma}} + \frac{1}{T_1-T_0} \right]^{1+\frac{1}{q}}
  \left[\iint_{Q_0}{v^{p } \frac{\dt\dx}{ |x|^{\gamma}} } \right]^{1+\frac{1}{q}},
\end{align}
where $q= r^*/(r^* - 2)$,  $a \in (1, \sr/2)$ is such that $\left(p+m-1\right)\left(a-1 \right)q=p $; moreover, $ c_1>0$ only depends on $m,p,N$ and is given in the energy inequality \eqref{sup.energy.inequality.upper}.\\
To prove the above inequality, we first recall the CKNI inequality \eqref{iterative_sobolev_inequality} with $ f^2 = v^{p+m-1} $: \begin{equation}\label{upper_iteration_first_step_1}
\begin{split}
 \iint_{Q_1}{v^{a \left(p+m-1 \right)} \frac{\dx\dt}{|x|^{\gamma}} } &\leq 2 \ S_{\gamma, \beta}^2 \left[\iint_{Q_1}\left(v^{p+m-1} |x|^{-\gamma} + \mu_\gamma(B_{R_1})^\frac{\sigma}{N-\gamma} |\nabla v^{\frac{p+m-1}{2}} |^2  |x|^{-\beta}\right) \dx \dt \right]  \\
 &\times \left[\sup_{t \in \left(T_1, T\right]}\mu_\gamma(B_{R_1})^{-1}\int_{B_{R_1}\left( x_0 \right)}{v^{\left(p+m-1\right)\left(a-1 \right)q}|x|^{-\gamma} \dx} \right]^{\frac{1}{q}}.
\end{split}
\end{equation}
Next, we estimate the first term in the right-hand side using the upper energy inequality \eqref{sup.energy.inequality.upper} applied to $v \geq 1$, so that $v^{p+m-1} \leq v^p $ and
\begin{align}\label{upper_iteration_first_step_2}
& \iint_{Q_1}{ \left(v^{p+m-1} \ |x|^{-\gamma} + \mu_\gamma(B_{R_1})^\frac{\sigma}{N-\gamma} \ |\nabla v^{\frac{p+m-1}{2}} |^2 \ |x|^{-\beta}\right) \dx \dt}
\leq  \tilde J    \iint_{Q_0}{v^{p } \ |x|^{-\gamma} \dx \dt},
\end{align}
where we have assumed (without loss of generality\footnote{Indeed, by \eqref{radii.times.upper.iteration}, it is clear that this is true at the k-th iteration step for large $k$. Indeed we could have done everything by replacing $k$ with $k+a_0$ for a suitable large $a_0$: after the whole iteration process, this would only affect only estimate \eqref{kappa11}, where $C_1$ and $C_2$ would depend on $a_0$, but a posteriori the dependence on $a_0$ can be easily eliminated. We have decided to omit this, to focus on the main ideas.}) that
\begin{equation}\label{Correzione111}
\tilde J:=2 c_1 \mu_\gamma(B_{R_1})^\frac{\sigma}{N-\gamma}   \left[\frac{h_{\sigma}\left(R_0, R_1, x_0 \right)}{\left(R_0 - R_1 \right)^{\sigma}} + \frac{1}{T_1-T_0} \right]\ge 1.
\end{equation}
We next estimate the  second term  in the right-hand side of \eqref{upper_iteration_first_step_1}, using again \eqref{sup.energy.inequality.upper} applied to $v\ge 1$.  Since we are assuming $p>p_c$, we can choose $a \in (1, \sr/2)$ such that $\left(p+m-1\right)\left(a-1 \right)q=p $, to get
\begin{align}\label{upper_iteration_first_step_3}
 \left[\sup_{t \in \left(T_1, T\right]}\int\limits_{B_{R_1}\left( x_0 \right)}v^{\left(p+m-1\right)\left(a-1 \right)q} \frac{\dx}{|x|^{\gamma}} \right]^{\frac{1}{q}} & \leq \left[2 c_1    \left(\frac{h_{\sigma}\left( R_0, R_1, x_0 \right)}{\left(R_0 - R_1 \right)^{\sigma}} + \frac{1}{T_1-T_0} \right) \iint\limits_{Q_0}{ v^{p}  \frac{\dt\dx}{ |x|^{\gamma}} } \right]^{\frac{1}{q}}.
\end{align}
Combining inequalities \eqref{upper_iteration_first_step_1}, \eqref{upper_iteration_first_step_2} and \eqref{upper_iteration_first_step_3} we finally obtain \eqref{upper_iteration_first_step_4}.

\noindent$\bullet~$\textsc{Step 2. }\noindent\textit{The $k^{\rm th}$ iteration step. }We define a sequence of increasing exponents $p_k\to +\infty$ and nested cylinders $Q_k\subset Q_{k+1}$ as follows. We define first the exponents $p_k$, recalling that $q= r^*/(r^* - 2)$:
\begin{equation}\label{exponents.pk.upper.iteration}
p_{k+1}=\left(1+\frac{1}{q} \right)p_k + m-1=\left(1+\frac{1}{q} \right)^{k+1} \left[p_0-q\left(1-m \right) \right] + q\left(1-m \right)\,.
\end{equation}
Notice that in the range $m \in (0,m_c]$, we have $p_k<p_{k+1}\to +\infty$ if and only if $p_0>p_c$, while  when $m\in (m_c,1)$ it suffices to have $p_0 \geq 1$; this justifies the assumption on the initial datum.\\
Next, we define the cylinders $Q_k\subset Q= \left(T_0, T \right]\times B_{R_0}$ as follows:
\begin{equation}\label{Def.Qk.upper.iteration}
Q_{k}:= \left(T_k, T \right] \times B_{R_k}\left( x_0 \right)\qquad\mbox{such that}\qquad Q=Q_0 \supset Q_{k}\supset Q_{k+1}\to
Q_\infty=\left(T_\infty, T \right] \times B_{R_\infty}
\end{equation}
where we have chosen a decreasing sequence of radii $R_0>R_k>R_{k+1}\to R_\infty$ and and increasing sequence of times $T_0 <  T_k < T_{k+1} \to T_{\infty} $ such that
\begin{equation}\label{radii.times.upper.iteration}
R_{k}- R_{k+1}= C_1\frac{ R_{0}- R_{\infty}}{\left(k+1 \right)^{\alpha}}\qquad\mbox{and}\qquad T_{k+1}-T_{k} = C_2\frac{ T_{\infty}-T_0 }{\left(k+1 \right)^{\alpha \sigma}} \,.
\end{equation}
where ( note that we choose $\alpha$ so that $C_1, C_2$ will be finite)
\begin{equation}\label{pho_tau_definition_upper}
 \alpha = 2 \vee \frac{1}{\sigma}, \qquad C_1 =  \left(\sum_{k=0}^{\infty}{\frac{1}{\left(k+1 \right)^{\alpha}}} \right)^{-1}
 \qquad\mbox{and}\qquad C_2= \left(\sum_{k=0}^{\infty}{\frac{1}{\left(k+1 \right)^{\alpha \sigma}}} \right)^{-1}\,.
\end{equation}
Plugging all the above defined quantities in inequality \eqref{upper_iteration_first_step_4},  observing   that letting $p=p_k$ implies $p_{k+1}=a(p+m-1)$, we can write the $k^{\rm th}$ iteration step as follows

\begin{equation}\label{upper_iteration_first_step_6}
\begin{split}
 \left[\iint_{Q_{k+1}}v^{p_{k+1}}|x|^{-\gamma} \dx\dt
 \right]^{\frac{1}{p_{k+1}}}&\leq \left[2 S_{\gamma, \beta}\right]^{\frac{2}{p_{k+1}}} \left[ 2  c_1   \right]^{\left(1+\frac{1}{q} \right)\frac{ 1}{p_{k+1}}}   \left[ \frac{h_{\sigma}\left( R_k, R_{k+1}, x_0 \right)}{\left(R_k - R_{k+1} \right)^{\sigma}} + \frac{1}{T_{k+1}-T_k} \right]^{\left(1+\frac{1}{q} \right)\frac{ 1}{p_{k+1}}}\\
  & \times \left[\iint_{Q_k}{v^{p_k}|x|^{-\gamma} \dx\dt} \right]^{\frac{p_{k}}{p_{k+1}}\left(1+\frac{1}{q}\right) \frac{1}{p_{k}}}.\vspace{-1mm}
\end{split}
\end{equation}
\noindent\textit{Bounds for the constants. }It is convenient to bound the constants appearing in \eqref{upper_iteration_first_step_6} by a quantity which does not depend on $p$, but only on $m, \gamma, \beta, R_{\infty}$ and $R_0$. Recall the expression of $c_1$, given in \eqref{sup.energy.inequality.upper}
\begin{equation*}
c_1= 2K_{\psi} \ c_{m,p_k}^{-1} \qquad\mbox{with}\qquad c_{m, p_k} =  \frac{p_k-1}{p_k} \wedge \frac{2m \left(p_k-1 \right)^2}{\left(p_k+m-1 \right)^2}\,,
\end{equation*}
and with $K_\psi>0$ depending only on $N$. The quantity $c_{m, p_k}$ needs to be bounded uniformly for all $k\ge 0$: since $ p_k > p_0 > 1 \vee p_c $ it is easy to show that
\[
c_{m, p_k} \geq  \left(1-\frac{1}{p_0}\right) \wedge 2m \left(1 - \frac{m}{p_0+m-1} \right)^2 \quad\mbox{so that}\quad c_1=c_1(p_k)\le \overline{c}=\overline{c}(m,N,p_0)<+\infty\,.
\]
As a consequence, $c_1(p_k)\le \overline{c}$ can be bounded uniformly by a constant that depends only on $N, m,p_0$.

On the other hand, $ h_{\sigma}(R_k, R_{k+1}, x_0)  $ may be bounded by a fixed quantity depending only on $\sigma, R_0$ and $R_{\infty}$,
\begin{equation}\label{bounds.hsigma.upper}
 h_{\sigma}\left(R_k, R_{k+1}, x_0\right) \leq  h_{\sigma}\left( R_0, R_{\infty}, x_0 \right) (k+1)^{\alpha\left(2-\sigma\right)_{+}} \, C_1^{-\left(2-\sigma\right)_{+}},
\end{equation}
the constant $C_1$  being  the one appearing in \eqref{pho_tau_definition_upper}. \\
Finally, we can rewrite the $k^{\rm th}$ iteration step \eqref{upper_iteration_first_step_6} as follows
\begin{equation}\label{upper_iteration_first_step_7}
\begin{split}
 \left[\iint_{Q_{k+1}}v^{p_{k+1}}\right.&\,\big.|x|^{-\gamma} \dx\dt\Bigg]^{\frac{1}{p_{k+1}}}\leq I_{k+1}^{\frac{1}{p_{k+1}}}\left[\iint_{Q_k}{v^{p_k}|x|^{-\gamma} \dx \dt} \right]^{\frac{1}{p_{k+1}}\left(1+\frac{1}{q}\right)},
\end{split}
\end{equation}
where the constant $ I_{k+1} $ is bounded by,
\begin{equation}
\begin{split}\label{bound.I.k}
 I_{k+1} \leq& \left[2 S_{\gamma, \beta} \right]^2
 \left[2 \overline{c}
 \left(\frac{h_{\sigma}\left( R_0, R_{\infty}, x_0 \right) C_1^{-\frac{b}{\alpha}}}{(R_0-R_\infty)^\sigma} + \frac{C_2^{-1}}{T_\infty-T_0} \right) \right]^{1+\frac{1}{q}}
 \left(k+1 \right)^{b\left(1+\frac{1}{q}\right)} \\
 :=& J_0 \ J_1^{1+\frac{1}{q}} \left(k+1 \right)^{b\left(1+\frac{1}{q}\right)}.
\end{split}
\end{equation}
where $b=2\alpha$ if $0 < \sigma < 2$, $b=\sigma \alpha$ otherwise. The assumption \eqref{pho_tau_definition_upper} on $\alpha$ ensures that $b> 1$.

\noindent$\bullet~$\textsc{Step 3. }\noindent\textit{The iteration. }We now iterate the inequalities \eqref{upper_iteration_first_step_7} and obtain
\begin{equation}\label{general_k_step4}
  \left[\iint_{Q_{k+1}}v^{p_{k+1}}|x|^{-\gamma} \dx\dt\right]^{\frac{1}{p_{k+1}}} \leq \left[\prod_{j=1}^{k+1}{I_{j}^{\frac{1}{p_{k+1}}\left(1+\frac{1}{q} \right)^{k+1-j}}}\right]\left[ \iint_{Q_{0}}{v^{p_{0}} |x|^{-\gamma} \dx\dt}\right]^{\frac{1}{p_{k+1}}\left(1+\frac{1}{q} \right)^{k+1}}.
\end{equation}
Using inequality \eqref{bound.I.k} we can  estimate the first term  appearing in right-hand side of \eqref{general_k_step4}\,:
\begin{align*}\label{constant_in_the_limit}
 \prod_{j=1}^{k+1}{I_{j}^{\frac{1}{p_{k+1}}\left(1+\frac{1}{q} \right)^{k+1-j}}} & \leq \left[J_0 J_1^{1+\frac{1}{q}} \right]^{\frac{1}{p_{k+1}}\sum_{j=0}^{k}{\left(1+\frac{1}{q}\right)^{j}}}  \left[\left(k+1 \right)^{b'}\right]^{\frac{1}{p_{k+1}}}
 \left[k^{b'}\right]^{\frac{1}{p_{k+1}}\left(1+\frac{1}{q}\right)}
 \cdots \left[2^{b'}\right]^{\frac{1}{p_{k+1}}\left(1+\frac{1}{q}\right)^{k-1}}  \nonumber \\
 & = \left[J_0 J_1^{1+\frac{1}{q}} \right]^{\frac{1}{p_{k+1}}\sum_{j=0}^{k}{\left(1+\frac{1}{q}\right)^{j}}} \ \prod_{j=1}^{k+1}{j^{\frac{b'}{p_{k+1}}\left(1+\frac{1}{q} \right)^{k+1-j}}}\,,
\end{align*}
 where $b'=b\left(1+\frac{1}{q}\right)$.  Notice that there is a  constant $  c'>0 $ depending on  $p_0$, $ N $ and $ \gamma, \beta $ such that
\begin{equation*}
\lim_{k\to\infty}  \prod_{j=1}^{k+1}{j^{\frac{b'}{p_{k+1}}\left(1+\frac{1}{q} \right)^{k+1-j}}}
  \leq  \exp{\left(\frac{b'}{p_0-q\left(1-m \right)} \sum_{j=1}^{\infty}{\left(\frac{q}{q+1} \right)^j \log{j}} \right)}\leq  (c')^{\frac{q}{p_0-q(1-m)}}<+\infty\,.
\end{equation*}
Using the expression \eqref{exponents.pk.upper.iteration} of $p_k$, we see that
\begin{equation*}\label{limits}
 \lim_{k \rightarrow \infty}{\frac{\sum_{j=0}^{k}{\left(1+\frac{1}{q} \right)^j}}{p_{k+1}}}=\frac{q}{p_0-q\left(1-m \right)} \qquad \mbox{and} \qquad  \lim_{k \rightarrow \infty}{\frac{\left(1+\frac{1}{q} \right)^{k+1}}{p_{k+1}}}= \frac{1}{p_0-q\left(1-m \right)}.\vspace{-2mm}
\end{equation*}
We can now take the limit in \eqref{general_k_step4} as $k \rightarrow \infty$ obtaining (recall that $\lim_{k\to\infty}\|f\|_{\LL^{p_k}_\gamma(Q_k)}\ge \|f\|_{\LL^\infty(Q_\infty)} $)
\begin{equation}\label{final_result_special_subsolution}
 \|v\|_{\LL^{\infty}\left(Q_{\infty} \right)} \leq  (c'\,J_0)^{\frac{q}{p_0-q\left(1-m \right)}} J_1^{\frac{q+1}{p_0-q\left(1-m \right)}}  \ \left( \iint_{Q_{0}}{v^{p_{0}} \ |x|^{-\gamma} \dx\dt}\right)^{\frac{1}{p_0-q\left(1-m \right)}}.\vspace{-2mm}
\end{equation}
Recalling that $J_0$ and $J_1$ are as in \eqref{bound.I.k}, and  that $ v(t, x)= u(t, x) \vee 1 $, so that $ u^{p_0} \leq v^{p_0} \leq u^{p_0}+1 $, we  obtain from \eqref{final_result_special_subsolution}
\begin{equation}\label{final_fast_diffusion_space_time_upper_estimates}\begin{split}
 \sup_{Q_{\infty}}{u} &
 \leq \ka_{11}\left[\frac{h_{\sigma}\left( R_0, R_{\infty}, x_0 \right)}{\left(R_0 - R_{\infty}\right)^{\sigma}} + \frac{1}{T_{\infty} - T_0} \right]^{\frac{q+1}{p_0-q\left(1-m \right)}} \left[ \iint_{Q_{0}}{\left(u^{p_{0}}+1\right) \ |x|^{-\gamma} \dx\dt}\right]^{\frac{1}{p_0-q\left(1-m \right)}}\, ,\vspace{-1mm}
\end{split}\end{equation}
where the $\ka_{11}>0$ depends only on $m$, $p_0$, $N$, $ \gamma $ and $\beta$. We see that this is exactly inequality \eqref{spacetime.smoothing}, recalling that  $r^*:=2(N-\gamma)/(N-2-\beta)$, $q= r^*/(r^* - 2)$,  $\sigma\vartheta_{p_0}=1/(p_0-q(1-m))$  and
$(N-\gamma+\sigma)\vartheta_{p_0}=\frac{q+1}{p_0-q\left(1-m \right)}$ and that
\begin{equation}\label{kappa11}
\ka_{11}= \left[2 S_{\gamma, \beta} \right]^{\frac{2q}{p_0-q\left(1-m \right)}} \left[2 \overline{c}\,(C_1^{-\frac{b}{\alpha}}\vee C_2^{-1}) \,\right]^{\frac{q+1}{p_0-q\left(1-m \right)}}\vspace{-1mm}
\end{equation}
where $S_{\gamma, \beta}$ is as in Proposition \eqref{prop.Sobolev.Balls}, $C_1,C_2$ are as in \eqref{pho_tau_definition_upper}, $b=2\alpha$ if $0 < \sigma < 2$, $b=\sigma \alpha$ otherwise,  and $\overline{c}=\overline{c}(m,N,p_0)$,  is as in Step 2.  The proof of Theorem \ref{spacetime.smoothing.theorem} in the case $p_0>1$ is concluded by letting $p_0=p$, $R_\infty = R_1$ and $T_{\infty} = T_1$.

\noindent$\bullet~$\textsc{Step 4. }\textit{The case $p=1$. }So far, we have proven the space time smoothing effect for solutions in  $\LL^{p_0}_\gamma$ for any $p_0>1$. Unfortunately we cannot  simply take the limit as $p_0\to 1$  in inequality \eqref{final_result_special_subsolution} since the constant $c'$ would blow up, being proportional to $c_1\sim  (p_0-1)^{-1}$, explicitly given in \eqref{sup.energy.inequality.upper}. We show how to deal with the limiting case $p=1$.  A standard way to proceed is to first prove the result for bounded initial data, for instance $u_{0,n}=u_0\wedge n$, then by a lengthy but straightforward approximation procedure the result holds for $\LL^1_\gamma$ solutions.
We are going to use inequality \eqref{final_result_special_subsolution} which holds true on any couple of cylinders of the form $Q_{\infty} \subset \underline{Q}:=(\underline{T}, T] \times B_{\underline{R}}(x_0) \subset \overline{Q}:= (\overline{T}, T] \times B_{\overline{R}}(x_0) \subset Q_0$ for any $p_0>1$ , and implies (by H\"older's and Young's inequalities)
\begin{equation}\label{final_result_special_subsolution.2}\begin{split}
 \|v\|_{\LL^{\infty}\left(\underline{Q} \right)} &\le \frac{1}{2} \|v\|_{\LL^{\infty}\left(\overline{Q} \right)}+ c'' J_1^{\frac{q+1}{1-q\left(1-m \right)}}\|v\|_{\LL^{1}\left(\overline{Q} \right)}^{\sigma\vartheta_1}
\end{split}\end{equation}
where $J_1>0$ is as in \eqref{bound.I.k} and can be estimated as follows,
\begin{equation}\label{J1111}
J_1\le \overline{c}'\left[\frac{h_{\sigma}\left( \overline{R}, \underline{R}, x_0 \right)}{\left(\overline{R} - \underline{R}\right)^{\sigma}} + \frac{1}{\underline{T} - \overline{T}} \right]
\end{equation}
with $\overline{c}'= 2\,\overline{c}\, C_1^{-b/\alpha}\vee C_2^{-1}$, with $\overline{c}', C_1,C_2>0$, as in the previous step. Moreover,
\begin{equation}\label{c2222}
c''=2^{\frac{\sigma(p-1)}{\sigma-(N-\gamma)(1-m)}} (c'\,J_0)^{\frac{q}{1-q\left(1-m \right)}}>0
\end{equation}
and $c', J_0 >0$ are as in the previous step, and depend only on $p_0$, $ N $ and $ \gamma, \beta $.\\
We are now in the position to iterate the above inequality, using ideas inspired by a classical Lemma due to DeGiorgi, that can be found in many different sources, for instance see Lemma 3.6 of \cite{BGV}. Fix $0< \tau <  1$ and define cylinders $Q_i:= (t_i , T] \times B_{r_i}(x_0)$ where
\begin{equation*}\begin{split}
r_0=\underline{R} \,\,\,\,\,\, \mbox{and} \,\,\,\,\,\,\ r_{i+1}:=r_i + (1-\tau)\tau^{i}(\overline{R}-\underline{R}) \\
t_0=\underline{T} \,\,\,\,\,\,\ \mbox{and} \,\,\,\,\,\,\ t_{i+1}:=t_i - (1-\tau^\sigma)\tau^{i \sigma}(\underline{T}-\overline{T})\,,
\end{split}\end{equation*}
we iterate \eqref{final_result_special_subsolution.2} as follows
\begin{equation}\begin{split}
\|v\|_{\LL^{\infty}\left(Q_{0}\right)}&\le \frac{1}{2} \|v\|_{\LL^{\infty}\left(Q_{1}\right)}+ c'''\left[\frac{h_{\sigma}\left( \overline{R}, \underline{R}, x_0 \right)}{\left(\overline{R} - \underline{R}\right)^{\sigma}} + \frac{1}{\underline{T} - \overline{T}}\right]^{(q+1)\vartheta_1}\|v\|_{\LL^{1}\left(\overline{Q} \right)}^{\sigma\vartheta_1} \tau^{-i\sigma(q+1)\vartheta_1}\\
 &\le \left(\frac{1}{2}\right)^{k} \|v\|_{\LL^{\infty}\left(Q_{k} \right)}+c'''\left[\frac{h_{\sigma}\left( \overline{R}, \underline{R}, x_0 \right)}{\left(\overline{R} - \underline{R}\right)^{\sigma}} + \frac{1}{\underline{T} - \overline{T}}\right]^{(q+1)\vartheta_1}\|v\|_{\LL^{1}\left(\overline{Q} \right)}^{\sigma\vartheta_1}\sum_{i=0}^{k-1}\left(2\tau^{\sigma(q+1)\vartheta_1}\right)^{-i}
\end{split}\end{equation}
where $c''':= c'' \left(\frac{\overline{c}'}{\left[(1-\tau^\sigma)\wedge(1-\tau)^\sigma\right]}\right)^{(q+1)\vartheta_1}$. Taking the limit $k \rightarrow \infty$ we get
\begin{equation}
  \|v\|_{\LL^{\infty}\left(Q_{0}\right)}= \|v\|_{\LL^{\infty}\left(\underline{Q}\right)} \leq \ka_{11} \left[\frac{h_{\sigma}\left( \overline{R}, \underline{R}, x_0 \right)}{\left(\overline{R} - \underline{R}\right)^{\sigma}} + \frac{1}{\underline{T} - \overline{T}}\right]^{(q+1)\vartheta_1}\|v\|_{\LL^{1}\left(\overline{Q} \right)}^{\sigma\vartheta_1}
\end{equation}
Finally, the constant $\ka_{11}=c''' \sum_{k=1}^\infty \left(2\tau^{\sigma(q+1)\vartheta_1}\right)^{-i}  <\infty$, whenever $2^{-\sigma(q+1)\vartheta_1}<\tau<1$. We have proven inequality \eqref{final_fast_diffusion_space_time_upper_estimates} also when $u_0\in \LL^1$, and the only thing that changes is the constant $\ka_{11}$, which in any case only depends on $p_0$, $ N $ and $ \gamma, \beta $ and we can even fix $p_0>1$ taking for instance $p_0=2$.\qed

\subsection{Proof of Theorem \ref{local_upper_bounds}}\label{Proof.Thm2.1}
The aim of this section is to prove Theorem \ref{local_upper_bounds}: we only sketch the main points, which are analogous to those in Section 2.3 of \cite{BV-ADV}, where a detailed proof is given,  just for simplicity here we take $t_0= 3^{-\sigma} t$.   Let $ u(t, x) $ be a weak solution to \ref{WFDE} in the cylinder $\left(0, T \right) \times B_{2R_0}\left(x_0 \right)$, and let $ 0 < \tau < T $, and define $\varepsilon=1/3$ and $ \rho=(3/2)R_0$. We apply  Theorem \ref{spacetime.smoothing.theorem}  to the rescaled function
\begin{equation*}\label{defintion_rescaled_function}
\hat{u}\left(t,x \right)= M u\left(\tau t, \rho x \right) \qquad \mbox{where} \qquad M= \left(\frac{\rho^\sigma}{\tau} \right)^\frac{1}{1-m},
\end{equation*}
which turns out to be a solution to the \ref{WFDE} over the cylinder $Q= (0, 1] \times B_{(1+\varepsilon)}\left( \rho^{-1} x_0 \right)$.
Consider the cylinders  $Q_0=\left(0, 1\right] \times B_1\left( \rho^{-1}x_0 \right) $  and $ Q_1=\left((1/3)^\sigma, 1 \right] \times B_{1-\varepsilon}\left( \rho^{-1}x_0 \right)$. Applying estimate \eqref{spacetime.smoothing} to $\hat{u}$ over these two cylinders we get
\begin{equation*}
 \sup_{Q_1}\hat{u} \leq \ka_{11} \left[2\,3^{\sigma}h_{\sigma}\left( (3/2)R_0, R_0, x_0 \right)\right]^{\frac{q+1}{p-q\left( 1-m \right)}} \left[\iint_{Q_0}(\hat{u}^{p}+1) \ |x|^{-\gamma} \dx \dt\right]^{\frac{1}{p - q\left(m-1 \right)}}.
\end{equation*}
 Applying the inequalities obtained in Proposition \ref{Herrero_Pierre_Thm} and Proposition \ref{lemma_time} to $ \hat{u} $, on the domains $ B_1\left( \rho^{-1}x_0 \right) \subset B_{1+\varepsilon}\left( \rho^{-1}x_0 \right) $, for times $ t \in \left[0, 1 \right] $ and integrating them in $ t $ over $ \left(0, 1 \right) $ we obtain
\begin{equation*}
\int_{0}^{1}{\int_{B_1}{\hat{u}^{p} \ |x|^{-\gamma} \dx} \dt} \leq   2^{\frac{p+m-1}{1-m}} \ \int_{B_{1+\varepsilon}}{\hat{u}\left(0, x \right)^{p} \ |x|^{-\gamma} \dx} + \frac{\left(1-m\right)2^{\frac{p+m-1}{1-m}}}{p+1-m}\,\mathcal{K},
\end{equation*}
which holds for any $p > p_c $ if $m \in (0, m_c] $ or $p \ge 1$ if $m \in (m_c, 1)$; notice that $\mathcal{K}>0$ is an upper bound for the two constants given in \eqref{Herrero_Pierre_inequality} and \eqref{lemma.time.ineq} when $p=1$ or $p>1$ respectively; the expression of $\mathcal{K}$ will be given below in \eqref{constants.upper.bounds}.
Rescaling back to $ u $, using   inequality $(a+b)^{s}\leq k_1 a^s + k_2 b^s$ and recalling that $h_{\sigma}\left( 1, 1-\varepsilon,  \rho^{-1}x_0 \right)= h_{\sigma}\left(  \rho, \rho(1-\varepsilon), x_0 \right)  $ and $\mu_\gamma (B_1( \rho^{-1} x_0)=  \rho^{\gamma - N} \mu_\gamma (B_{ \rho}( x_0))$, we finally obtain
\begin{equation*}\label{inequality.1.thm2.2}
\sup_{\left((1/3)^\sigma \tau, \tau \right] \times B_{R_0}(x_0)}{u} \leq \frac{C_1}{\tau^{\frac{q}{p+q\left(m-1 \right)}}} \left[\int_{B_{2R_0}(x_0)}{u\left(0, y \right)^{p} |y|^{-\gamma} \dy} \right]^{\frac{1}{p - q\left(m-1 \right)}} + C_2\left[\frac{\tau}{R_0^{\sigma}} \right]^{\frac{1}{1-m}},
\end{equation*}
where the constants are (recall that $\varepsilon=1/3$ and $ \rho=(3/2)R_0$):
\begin{equation}\label{constants.upper.bounds}\begin{split}
 C_1 &= k_1  \ka_{11} 2^{\frac{p+m-1}{(1-m)(p + q\left(m-1 \right))}} \left[2\,3^{\sigma}h_{\sigma}\left(  (3/2)R_0, R_0,x_0 \right)\right]^{\frac{\left(q+1\right)}{p-q\left(1-m \right)}}, \\
 C_2&=k_2 \ka_{11} \left[ \frac{\left(1-m\right)2^{\frac{p+m-1}{1-m}}}{p+1-m}\,\mathcal{K} + \frac{\mu_\gamma (B_{ (3/2)R_0}( x_0))}{((3/2)R_0)^{N - \gamma }} \right]^{\frac{1}{p + q\left(m-1 \right)}} \left[2\,3^{\sigma}h_{\sigma}\left( (3/2)R_0, R_0 ,x_0 \right)\right]^{\frac{\left(q+1\right)}{p-q\left(1-m \right)}}, \\
 \mathcal{K}&= \left\{
\begin{array}{lll}
 c_{p} \, \left(3^{\sigma} h_{\sigma}\left( 2R_0 , (3/2)R_0, x_0 \right)\right)^{\frac{p}{1-m}}  (2R_0)^{\gamma-N}\mu_\gamma\left(B_{2R_0}(x_0)\right) &\mbox{if $ p > 1$ }\\
 \ka_{10}^{\frac{1}{1-m}}(2R_0)^{\gamma-N}\mu_\gamma\left(B_{2R_0}(x_0)\right) \left(\rho^{\gamma,\beta}_{x_0}(R_0/3)\right)^{-1}  &\mbox{if $p = 1$}, \,\\
\end{array}
\right.
\end{split}\end{equation}
where $ c_{p}$ is defined in \eqref{lemma.time.ineq}. This concludes the proof. \qed

\section{Part II. Positivity estimates}\label{sec:PART.II}
This second part of the paper is devoted to the proof of our main positivity result, Theorem \ref{LOCAL.LOWER.BOUNDS}. The proof is delicate, quite long and technical, and represents the major novelty of this paper, as already explained.
The strategy of the proof of our local lower bounds, Theorem \ref{LOCAL.LOWER.BOUNDS},  relies  on the study of the worst-case scenario: we will prove lower bounds for solutions to a ``smaller problem'', that we will call the \textit{Minimal Dirichlet Problem} (MDP), following \cite{BV-ADV}. Then, by  nowadays standard comparison arguments, see e.g.  \cite{HerreroPierre, JLVPorousMedioum}, we can extend the result to local nonnegative solutions.  Let us consider a the Minimal Dirichlet Problem, i.e. a homogeneous Dirichlet problem localized on a ball, with  a  smaller initial datum:\vspace{-2mm}
\begin{equation}\label{MINIMAL.DIRICHLET.PROBLEM}\tag{MDP}
\begin{cases}
\begin{aligned}
&\partial_t u  = |x|^{\gamma}\nabla \cdot \left(|x|^{-\beta}\nabla u^m \right)  &\qquad\, \text{in } Q_T=\left(0, T\right)\times B_{R_0}(x_0) \\
&u\left(t,x\right) = 0 &\qquad\, \text{for } t>0 \, \text{and}\, x \in \partial B_{R_0}(x_0) \\
&u\left(0,x\right) =u_0\chi_{B_R(x_0)}\ge 0 &\qquad\, \text{in $B_{R_0}(x_0)$, with $4R<R_0$}. \\
\end{aligned}
\end{cases}
\end{equation}
\noindent\textit{Extinction time for MDP and minimal life time. }We will show that nonnegative solutions to MDP extinguish in finite time $T=T(u_0)>0$. Moreover,  $T$ (hence its lower bound $t_*$) provides an estimate of the time interval for which any non-negative super-solution is strictly positive: recall that also super-solutions can extinguish in finite time. For this reason we call $t_*$ \textit{minimal life time }of the (super)solution $u$, following \cite{BV-ADV}. Estimating $T$ in terms of the initial datum (or of the solution at a reference time) will provide estimates on the size of the intrinsic cylinders (the natural domains for positivity and Harnack estimates, whose size depends on $u$) for any local super-solution.
 Let us state the main result of this part.\vspace{-2mm}
\begin{thm}[Interior Lower Bounds for MDP]\label{positivity.MDP.thm.general}
Let $0<4R= R_0$ and $u$ be the solution to \ref{MINIMAL.DIRICHLET.PROBLEM} corresponding to the initial datum $u_0\chi_{B_R(x_0)} \in \LL^p_\gamma(B_{R_0}(x_0))$  with $ p > p_c $ when $m\in (0,m_c]$ and with $p \geq 1  $ when $m\in (m_c,1)$, moreover assume that $B_{R_0}(x_0)$ satisfies either \textrm{(1)}, (2) or (3). Then, there exist $\kappa_*, \ka_{p,0}>0$ depending on $N,m,\gamma,\beta$, given in \eqref{t.star.MDP.L1} and  \eqref{Kp0} respectively,  such that we have the following estimates for the extinction time $T=T(u_0)$:
\begin{equation}\label{Bounds.FET}
t_{\ast}:=  \kappa_*R^{\sigma}\left(\frac{\|u_0\|_{\LL^1_\gamma(B_R)}}{\mu_\gamma(B_{R}(x_0))}\right)^{1-m}
\le T \le  \frac{ \mu_\gamma(B_{R_0}(x_0))^{\frac{\sigma}{(N-\gamma)}}}{\ka_{p,0}} \left(\frac{\|u_0\|_{\LL^p_\gamma(B_{R_0}(x_0))}}{\mu_\gamma(B_{R_0}(x_0))^{\frac{1}{p}}}\right)^{1-m}
\end{equation}
Moreover, there exists $\kb>0$ such that
\begin{equation}\label{LOCAL.LOWER.INEQUALITY.MDP}
    \inf_{x \in B_{2R}(x_0)} u(t,x) \ge \kb \left[\frac{\mu_\beta(B_R(x_0))}{\mu_\gamma(B_R(x_0))}\,\frac{t}{R^{2}}\right]^{\frac{1}{1-m}}\qquad\mbox{for any $t \in[0, t_*]$.}
\end{equation}
The constant $\kb>0$ depends on $N,m,\gamma,\beta$ and has an (almost) explicit expression given in \eqref{const.lower.cor.L1},  and only  when $0<m\le m_c$  it   depends also on $ \widetilde{H}_p$, defined in \eqref{def.HP}.
\end{thm}
\noindent\textbf{Remarks. }(i) Theorem \ref{positivity.MDP.thm.general}  holds   for more general scenarios, as it will be clear by a close inspection of the proof presented in this Part. However, to simplify the presentation we have decided to state the result  only  under assumptions (1), (2) or (3), since as already remarked, they represent the most relevant cases.\\
(ii) We recall that $\kb$ has a precise behaviour given in term of $\widetilde{H}_p$, and that when $m>m_c$ it does not depend on $\widetilde{H}_p$, see Remark \ref{rem.low.bdds.intro} (ii) and (iii) for more details. \\
(iii) A closer inspection of the proofs reveals that analogous results hold for solutions of the Dirichlet problem on arbitrary bounded domains $\Omega$ and general initial data. Bounds  on the extinction time similar to \eqref{Bounds.FET} have been obtained in \cite{BV-ADV} for the model equation, and in \cite{BGV-JEE} on Riemannian manifolds.

\medskip

\noindent\textbf{Strategy of the proof of positivity estimates. }As already explained above, it is sufficient to prove our lower estimates for solutions to the reduced problem MDP, and to avoid unnecessary technicalities, we will work with strictly positive solutions which solve a ``lifted problem'', cf. Subsection \ref{ssec.lifted.problem}. The proof of our main positivity result is quite complex, as already mentioned in Subsection \ref{ssec.111}, indeed, more standard techniques seem to fail to give quantitative estimates, hence we develop a new method, that we split it into four steps:
\[
\underbrace{\LL^{-\infty}\xrightarrow[\mbox{\footnotesize Moser iteration}]{\mbox{\footnotesize Step 1, Section \ref{sec:proof.lower.bounds}}}\LL^{-s}_\gamma
\xrightarrow[\mbox{\footnotesize Smoothing}]{\mbox{\footnotesize Step 2, Section \ref{ssec.Step2}}}\LL^{-\varepsilon}_\gamma
\xrightarrow[\mbox{\footnotesize Parabolic John-Nirenberg}]{\mbox{\footnotesize Step 3, Section \ref{ssec.parab.JN}}}
}_{\mbox{Prop. \ref{lower.smoothing.final.form} for Lifted Problem, then let $\delta\to 0^+$ for MDP (Cor. \ref{weak.positivity.MDP})}}
\LL^{\varepsilon}_\gamma\xrightarrow[\mbox{\footnotesize Smoothing}]{\mbox{\footnotesize Step 4, Section \ref{ssec.step4}}}\LL^1_\gamma\,.
\]
The first step (Subsection \ref{sec:proof.lower.bounds}) consists in proving a $\LL^{-\infty}-\LL^{-s}_\gamma$ estimate through a lower Moser-type iteration with negative exponents; due to the nonlinear character of our equation, such iteration does not allow  one  to reach all negative exponents, in contrast with  what happens in the linear case here we can only reach $-s<-(1-m)$. The second step (Subsection \ref{ssec.Step2}) consists in proving
quantitative $\LL^{-s}_\gamma-\LL^{-\varepsilon}_\gamma$ estimates for any $\varepsilon\in (0,1-m)$.
 Subsection \ref{ssec.parab.JN}  contains the third step,  a parabolic analog of the celebrated John-Nirenberg Lemma: our Lemma provides a Reverse H\"older inequality for small exponents, in the form of a $\LL^{-\varepsilon}_\gamma-\LL^{\varepsilon}_\gamma$ estimate,  and it holds   for solutions to the (lifted) MDP; the proof relies on the monotonicity properties of the solutions to the MDP, combined with a weighted version of the celebrated John-Nirenberg Lemma, which we borrow from \cite{HKM}. Corollary \ref{lower.smoothing.final.form} collects all the results of the first three steps, in the form of $\LL^{-\infty}-\LL^{\varepsilon}_\gamma$ smoothing effect for solutions $u_\delta$ to the \ref{LARGE.DIRICHLET.PROBLEM}. Next, letting $\delta\to 0$, we prove Corollary \ref{weak.positivity.MDP}, which is the analogous result of Corollary \ref{lower.smoothing.final.form} for solutions of the MDP. Subsection \ref{ssec.step4} contains the fourth and last step, namely $\LL^{\varepsilon}_\gamma-\LL^1_\gamma$ estimates, see Lemma \ref{positivity.MDP.thm}; gathering all the previous results, we finally obtain the $\LL^{-\infty}-\LL^1_\gamma$ estimates in Corollary \ref{weak.positivity.MDP.L1}. The proof of Theorems  \ref{positivity.MDP.thm.general} and \ref{LOCAL.LOWER.BOUNDS}  is contained in Section \ref{sec.final.proofs.positivity}.
\subsection{Basic properties of solutions to the Minimal Dirichlet Problem}
We summarize here the standard properties of solutions to the \ref{MINIMAL.DIRICHLET.PROBLEM} which will be used in what follows.
\begin{prop}\label{SmoothingEffectDirichletProblem}
Let $u_0 \in \LL^p_{\gamma}(B_{R_0}(x_0))$ with $ p > p_c $ when $m\in (0,m_c]$ and with $p \geq 1  $ when $m\in (m_c,1)$. Then there exists a unique strong solution to the problem \ref{MINIMAL.DIRICHLET.PROBLEM} and the following properties  hold:
  \begin{itemize}[leftmargin=7mm]\itemsep1pt \parskip1pt \parsep0pt
    \item[\rm (i)]  There exists $\ka_{12}=\ka_{12}(\gamma, \beta, N, m, p)>0$ such that for any $ t > 0$
    \begin{equation}\label{SmoothingEffectInequality}
        \|u(t)\|_{\LL^{\infty}(B_{R_0}(x_0))}\leq \ka_{12}\frac{\|u_0\|^{\sigma p \vartheta_p}_{\LL^p_{\gamma}(B_R(x_0))}}{t^{(N-\gamma)\vartheta_p}},
    \end{equation}
    where $\vartheta_p= \frac{1 }{\sigma p - \left( N-\gamma \right)\left( 1-m \right)}$ and $\sigma=2+\beta-\gamma$.
    \item[\rm (ii)] For all $t > 0$ we have $\|u(t)\|_{\LL^p_{\gamma}(B_{R_0}(x_0))} \leq \|u_0\|_{\LL^p_{\gamma}(B_{R_0}(x_0))}$.
    \item[\rm (iii)] For all $t>0$, the function $t \rightarrow u(t,x)t^{-\frac{1}{1-m}}$ is non-increasing, for almost any $x \in B_{R_0}(x_0)$.
  \end{itemize}
\end{prop}
\noindent\textbf{Remark. }The constant $\ka_{12}$ may be chosen to be independent of $B_{R_0}(x_0)$\,, as explained in the proof.

\noindent {\bf Proof.~}Existence and uniqueness of strong solutions follow by minor modifications  of standard  arguments, cf. \cite{JLVPorousMedioum}. Also properties (ii) and (iii) follow by standard arguments that can be found in \cite{JLVPorousMedioum}. The upper bounds (i) hold as a consequence of the smoothing effects for the Cauchy problem on the whole space, namely inequality \eqref{SmoothingEffectInequality} with the norms taken on $\RR^N$. The proof of such global estimates is easier than its local counterpart: on one hand, such bounds can be proven directly by doing a Moser iteration, and then noticing that the solution the MDP (extended to zero outside of $B_{R_0}(x_0)$) is a subsolution to the Cauchy problem for \ref{WFDE} posed on $\RR^N$. On the other hand, such upper bounds can be deduced from the local upper estimates of Theorem \ref{local_upper_bounds} by letting $R_0\to \infty$ in inequality \eqref{upper.provisional.2} (in such limit the constant $\ka_8$ becomes independent of $R_0$ and the second term vanish), so we obtain \eqref{SmoothingEffectInequality}; in the latter case, assumptions (1),(2), (3) do not play an essential role, since we can always consider $x_0=0$.\qed
\subsection{Proof of the bounds \ref{Bounds.FET} on the extinction time for MDP. }
We now prove the two sided estimate \ref{Bounds.FET} on the extinction time. While the lower bounds follows from Proposition \ref{Herrero_Pierre_Thm}, the upper bounds requires the following Proposition, in which we show that Sobolev and Poincar\'e inequalities imply extinction in finite time for solutions to the MDP, as already observed in \cite{BGV-JEE,BV-ADV} for the model equation ($\beta=\gamma=0$) both in Euclidean and Riemannian manifolds settings.
\begin{prop}Let $u$ be the solution to the \ref{MINIMAL.DIRICHLET.PROBLEM}, corresponding to $u_0 \in \LL^p_{\gamma}(B_{R_{0}}(x_0))$ with $ p \ge p_c\vee 1 $. Then, for all $q>1$ there exists a constant $\ka_q>0$, such that for all $0\le \tau \le t $ we have
\begin{equation}\label{formual.Lp}
    \|u(t)\|_{\LL^q_{\gamma}(B_{R_{0}}(x_0))}^{1-m}   \le \|u(\tau)\|_{\LL^q_{\gamma}(B_{R_{0}}(x_0))}^{1-m} -\ka_q (t-\tau),
\end{equation}
as long as the right-hand side is nonnegative  where the constant $\ka_q$ is given by
\begin{equation}\label{Kp0}
  \ka_q =\ka_{q,0}\mu_\gamma(B_{R_{0}}(x_0))^{-\frac{\sigma}{(N-\gamma)}\left(1-\frac{p_c}{q}\right)}\qquad\mbox{with}\qquad \ka_{q,0}=\frac{4(q-1)(1-m)}{\ka_{13}^{2} \,(q+m-1)^2}\,,
\end{equation}
where $\ka_{13}$ depends only on geometrical quantities, but not on ${R_{0}}$; we recall that $\sigma=2+\beta-\gamma$.\end{prop}
\noindent {\bf Proof.~}We will write $B_{R_0}$ instead of $B_{R_0}(x_0)$ and $\LL^q_\alpha$ instead of $\LL^q_\alpha(B_{R_0})$ when no confusion arises. We combine weighted Sobolev and Poincar\'e inequalities \eqref{w_sobolev} via   H\"older's inequality  as follows: \begin{equation}\label{sobolev.on.domains}
  \|f\|_{\LL^s_{\gamma}}  \le \|f\|_{\LL^2_\gamma}^{1-\theta} \|f\|_{\LL^{\sr}_\gamma}^{\theta} \le \ka_{13} \mu_\gamma(B_{R_{0}})^{\frac{\sigma}{2(N-\gamma)}(1-\theta)} \|\nabla f\|_{\LL^2_{\beta}}.
\end{equation}
for any $s \in (2, \sr)$ for any function $f \in \D_{\gamma, \beta} $. Letting now $s=\frac{2q}{q+m-1}$, $f= u^{\frac{q+m-1}{2}}$ and $\theta=p_c/q$, and recalling that $s < \sr $ if and only if $q  > p_c$, we have that
inequality \eqref{sobolev.on.domains} implies
\begin{equation*}
 \|f\|_{\LL^{s}_{\gamma}}^2= \|u\|_{\LL^q_{\gamma}}^{q\left[1-\frac{1-m}{q}\right]} \le \ka_{13}^2 \mu_\gamma(B_{R_{0}})^{\frac{\sigma}{(N-\gamma)}(1-\frac{p_c}{q})} \|\nabla u^{\frac{q+m-1}{2}}\|_{\LL^2_{\beta}}^{2}
\end{equation*}
Next, we formally  take the derivative  $\LL^q_\gamma$ norm of $u$ to get the following differential inequality
\begin{equation*}
\frac{d}{dt}\|u\|_{\LL^q_{\gamma}}^q= -\frac{4q(q-1)}{(q+m-1)^2} \|\nabla u^{\frac{q+m-1}{2}}\|_{\LL^2_{\beta}}^{2} \leq -\frac{4q(q-1)}{(q+m-1)^2} \ka_{13}^2 \mu_\gamma(B_{R_{0}})^{-\frac{\sigma}{(N-\gamma)}(1-\theta)} \|u\|_{\LL^q_{\gamma}}^{q\left(1-\frac{1-m}{q}\right)}
\end{equation*}
Integrating this last inequality in $[\tau,t]\subset [0,T]$  we get \eqref{formual.Lp}. A rigorous proof (long and technical, but nowadays standard)  can be done by using energy inequalities and Grownwall-type arguments.\qed
\noindent {\bf Proof of Inequalities \ref{Bounds.FET}.~}As an immediate corollary of inequalities \eqref{formual.Lp}, there exist the extinction time $T>0$; moreover, it satisfies the upper bound \ref{Bounds.FET}, which is nothing but inequality \eqref{formual.Lp} with $\tau=0$ and $t=T$.
The  lower bound  follows by letting $t=0$ and $\tau=T$ in formula \eqref{Herrero.Pierre.123}.\qed
\subsection{Lifted problem and a first positivity result}\label{ssec.lifted.problem}
In this section we address the question of proving quantitative estimates of positivity for the \ref{MINIMAL.DIRICHLET.PROBLEM}.
We begin by introducing the following  ``lifted'' problem: let $\delta>0$, $0<R\le R_0$
\begin{equation}\label{LARGE.DIRICHLET.PROBLEM}\tag{$\delta$-MDP}
\begin{cases}
\begin{aligned}
&\partial_t u_{\delta}  = |x|^{\gamma}\nabla \cdot \left(|x|^{-\beta}\nabla u^m_{\delta} \right)  &\qquad\, \text{in }\left(0, \infty\right)\times B_{R_0}(x_0) \\
&u_{\delta} = \delta &\qquad\, \text{on }\left(0, \infty\right)\times \partial B_{R_0}(x_0) \\
&u_{\delta}\left(0,x\right) =u_0(x)\chi_{B_{R}(x_0)}(x) + \delta &\qquad\, \text{for }x\in B_{R_0}(x_0). \\
\end{aligned}
\end{cases}
\end{equation}
The results of Proposition \ref{SmoothingEffectDirichletProblem}  hold  also for solutions to the \ref{LARGE.DIRICHLET.PROBLEM}, more precisely, (ii) and (iii) hold in the same form, while (i)  holds   with an extra $\delta$ factor on the right-hand side, as explained below; the proofs of the latter facts  are a straightforward modification of the case $\delta=0$, indeed, it is clear that $v_\delta=u_\delta-\delta $, solves a homogeneous Dirichlet problem with a regularized nonlinearity: $\partial_t v_\delta=|x|^{\gamma}\nabla \cdot \left(|x|^{-\beta}\nabla (v_\delta^m + \delta^m)-\delta^m \right)$, as in Appendix B of \cite{BFR}, or as in Chapter 5 of \cite{JLVPorousMedioum}, where even more general nonlinearities are treated. However, we recall explicitly the two main estimates that we will use in what follows: the time monotonicity property (iii) of Proposition \ref{SmoothingEffectDirichletProblem},  namely the fact that  $t \rightarrow u_\delta(\cdot,t)t^{-\frac{1}{1-m}}$ is non-increasing in $t$; the upper bounds, namely, there exists $\ka_{12}$, given in \eqref{SmoothingEffectInequality}, such that for any $(t, x) \in (0, \infty) \times B_{R_0}(x_0)$
\begin{equation}\label{uniform.boundedness.approximate.solutions}
\delta \le u_{\delta}(t,x) \le\frac{ 2\ka_{12}}{t^{(N-\gamma)\vartheta_p}}\|u_0\chi_{B_{R}(x_0)} + \delta\|^{p\sigma \vartheta_p}_{\LL^p_{\gamma}(B_{R_0}(x_0))} + 2\delta := M_p(u_0, \delta, t).
\end{equation}
The next Proposition proves $\LL^{-\infty}-\LL^{\varepsilon}$ estimates for solutions to \ref{LARGE.DIRICHLET.PROBLEM}; this is the core of the proof of the main Theorem \ref{positivity.MDP.thm.general} and it contains the first three steps explained above.
In what follows, we need to assume a certain weighted $\LL^p$ integrability on the initial datum,  otherwise the above inequality may  fail, as  thoroughly explained before. It is convenient, even if not strictly needed, to assume  $u_\delta$ to be bounded: this will simplify the proofs of the lower iteration.
\begin{prop}[$\LL^{-\infty}-\LL^{\varepsilon}$ estimates for \ref{LARGE.DIRICHLET.PROBLEM}]\label{lower.smoothing.final.form}
Let $u_\delta$ be a strong (super)solution to \ref{LARGE.DIRICHLET.PROBLEM} on $(0, \infty) \times B_{R_0}(x_0)$, corresponding to  $u_0\chi_{B_R(x_0)} \in \LL^p_\gamma(B_{R_0}(x_0))$ with $ p > p_c $ when $m\in (0,m_c]$ and with $p \geq 1  $ when $m\in (m_c,1)$.
For any $0 \le T_0 < T_1 < T_2<T_3 $, such that  $T_3-T_2\ge T_1-T_0$,  for any  $x_0 \in \RR^N$ and for any $0<R_2<R_1<R$ such that $0 \not \in \overline{B_{R_1}(x_0)\setminus B_{R_2}(x_0)}$, and for any $\varepsilon>0$  such that
\begin{equation}\label{nu.delta.def.Final3}
0 < \varepsilon <   \nu_\delta\wedge (1-m)\qquad\mbox{with}\qquad\nu_\delta = \frac{1}{\ka_{14} \, \ka_{6}} \left[ 1 + \frac{R_1^\sigma}{T_0}\left(\frac{|x_0|}{R_1}\vee 1\right)^{\beta-\gamma} M_p(u_0, \delta, T_0)^{1-m}\right]^{-1/2},
\end{equation}
there exist $s_\varepsilon>1-m$ and $\kb_2>0$ such that
\begin{align}\label{formula.inf.finalBBBB}
  \inf_{(T_1, T_2]\times B_{R_2}(x_0)}{u_\delta}
  &\geq \kb_2\left[\left(1+M_{p}(u_0,\delta, T_0)^{1-m}\right)\right]^{\eta_{\varepsilon}+\frac{s_\varepsilon}{s_\varepsilon+m-1}\zeta_\varepsilon }\left(T_2-T_0\right)^{-\frac{1}{s_\varepsilon+m-1}} \nonumber \\
  &\times\bigg[
  \left(\mu_\gamma(B_{R_2}(x_0))^{-\frac{\sigma}{N-\gamma}}\vee\left(\frac{h_\sigma(R_1, R_2, x_0)}{(R_1-R_2)^\sigma}+ \frac{1}{T_1-T_0}\right)\right) \bigg]^{\eta_{\varepsilon}+\frac{s_\varepsilon}{s_\varepsilon+m-1}\zeta_\varepsilon }\nonumber \\
  & \times \left[\left(\frac{T_0}{T_3}\right)^{\frac{1}{1-m}}
   \mu_\gamma(B_{R_1}(x_0))^{-\frac{2}{s_\varepsilon}}
  \left( \int_{B_{R_1}(x_0)}u_\delta(T_3, x)^{\varepsilon}|x|^{-\gamma} \dx \right)^{\frac{1}{\varepsilon}}\right]^{\frac{s_\varepsilon}{s_\varepsilon+m-1}}
\end{align}
where $\ke$ is the smallest $k\in \NN$ such that $(r^*/2)^{k}\varepsilon>1-m$ and
\begin{align}\label{exponents.lower.smoothing.Final3}
 s_\varepsilon& :=\left(\frac{r^*}{2}\right)^{\ke}\varepsilon  \,, \quad  \eta_{\varepsilon}:={-\frac{1}{s_\varepsilon+m-1}\left(\frac{N-\gamma}{2+\beta-\gamma}+1 \right)}<0\,,\quad
\zeta_\varepsilon:={-\frac{1}{\varepsilon}\frac{1-\left(\frac{2}{r^*}\right)^{\ke}}{1-\frac{2}{r^*}}}<0\,,\
\end{align}
with $h_\sigma$, $M_p$ and $\sr$   defined in \eqref{function_h}, \eqref{uniform.boundedness.approximate.solutions} and  \eqref{w_sobolev} respectively;
we also provide the expression of  the constant  $\kb_2=\kb_3\kb_4^{\frac{s_\varepsilon}{s_\varepsilon+m-1}}\ka_{7}^{-\frac{2}{s_\varepsilon+m-1}}2^{(2\vee\sigma)(\eta_{\varepsilon}+\frac{s_\varepsilon}{s_\varepsilon+m-1}\zeta_\varepsilon)} $ with $\kb_3,\kb_4,\ka_7$ depending  only on $N,m,p,\beta,\gamma,\varepsilon$ and given in \eqref{formula.inf}, \eqref{const.finite.iter.negative.norms}, \eqref{rev.hold.cor.intro} respectively;   $\ka_{6}$  given in Lemma \ref{JOHN.NIRENBERG} depends on $N, \gamma$, $\ka_{14}$ given in Lemma \ref{Caccioppoli.JohnNirenberg.lemma.LDP} depends on $N, \gamma, \beta, m$.
\end{prop}
The next four subsections will be devoted to the proof of the above Proposition.\vspace{-2mm}
\subsection{Step 1. Lower Moser Iteration. }  \label{sec:proof.lower.bounds}
In this section we prove $\LL^{-\infty}-\LL^{-s}$ smoothing effects by means of a lower Moser-type iteration. We will cover more general cases than (1), (2) or (3). We are going to prove a priori estimates of the solution to the \ref{LARGE.DIRICHLET.PROBLEM}, and such bounds involve the quantity $ h_{\sigma}$ that under our assumptions will always be bounded.\vspace{-2mm}
\begin{prop}[Nonlinear case]\label{lower.spacetime.smoothing}Let $u$ be a strong (super)solution to \ref{LARGE.DIRICHLET.PROBLEM} on $(T_0, T_2) \times B_{R_0}(x_0)$, corresponding to $u_0 \in \LL^{p}_{\gamma}(B_{R_0}(x_0))$ with $ p > p_c $ when $m\in (0,m_c]$ and with $p \geq 1  $ when $m\in (m_c,1)$.
Then, for any $s>1-m$, for any $0 < R_2 < R_1 \le R_0$  such that $0 \not \in \overline{B_{R_1}(x_0)\setminus B_{R_2}(x_0)}$  and for any $T_1\in(T_0, T_2)$, there exists $\kb_3>0 $ such that\vspace{-2mm}
\begin{align}\label{formula.inf}
  \inf_{(T_1, T_2]\times B_{R_2}(x_0)}{u} &\geq \kb_3 \left[\left(1+M_{p}(u_0,\delta, T_0)^{1-m}\right)
  \left(\frac{h_{\sigma}\left( R_1, R_2, x_0 \right)}{\left(R_1-R_{2}\right)^{\sigma}} + \frac{1}{T_{1}-T_0}\right) \right]^{\eta_s}\nonumber \\
  & \hskip 0.5cm \times \left(\frac{T_0}{T_2}\right)^{\frac{s}{(1-m)(s+m-1)}}\left(T_2-T_0\right)^{-\frac{1}{s+m-1}}\left(\int_{B_{R_1}(x_0)}u^{-s}(T_2, x) |x|^{-\gamma} \dx\right)^{-\frac{1}{s+m-1}}
\end{align}
where $\eta_s:={-\frac{1}{s+m-1}\left(\frac{N-\gamma}{2+\beta-\gamma}+1 \right)}$, $h_\sigma$ and $M_p$ are defined in \eqref{function_h} and \eqref{uniform.boundedness.approximate.solutions} respectively, and $\kb_3>0$ depends on $s, p, N, \gamma, \beta, m$.
\end{prop}
\noindent\textbf{Remarks. }(i) We recall that the technical assumption $0 \not \in \overline{B_{R_1}(x_0)\setminus B_{R_2}(x_0)}$ is needed only to guarantee that the quantity $h_\sigma$ is finite.\\
(ii) The above estimate also holds when $m=1$, in which case we can take any $s>0$, see Proposition \ref{cor.lin.lower}. This fact considerably simplifies the proof of lower bounds in the linear case:  the lower bound is formulated in terms of a  space-time  integral on a bigger parabolic cylinder, and this is the classical result in linear parabolic equations, see for instance \cite{CS-AA87,CS-RSMUP85,CS-AMPA84,CS-CPDE84,GW,GW2,Moser,MoserCpam71}. We refer to Section \ref{sec.linear.regularity} for further details.

\smallskip

\noindent {\bf Proof of Proposition \ref{lower.spacetime.smoothing}.~}  Throughout the proof $u_\delta$ will be a strong (super)solution to \ref{LARGE.DIRICHLET.PROBLEM}  on a generic cylinder $Q:=(T_0, T) \times B_R(x_0)$ where $T$ and $R$ are fixed, and will be chosen at the end of the proof; we will write $u=u_\delta$, and $B_R=B_R(x_0)$ ($B_{R_1}=B_{R_1}(x_0)$ resp.) when no confusion arises.  We recall that $\sigma=2+\beta-\gamma>0$, $\sr=2(N-\gamma)/(N-2-\beta)$, and we set $q=(N-\gamma)/\sigma$. We will split the proof into several steps.

\noindent$\bullet~$\textsc{Step 1. }\textit{Preparation of the iteration step.}
Let $ 0 < R_1 < R_0 $ and $ 0 \leq T_0 < T_1 < T $ and define $ Q_0 := Q $ and $ Q_1 := \left(T_1,  T\right] \times B_{R_1}\left( x_0 \right) $.
We are going to prove the following inequality, with $a \in (1, \sr/2)$:\vspace{-1mm}
\begin{equation}\label{lower.first.inequality.4}
\begin{split}
\iint_{Q_1}u^{a(m-p-1)}|x|^{-\gamma} \dx \dt  & \leq  \left[2 S_{\gamma, \beta}\right]^2 \left[2 c(m, p)\right]^{1+\frac{1}{q}}   \left[\frac{h_{\sigma}\left( R_0, R_1, x_0 \right)}{\left(R_0 - R_1 \right)^{\sigma}} + \frac{1}{T_1-T_0} \right]^{1+\frac{1}{q}}  \\
  & \times  \left[1+M_{\widetilde{p}}(u_0, \delta, T_0)^{1-m}\right]^{1+\frac{1}{q}} \left[ \iint_{Q_0}{u^{-p+m-1}|x|^{-\gamma} \dx \dt} \right]^{1+\frac{1}{q}}.\vspace{-1mm}
\end{split}
\end{equation}
To prove the above inequality, we first recall the CKNI inequality \eqref{iterative_sobolev_inequality} with  $ f^2 = u^{m-p-1}$ and  $p>0$:\vspace{-1mm}
\begin{equation}\label{lower.first.inequality.1}
\begin{split}
 \iint_{Q_1}{u^{a \left(m-p-1 \right)} \frac{\dx\dt}{|x|^{\gamma}}  } &\leq 2  S_{\gamma, \beta}^2  \left[  \iint_{Q_1}{ \left(u^{ m-p-1}  |x|^{-\gamma} + \mu_\gamma(B_{R_1})^\frac{\sigma}{N-\gamma} \left|\nabla u^{\frac{m-p-1}{2}} \right|^2 \ |x|^{-\beta}\right) \dx \dt} \right]  \\
 & \times \mu_\gamma(B_{R_1})^{-\frac{1}{q}}\left[\sup_{t \in \left(T_1, T\right]} \int_{B_{R_1}}{u^{\left(m-p-1\right)\left(a-1 \right)q} |x|^{-\gamma} \dx} \right]^{\frac{1}{q}}.\vspace{-2mm}
 \end{split}
\end{equation}
We are going to estimate the two terms on the right-hand side of the above inequality, by means of a modified form of the lower energy estimates \eqref{sup.energy.inequality.lower}, which easily follows by using \eqref{uniform.boundedness.approximate.solutions}: for any $p > 0$ and any $T_1 \in (T_0, T)$ and $R_1 \in (0, R)$ we have\vspace{-1mm}
\begin{equation}\label{sup.energy.inequality.lower.approximate}
\begin{split}
 \sup\limits_{\tau \in \left[T_1, T \right]} \int\limits_{B_{R_1}(x_0)}{u(\tau, x)^{-p} \frac{\dx}{|x|^{\gamma}}}  + \int\limits_{T_1}^T\int\limits_{B_{R_1}(x_0)}{\left|\nabla u^{\frac{-p+m-1}{2}} \right|^2  \frac{\dx \dt}{|x|^{\beta}}}  \leq c_3 & \left[1+M_{\widetilde{p}}(u_0, \delta, T_0)^{1-m}\right] \\
   \times\left[\frac{h_\sigma(R, R_1, x_0)}{\left(R-R_1 \right)^\sigma} + \frac{1}{T_1-T_0} \right]
  \int\limits_{T_0}^T\int\limits_{B_R(x_0)}{u^{-p+m-1} \frac{\dx \dt}{|x|^{\gamma}}}&,
\end{split}
\end{equation}
where $M_{\widetilde{p}}(u_0, \delta, T_0)$ is defined in \eqref{uniform.boundedness.approximate.solutions}, and $ \widetilde{p} > p_c $ when $m\in (0,m_c]$ and with $\widetilde{p} \geq 1  $ when $m\in (m_c,1)$;  $h_\sigma$ is defined in \eqref{function_h}, and $c_3>0$ depends on $m,p,N$, as in and is given in the energy inequality \eqref{sup.energy.inequality.lower}.
We estimate now the first term of \eqref{lower.first.inequality.1} using \eqref{sup.energy.inequality.lower.approximate}:
\begin{equation}\label{lower.first.inequality.2}
\begin{split}
\iint_{Q_1} \Bigl(u^{ m-p-1}  |x|^{-\gamma}  &+  \mu_\gamma(B_{R_1})^\frac{\sigma}{N-\gamma}  \Bigl|\nabla u^{\frac{m-p-1}{2}} \Bigr|^2 \ |x|^{-\beta}\Bigr) \dx \dt
\leq J\iint_{Q_0}{u^{m-p-1} |x|^{-\gamma} \dx \dt}.
\end{split}
\end{equation}
where
\[
J:=2 c_3  \mu_\gamma(B_{R_1})^\frac{\sigma}{N-\gamma} \left[1+M_{\widetilde{p}}(u_0, \delta, T_0)^{1-m}\right] \left[\frac{ h_\sigma(R, R_1, x_0)}{\left(R-R_1 \right)^\sigma} + \frac{1}{T_1-T_0} \right]\ge 1.
\]
Note that it is not restrictive to assume that $J\ge 1$, by an argument  similar to the footnote related to formula \eqref{Correzione111}.

We estimate the  second term  in the right-hand side of \eqref{lower.first.inequality.1} using again \eqref{sup.energy.inequality.lower.approximate}, observing that we can always choose $a \in (1, \sr/2)$ such that $(m-p-1 )(a-1)q= -p$ to get
\begin{equation}\label{lower.first.inequality.3}
\begin{split}
\left[\sup_{t \in \left(T_1, T\right]}\int_{B_{R_1}}{u^{-p} |x|^{-\gamma} \dx}\right]^{\frac{1}{q}} &\leq \left[  c_3  \left(1+M_{\widetilde{p}}(u_0, \delta, T_0)^{1-m}\right) \right]^{\frac{1}{q}} \\
&  \times \left[\frac{h_\sigma(R, R_1, x_0)}{\left(R-R_1 \right)^\sigma} + \frac{1}{T_1-T_0} \right]^{\frac{1}{q}} \left[ \iint_{Q_0}{u^{-p+m-1}|x|^{-\gamma} \dx \dt} \right]^{\frac{1}{q}}.
\end{split}
\end{equation}
Combining inequalities \eqref{lower.first.inequality.1}, \eqref{lower.first.inequality.2} and \eqref{lower.first.inequality.3} we get \eqref{lower.first.inequality.4}.

\noindent$\bullet~$\textsc{Step 2. }\noindent\textit{The $k^{\rm th}$ iteration step. }We first define an increasing sequence of exponents $p_k$, recalling that $q= r^*/(r^* - 2)$ and $p_0>0$,  set
\begin{equation*}\label{exponents.pk.lower.iteration}
p_{k+1}:=\left(1+\frac{1}{q} \right)p_k =\left(1+\frac{1}{q} \right)^{k+1} p_0\xrightarrow[k\to+\infty]{}+\infty \,.
\end{equation*}
Next, we define the cylinders $Q_{k}:= \left(T_k, T \right] \times B_{R_k}\left( x_0 \right)$ as in \eqref{Def.Qk.upper.iteration}, so that
$Q\supset Q_{k}\supset Q_{k+1}\to Q_\infty$; we have chosen a decreasing sequence of radii $R_0>R_k>R_{k+1}\to R_\infty$  and an  increasing sequence of times $T_0 <  T_k < T_{k+1} \to T_{\infty} $ as in \eqref{radii.times.upper.iteration}, namely $R_{k}- R_{k+1}= C_1(R_{0}- R_{\infty})(k+1)^{-\alpha}$ and $T_{k+1}-T_{k} = C_2 (T_{\infty}-T_0)(k+1)^{-\alpha \sigma}$,
where $\alpha=2 \vee \sigma^{-1}$ and  $C_1, C_2>0$ are as in \eqref{pho_tau_definition_upper}.\\
Plugging all the above defined quantities in inequality \eqref{lower.first.inequality.4}, the $k^{\rm th}$ iteration step reads
\begin{align}\label{lower.first.inequality.6}
&\left[ \iint_{Q_{k+1}}{u^{-(p_{k+1}+1-m)}|x|^{-\gamma} \dx\dt} \right]^{-\frac{1}{p_{k+1}+1-m}}\\
&\geq \left[2 S_{\gamma, \beta}\right]^{-\frac{2}{p_{k+1}+1-m }}
\left[2c_3\,(1+M_{\widetilde{p}}(u_0, \delta, T_k)^{1-m})\right]^{-\frac{1}{p_{k+1}+1-m}\left(1+\frac{1}{q}\right)}  \nonumber \\
 & \times \left[\frac{h_{\sigma}( R_k, R_{k+1}, x_0 )}{(R_k-R_{k+1})^\sigma} + \frac{1}{T_{k+1}-T_k} \right]^{-\frac{1}{p_{k+1}+1-m}\left(1+\frac{1}{q}\right)}
\left[ \iint_{Q_k}{u^{-(p_k+1-m)}|x|^{-\gamma} \dx\dt} \right]^{-\frac{1}{p_{k+1}+1-m}\left(1+\frac{1}{q}\right)}\,.\nonumber
\end{align}

\noindent\textit{Bounds for the constants. }It is convenient to bound the constants appearing in \eqref{lower.first.inequality.6} by a quantity which does not depend on $p$, but only on $m, \gamma, \beta, R_{\infty}$ and $R_0$. Recall that $c_3:= 2K_{\psi} \ c_{m,p_k}$ is given in the energy inequality  \eqref{sup.energy.inequality.lower}: while $K_\psi>0$ depends only on $N$, the quantity
$c_{m, p_k} =    \frac{p_k+1}{p_k} \left( \frac{p_k+1}{p_k}\wedge \frac{2 m  (p_k+1)^2}{ (m-p_k-1)^2} \right)^{-1}$  needs to be bounded uniformly for all $k\ge 0$; since $ p_k > p_0>0$ it is easy to show that we have $ c_{m, p_k} \leq (1+p_0^{-1})/(1\wedge 4m)$,  so that  $2c_3=2c_3(p_k)\le \underline{c}=\underline{c}(m,N,p_0)<+\infty$\,; hence  $2c_3(p_k)\le \underline{c}$   can be bounded uniformly by a constant that depends only on $N, m,p_0$. As in \eqref{bounds.hsigma.upper}, also $ h_{\sigma}(R_k, R_{k+1}, x_0) $ can be bounded by a quantity depending only on $\sigma, R_0$ and $R_{\infty}$, as follows:
\begin{equation*}
 h_{\sigma}\left(R_k, R_{k+1}, x_0\right) \leq  h_{\sigma}\left( R_0, R_{\infty}, x_0 \right) (k+1)^{\alpha\left(2-\sigma\right)_{+}} \, C_1^{-\left(2-\sigma\right)_{+}},
\end{equation*}
where $C_1>0$ is as in  \eqref{pho_tau_definition_upper}. Moreover, by \eqref{uniform.boundedness.approximate.solutions} we have that $M_{\widetilde{p}}(u_0, \delta, T_k)\le M_{\widetilde{p}}(u_0, \delta, T_0)$, since $T_0<T_k$.
Finally, we can rewrite the $k^{\rm th}$ iteration step \eqref{lower.first.inequality.6} as follows: \vspace{-2mm}
\begin{equation}\label{lower_iteration_first_step_7}
\begin{split}
\biggl[ \iint_{Q_{k+1}}{u^{-(p_{k+1}+1-m)}\frac{ \dx\dt}{|x|^{\gamma}}} \biggr]^{\frac{-1}{p_{k+1}+1-m}}\geq I_{k+1}^{\frac{-1}{p_{k+1}+1-m}}
\left[ \iint_{Q_k}{u^{-(p_k+1-m)} \frac{\dx\dt}{|x|^{\gamma}}} \right]^{\frac{-\left(1+\frac{1}{q}\right)}{p_{k+1}+1-m}}\,,\vspace{-2mm}
\end{split}
\end{equation}
where the constant $ I_{k+1} $ is bounded by\vspace{-1mm}
\begin{equation}\label{lower.first.inequality.7}\begin{split}
&I_{k+1}
 \le  \left[ 2 S_{\gamma, \beta}\right]^2  \left[\underline{c}\left(\frac{h_{\sigma}\left( R_0, R_{\infty}, x_0 \right) C_1^{-\frac{b}{\alpha}}}{(R_0-R_\infty)^\sigma} + \frac{C_2^{-1}}{T_\infty-T_0} \right)  (1+M_{\widetilde{p}}(u_0, \delta, T_0)^{1-m})\right]^{1+\frac{1}{q}} (k+1)^{b\left(1+\frac{1}{q}\right)}  \\
 &\leq \left[ 2 S_{\gamma, \beta}\right]^2 \left[\underline{c}(C_1^{-\frac{b}{\alpha}}\vee C_2^{-1}) \left(\frac{h_{\sigma}\left( R_0, R_{\infty}, x_0 \right) }{(R_0-R_\infty)^\sigma} + \frac{1}{T_\infty-T_0} \right)  (1+M_{\widetilde{p}}(u_0, \delta, T_0)^{1-m})\right]^{1+\frac{1}{q}} (k+1)^{b\left(1+\frac{1}{q}\right)} \\
&:= J_0 J_1^{1+\frac{1}{q}} (k+1)^{b\left(1+\frac{1}{q}\right)}\,,
\end{split}\end{equation}
where $b=2\alpha$ if $0 < \sigma < 2$, $b=\sigma \alpha$ otherwise.  The assumption \eqref{pho_tau_definition_upper} on $\alpha$ ensures that $b> 1$.

\noindent$\bullet~$\textsc{Step 3. }\noindent\textit{The iteration. }
We now iterate inequalities \eqref{lower_iteration_first_step_7} to get
\begin{equation}\label{lower.first.inequality.8}
\biggl[ \iint_{Q_{k+1}}{u^{-(p_{k+1}+1-m)}\frac{ \dx\dt}{|x|^{\gamma}}} \biggr]^{\frac{-1}{p_{k+1}+1-m}}
 \geq  \prod_{j=1}^{k+1}  I_{j}^{\frac{-\left(1+\frac{1}{q} \right)^{k+1-j}}{p_{k+1}+1-m}}
 \biggl[ \iint_{Q_0}{u^{-(p_0+1-m)}\frac{ \dx\dt}{|x|^{\gamma}}} \biggr]^{\frac{-\left(1+\frac{1}{q} \right)^{k+1}}{p_{k+1}+1-m}}
\end{equation}
We  estimate  the product appearing in the above inequality as follows:
\begin{align*}\label{constant.inf.limit}
 \prod_{j=1}^{k+1}{I_{j}^{\frac{-1}{p_{k+1}+1-m}\left(1+\frac{1}{q} \right)^{k+1-j}}}& \geq \left[J_0 J_1^{1+\frac{1}{q}} \right]^{\frac{-1}{p_{k+1}+1-m}\sum\limits_{j=0}^{k}{\left(1+\frac{1}{q}\right)^{j}}}
    \left[\left(k+1 \right)^{b'}\right]^{\frac{-1}{p_{k+1}+1-m}}
    \cdots \left[2^{b'}\right]^{\frac{-1}{p_{k+1}+1-m}\left(1+\frac{1}{q}\right)^{k-1}}  \nonumber \\
 & = \left[J_0 J_1^{1+\frac{1}{q}} \right]^{\frac{-1}{p_{k+1}+1-m}\sum\limits_{j=0}^{k}{\left(1+\frac{1}{q}\right)^{j}}} \  \prod_{j=1}^{k+1}{j^{\frac{-b'}{p_{k+1}+1-m}\left(1+\frac{1}{q} \right)^{k+1-j}}}\,,
\end{align*}
 where $b'=b\left(1+\frac{1}{q}\right)$.  Recalling that $p_k=(1+\frac{1}{q})^k p_0$ we get
\begin{equation*}
  \frac{-1}{p_{k+1}+1-m}\sum_{j=0}^{k}{\left(1+\frac{1}{q}\right)^{j}}\xrightarrow[k\to+\infty]{}  -\frac{q}{p_0} \qquad \mbox{and} \qquad \frac{p_{k+1}+1-m}{\left(1+\frac{1}{q} \right)^{k+1}} \xrightarrow[k\to+\infty]{} p_0\,.
\end{equation*}
Moreover, it is easy to show that
\begin{equation*}
\lim_{k\rightarrow +\infty}\prod_{j=1}^{k+1}j^{\frac{b'}{p_0(1+\frac{1}{q})^{k+1}+1-m}\left(1+\frac{1}{q}\right)^{k+1-j}} \le \lim_{k\rightarrow +\infty}\exp\left(\frac{b'}{p_0 + \frac{(1-m)}{(1+\frac{1}{q})^{k+1}}}\sum_{j=1}^{\infty}\left(\frac{q}{q+1}\right)^j \log{j}\right)
\le (c'')^{\frac{q}{p_0}}<+\infty\,.
\end{equation*}
Taking the limit in \eqref{lower.first.inequality.8} as $k \rightarrow \infty$ we obtain
\begin{align}\label{lower.first.inequality.9}
 \inf_{Q_{\infty}}{u} \geq   \left[  c''\,J_0 J_1^{1+\frac{1}{q}} \right]^{-\frac{q}{p_0}} \left(\iint_{Q_0}{u^{-\left(p_{0}+1-m \right)} |x|^{-\gamma} \dx \dt} \right)^{-\frac{1}{p_{0}}}.
\end{align}
Note that $ J_0 $ and $   c'' $ depend only on $ m $, $ p_0 $, $N$, $ \gamma $  and $ \beta $, while $ J_1 $ depends also on  $ R_0 - R_{\infty}$ and $ T_{\infty}-T_{0} $.
\noindent$\bullet~$\textsc{Step 4. }The goal of this step is show that estimate \eqref{lower.first.inequality.9} implies estimate \eqref{formula.inf}. To this end, we use the time monotonicity properties of the solution, namely  the fact  that $t \rightarrow u(t,x)t^{-\frac{1}{1-m}}$ is non-increasing in time for almost every $x \in B_{R_0}(x_0)$. Recalling that $p_0 > 0$, we get for any $t \in [T_0, T]$
\begin{equation*}
  \int_{B_{R_0}(x_0)}  u(t,x)^{m-1-p_0} |x|^{-\gamma}  \dx \leq \left(\frac{T}{T_0}\right)^{\frac{p_0+1-m}{1-m}}\int_{B_{R_0}(x_0)}u(T,x)^{m-p_0-1} |x|^{-\gamma} \dx.
\end{equation*}
Hence, \eqref{lower.first.inequality.9} can be estimated as follows (recall that $J_0, J_1$ are defined in \eqref{lower.first.inequality.7})
 \begin{equation*}\begin{split}
 \inf_{Q_{\infty}}{u}
 & \ge \kb \left[ J_0 J_1^{1+\frac{1}{q}} \right]^{-\frac{q}{p_0}}\left(\frac{T_0}{T}\right)^{\frac{p_0+1-m}{p_0(1-m)}}\left(T-T_0\right)^{-\frac{1}{p_0}}\left(\int_{B_{R_0}(x_0)}u^{m-p_0-1}(T, x) |x|^{-\gamma} \dx\right)^{-\frac{1}{p_0}}
\end{split}\end{equation*}
This is exactly \eqref{formula.inf}, with $R_0 \rightsquigarrow R_1$, $R_{\infty}\rightsquigarrow R_2$,  $T_{\infty}\rightsquigarrow T_1$,  $T \rightsquigarrow T_2$ and $s:=p_0+1-m$, $p:=\widetilde{p}$.\qed
\noindent\textbf{Remark. }Note that the estimate degenerates in the limits $m \to 1^-$ or $m\to 0^+$, indeed the term $(T_0/T)^{\frac{p_0+1-m}{p_0(1-m)}}\xrightarrow[]{m \to 1^-} 0$; also,  by \eqref{lower.first.inequality.7} and previous discussions  we have $J_1 \sim 1/m $, so that $J_1^{-\frac{q+1}{p_0}}\xrightarrow[]{m \to 0^+} 0$.

\subsection{Step 2. Smoothing effects for negative norms.}\label{ssec.Step2}
In the previous step we have proved an estimate of type $\LL^{-\infty} - \LL^{-s}$, which holds for any $s > 1-m$, but this is not sufficient in order to use the reverse H\"older inequalities of Corollary \ref{Reverse.Holder.Corollay}, which may  hold  only for exponents close to $0$. In this step we solve this issue by proving $\LL^{-s}- \LL^{-\varepsilon}$ estimates, for any $\varepsilon \in (0,1-m)$, through a finite iteration.

\begin{prop}[$\LL^{-s}- \LL^{-\varepsilon}$ Smoothing Effects]\label{prop.finite.iter.negative.norms}
Let $u_\delta$ be a strong (super)solution to \ref{LARGE.DIRICHLET.PROBLEM} on $(0, \infty) \times B_{R_0}(x_0)$, corresponding to $u_0\chi_{B_R(x_0)} \in \LL^p_\gamma(B_{R_0}(x_0))$ with $ p > p_c $ when $m\in (0,m_c]$ and with $p \geq 1  $ when $m\in (m_c,1)$. Let $0<R_2 < R_1 \leq R$  be  such that $0\not\in \overline{B_{R_1}(x_0)\setminus B_{R_2}(x_0)}$ and let $0 < T_0 < T_1$. For any $\varepsilon\in (0,1-m)$  there exists $s_\varepsilon\in \big(1-m,\frac{r^*}{2}(1-m)\big]$ and $\kb_4>0$  such that
\begin{align}\label{ineq.finite.iter.negative.norms}
&\left[ \int_{B_{R_2}(x_0)}u_\delta(T_0, x)^{-s_\varepsilon} |x|^{-\gamma} \dx  \right]^{-\frac{1}{s_\varepsilon}}
\geq \kb_4  \left(\frac{T_0}{T_1}\right)^{\frac{1}{1-m}}\left[ \int_{B_{R_1}(x_0)}u_\delta(T_1, x)^{-\varepsilon}|x|^{-\gamma} \dx \right]^{-\frac{1}{\varepsilon}}\,    \nonumber\\
&\times\bigg[  (1 + M_p(u_0, \delta, T_0)^{1-m})
 \left(\mu_\gamma(B_{R_2}(x_0))^{-\frac{\sigma}{N-\gamma}}\vee\left(\frac{h_\sigma(R_1, R_2, x_0)}{(R_1-R_2)^\sigma}+ \frac{1}{T_1-T_0}\right)\right)\bigg]^{\zeta_\varepsilon}\,,
\end{align}
where $\ke$ is the smallest $k\in \NN$ such that $(r^*/2)^{k}\varepsilon>1-m$ and
\begin{align}\label{const.finite.iter.negative.norms}
 s_\varepsilon& :=(r^*/2)^{\ke}\varepsilon>1-m \qquad\mbox{and}\qquad
\zeta_\varepsilon:={-\frac{1}{\varepsilon}\frac{1-\left(\frac{2}{r^*}\right)^{\ke}}{1-\frac{2}{r^*}}}\nonumber \\
\kb_4&:= \kb_4' \left\{
\begin{array}{lll}
1 &\mbox{if $\varepsilon\neq (1-m)(2/r^*)^k$ for all $k\in \NN$}\\
\mu_\gamma(B_{R_1})^{-\frac{1}{1-m}\big(\frac{r^*}{2}\big)^k\frac{r^*-2}{r^*+2}}  &\mbox{if $\varepsilon= (1-m)(2/r^*)^k$ for some $k\in \NN$}\,,\\
\end{array}
\right.
\end{align}
where $h_\sigma$, $M_p$ and $\sr$ are defined in \eqref{function_h}, \eqref{uniform.boundedness.approximate.solutions} and \eqref{w_sobolev} respectively, while $\kb_4'$ only depends on $N,m,\beta,\gamma,\varepsilon$.
\end{prop}
\noindent {\bf Proof.~}The proof relies on a finite iteration and it is split into two steps. Let us fix $x_0\in \RR^N$\,, and simply denote $B_R=B_R(x_0)$ and $u=u_\delta\ge \delta>0$, when there is no ambiguity.

\noindent$\bullet~$\textsc{Step 1. }\textit{Preparation of the iteration step. }We are going to prove that for any $s \in (0, 1-m) $, for any  $0 < R _2 < R_1 \le R \le R_0$  and $0 < t_0 < t_1$, there exists a constants $\kb_4>0$, depending on $q,  \gamma, \beta, N, m, x_0 $ and $R_1, R_2, t_0, t_1$ and the norm of the initial data $\|u_0 \|_{\LL^p_\gamma (B_R(x_0))}$ such that
\begin{equation}\label{inequality.iteration.again}
\left[ \int_{B_{R_2}(x_0)}u(t_0, x)^{\frac{-\sr s}{2}} |x|^{-\gamma} \dx  \right]^{-\frac{2}{\sr s}}
\geq \underline{c}     \left[ \int_{B_{R_1}(x_0)}u(t_1, x)^{-s}|x|^{-\gamma} \dx \right]^{-\frac{1}{s}},
\end{equation}
where $\underline{c} =\left[\left( \mu_\gamma(B_{R_2})^{\frac{-\sigma}{(N-\gamma)}}\vee \overline{c}_1 \right)4S_{\gamma, \beta}^2  \right]^{-1/s} (t_1/t_0)^{1/(1-m)}$, $\overline{c}_1$ is as in \eqref{cbar} and $S_{\gamma, \beta}$ is as in \eqref{WSIcomplete}. \\
 Note  that while $s\in(0,1-m)$\,, in general $r^*s/2$ may be bigger than  $1-m$:  we will exploit this fact in the next step. We first integrate in time inequality \eqref{WSIcomplete} applied to  $f=u^{\frac{\tilde{p}+m-1}{2}}$,  with $\tilde{p}\in (0,1-m)$, to get
\[
\int\limits_{t_0}^{t}\Bigg(\int\limits_{B_{R_2}} u^{\frac{\sr}{2}(\tilde{p}+m-1 )}  \frac{\dx}{|x|^{\gamma}}\Bigg)^{\frac{2}{\sr}}\dtau
\leq  2 \overline{S}_{\gamma, \beta}^2  \Bigg[  \mu_\gamma(B_{R_2})^{\frac{-\sigma}{(N-\gamma)}}\int\limits_{t_0}^{t}\int\limits_{B_{R_2}} u^{\tilde{p}+m-1} \frac{\dx\dtau}{|x|^{\gamma}} + \int\limits_{t_0}^{t}\int\limits_{B_{R_2}} \left|\nabla u^{\frac{\tilde{p}+m-1}{2}}\right|^2  \frac{\dx\dtau}{|x|^{\beta}} \Bigg]
\]
Recalling that by \eqref{uniform.boundedness.approximate.solutions} we have $u(t,x)\le M_p(u_0, \delta, t_0)$, for all $t\ge t_0$ and $x\in B_R(x_0)$, then the energy inequality \eqref{energy.inequality.last} implies for all $\tilde{p}\in (0,1-m)$ (recall that $\tilde{p}+m-1<0$ and $u\ge \delta>0$)
\begin{equation*}\label{modified.energy.inequality.last}\begin{split}
   \int_{t_0}^{t}\int_{B_{R_2}} \left|\nabla u^{\frac{\tilde{p}+m-1}{2}} \right|^2  \frac{\dx \dt}{|x|^{\beta} }
 \leq \overline{c} \int_{t_0}^{t_1}\int_{B_{R_1}}u^{\tilde{p}+m-1}  \frac{\dx \dt}{|x|^{\gamma} }.
\end{split}\end{equation*}
where $t \in (t_0, t_1)$ and
\begin{equation*}\label{cbar0}
\overline{c}=c_2 \left[1+M_p(u_0, \delta, t_0)^{1-m}\right]\left[\frac{ h_\sigma(R_1, R_2, x_0)}{(R_1-R_2)^\sigma}+ \frac{1}{t_1-t}\right],
\end{equation*}
$h_\sigma$ and $M_p$ are defined in \eqref{function_h} and \eqref{uniform.boundedness.approximate.solutions} and  $c_2=c_2(m,\tilde{p})>0$ is as in \eqref{energy.inequality.last}.
Combining the two above inequalities we obtain
\begin{equation}\label{smoothing.negative.spacetime}\begin{split}
\int_{t_0}^{t}\Bigg(\int_{B_{R_2}} u^{\frac{\sr(\tilde{p}+m-1)}{2}}  \frac{\dx}{|x|^{\gamma}}\Bigg)^{\frac{2}{\sr}}\dtau
& \leq 2\overline{S}_{\gamma, \beta}^2 \left(\mu_\gamma(B_{R_2})^{\frac{-\sigma}{(N-\gamma)}}\vee \overline{c} \right)
 \int_{t_0}^{t_1}\int_{B_{R_1}(x_0)}u^{\tilde{p}+m-1}   \frac{\dx\dtau}{|x|^{\gamma}}
\end{split}\end{equation}
We recall that $t \rightarrow u(t,x)t^{-\frac{1}{1-m}}$ is non-increasing in time for almost every $x \in B_{R_0}(x_0)$, hence we can estimate the two sides of the above inequality: the left-hand side can be estimated from below
\begin{equation}\label{smoothing.negative.spacetime.lll}\begin{split}
  \int_{t_0}^{t}\left[\int_{B_{R_2}(x_0)}u(\tau, x)^{\frac{\sr(\tilde{p}+m-1)}{2}} \frac{\dx }{|x|^{\gamma}}  \right]^{\frac{2}{\sr}} \dtau &\geq \tilde c\, \frac{t^{\frac{\tilde{p} }{1-m}}-t_0^{\frac{\tilde{p} }{1-m}}}{t_0^{\frac{ \tilde{p}+m-1}{1-m}}}
   \left[ \int_{B_{R_2}(x_0)}u(t_0, x)^{\frac{\sr(\tilde{p}+m-1)}{2}} \frac{\dx}{|x|^{\gamma}} \right]^{\frac{2}{\sr}},
\end{split}\end{equation}
where $\tilde c=(1-m)/\tilde{p} $. Analogously, we can estimate the right-hand side of \eqref{smoothing.negative.spacetime} from above,
\begin{equation}\label{smoothing.negative.spacetime.rrr}
\int_{t_0}^{t_1}\int_{B_{R_1}(x_0)}u(t,x)^{\tilde{p}+m-1} \frac{\dx \dt }{|x|^{\gamma}} \leq \tilde c  \frac{t_1^{\frac{\tilde{p} }{1-m}}-t_0^{\frac{\tilde{p} }{1-m}}}{t_1^{\frac{ \tilde{p}+m-1}{1-m}} }  \int_{B_{R_1}(x_0)}u(t_1, x)^{\tilde{p}+m-1}\frac{\dx}{|x|^{\gamma}} .
\end{equation}
Finally, letting $s=  -\tilde{p}+1-m$, we have $s \in (0, 1-m)$; taking $t^{\tilde{p}/(1-m)}=\bigl(t_1^{\tilde{p} /(1-m)}+t_0^{\tilde{p} /(1-m)}\bigr)/2$, we have $2=\bigl(t_1^{\tilde{p} /(1-m)}-t_0^{\tilde{p} /(1-m)}\bigr)/\bigl(t^{\tilde{p} /(1-m)}-t_0^{\tilde{p} /(1-m)}\bigr)$, so that \eqref{smoothing.negative.spacetime}, \eqref{smoothing.negative.spacetime.lll} and \eqref{smoothing.negative.spacetime.rrr}  give  \eqref{inequality.iteration.again}.  Note  that we also have
\begin{equation}\label{cbar}
\overline{c}\le 2c_2   \left[1+M_p(u_0, \delta, t_0)^{1-m}\right]
\left[\frac{h_\sigma(R_1, R_2, x_0)}{(R_1-R_2)^\sigma}+ \frac{1}{t_1-t_0}\right]:= \overline{c}_1,
\end{equation}
since we have that $\tilde{p} <1-m$ implies $t_1-t= t_1-\big[\bigl(t_1^{\tilde{p} /(1-m)}+t_0^{\tilde{p} /(1-m)}\bigr)/2\big]^{(1-m)/\tilde{p}}\ge (t_1-t_0)/2$.

\noindent$\bullet~$\textsc{Step 2. }\textit{The finite iteration. }Fix $\varepsilon\in (0,1-m)$ and assume that $\varepsilon\neq (1-m)(2/r^*)^k$ for all $k\in \NN$, to avoid that $(r^*/2)^{k}\varepsilon=1-m$ for some $k\in \NN$; the remaining cases are similar and will be discussed at the end of the proof. Let $\ke$ be the smallest positive integer such that  $(r^*/2)^{\ke}\varepsilon>1-m$: note  that $\ke>\log[(1-m)/\varepsilon]/\log[r^*/2]$. We are going to iterate inequality \eqref{inequality.iteration.again} $\ke$ times; let us define a decreasing sequence of exponents, and increasing sequences of radii and times for all $0\le i\le \ke$
\begin{equation}\label{finite.iter.s.r.t}
s_i:=(r^*/2)^{\ke-i}\varepsilon\,,\qquad r_i:=R_2+\frac{i}{\ke}(R_1-R_2)\qquad\mbox{and}\qquad t_i:= T_0+\frac{i}{\ke}(T_1-T_0)
\end{equation}
 Note  that by construction $s_i\in (0,1-m)$ for all  $i=0,\dots,\ke$ , $s_{\ke}=\varepsilon$ and $s_0=(r^*/2)^{\ke}\varepsilon>1-m$, so that we can rewrite inequality \eqref{inequality.iteration.again} as follows:
\begin{equation}\label{inequality.iteration.again2}
 \|u (t_i)\|_{\LL^{-s_i}_\gamma(B_{r_i}(x_0))} \geq c_i  \|u (t_{i+1})\|_{\LL^{-s_{i+1}}_{\gamma}(B_{r_{i+1}}(x_0))},
\end{equation}
where $c_i=\underline{c}(s_i,t_i,r_i)$ has the expression
\begin{equation}\label{formula.ci}\begin{split}
  c_i = &  \left[4 S_{\gamma, \beta}^2 \left(1\vee  c(m,1-m-s_{i+1}) \left(1+M_p(u_0, \delta, t_{i})^{1-m}\right)
  \right)\right]^{-\frac{1}{s_{i+1}}}\\
& \times \left[\mu_\gamma(B_{r_i}(x_0))^{-\frac{\sigma}{N-\gamma}}\vee\left(\frac{h_\sigma(r_{i+1}, r_{i}, x_0)}{(r_{i+1}-r_i)^\sigma}+ \frac{1}{t_{i+1}-t_{i}}\right)  \right]^{-\frac{1}{s_{i+1}}}   \left[\frac{t_i}{t_{i+1}}\right]^{\frac{1}{1-m}},
\end{split}
\end{equation}
where $c=c(m, 1-m-s_{i+1})$ is as in \eqref{energy.inequality.last}. Iterating $\ke$-times inequality \eqref{inequality.iteration.again2} we get
\begin{equation}\label{finte.iteration}\begin{split}
 \|u(t_0)\|_{\LL^{-s_0}_\gamma (B_{r_0})} &\ge c_0 \|u(t_1)\|_{\LL^{s_1}(B_{r_1})}
 \ge\dots\ge \left(\prod_{i=0}^{\ke-1}c_i\right)  \|u(t_{\ke})\|_{\LL^{-s_{\ke}}_\gamma (B_{r_{\ke}})}.
\end{split}\end{equation}
Recalling that $t_0=T_0<t_{\ke}=T_1$, $r_0=R_2<r_{\ke}=R_1\le R$\,, and that $s_0=(r^*/2)^{\ke}\varepsilon>1-m$ we have \begin{equation*}\label{finte.iteration.2}\begin{split}
 \|u(T_0)\|_{\LL^{-s_0}_\gamma (B_{R_2})}
 \ge C  \|u(T_1)\|_{\LL^{-\varepsilon}_\gamma (B_{R_1})},
\end{split}\end{equation*}
where $C>0$ is the lower bound of the product $\prod_{i=0}^{\ke-1}c_i$ that we are going to estimate explicitly below.
From formulae \eqref{formula.ci}  and \eqref{finite.iter.s.r.t} we deduce that
\begin{equation*}\label{bound.c1}\begin{split}
  c_i &\ge \underline{c}_i  \left(\frac{t_i}{t_{i+1}}\right)^{\frac{1}{1-m}}
  \left[\mu_\gamma(B_{R_2}(x_0))^{-\frac{\sigma}{N-\gamma}}\vee\left(\frac{h_\sigma(R_1, R_2, x_0)}{(R_1-R_2)^\sigma}+ \frac{1}{T_1-T_0}\right)\right]^{-\frac{1}{s_{i+1}}}    \left[\,\big(1 + M_p(u_0, \delta, T_0)^{1-m}\big)\right]^{-\frac{1}{s_{i+1}}}
\end{split}\end{equation*}
 since $M_p(u_0, \delta, t_{i})\leq M_p(u_0, \delta, T_0)$;  moreover,  $h_\sigma(r_{i+1}, r_i, x_0) \leq \ke^{(2-\sigma)_{+}} h_\sigma(R_1, R_2, x_0) $ and $\mu_\gamma(B_{R_2})\leq \mu_\gamma(B_{r_i})$;   finally we have set  $\underline{c}_i=\left[2^4 S_{\gamma, \beta}^2 c(m,1-m-s_{i+1})(\ke^{2\vee \sigma}+\ke)\right]^{\frac{-1}{s_{i+1}}} $.
Finally, we can estimate $C$ as  follows:
\begin{align*}
 &\prod_{i=0}^{\ke-1}c_i  \ge \left(\frac{T_0}{T_1}\right)^{\frac{1}{1-m}} \left[\prod_{i=0}^{\ke-1}\underline{c}_i\right]
 \left[  (1 + M_p(u_0, \delta, T_0)^{1-m})
 \left(\underline{\mu}\vee\left(\frac{h_\sigma(R_1, R_2, x_0)}{(R_1-R_2)^\sigma}+ \frac{1}{T_1-T_0}\right)\right)\right]^{\sum\limits_{i=0}^{\ke-1}\frac{-1}{s_{i+1}}} \nonumber\\
&\ge \kb_4' \left(\frac{T_0}{T_1}\right)^{\frac{1}{1-m}}
 \left[  (1 + M_p(u_0, \delta, T_0)^{1-m})
 \left(\underline{\mu}\vee\left(\frac{h_\sigma(R_1, R_2, x_0)}{(R_1-R_2)^\sigma}+ \frac{1}{T_1-T_0}\right)\right)\right]^{-\frac{1}{\varepsilon}\frac{1-\left(\frac{2}{r^*}\right)^{\ke}}{1-\frac{2}{r^*}}}\,,
 \end{align*}
where we have put $\underline{\mu}:= \mu_\gamma(B_{R_2}(x_0))^{-\frac{\sigma}{N-\gamma}}$ and   we have used that  $\sum_{i=0}^{\ke-1}\frac{1}{s_{i+1}}=\frac{1}{\varepsilon}\sum_{j=0}^{\ke-1}\left(\frac{2}{r^*}\right)^j
=\frac{1}{\varepsilon}\big[1-\left(\frac{2}{r^*}\right)^{\ke}\big]/\left[ 1-\frac{2}{r^*} \right]$, and that $\prod_{i=0}^{\ke-1} \left(t_i/t_{i+1}\right)^{1/(1-m)}=\left(T_0/T_1\right)^{1/(1-m)}$ and we have defined $\kb_4'=\prod_{i=0}^{\ke-1}\underline{c}_i>0$ so that it only depends on $N,m,\beta,\gamma,\varepsilon$.

\noindent$\bullet$ \textit{The cases when $\varepsilon=(1-m)(2/r^*)^k$ for some $k\in \NN$. }When $\varepsilon=(1-m)(2/r^*)^k\in (0,1-m)$, one can start by a slightly smaller value, say $\tilde\varepsilon=\frac{1-m}{2}\big[\big(\frac{2}{r^*}\big)^k+\big(\frac{2}{r^*}\big)^{k+1}\big]<\varepsilon$, proceed as above and obtain  \eqref{finte.iteration} with $k_{\tilde\varepsilon}=k+1$, namely $\|u(T_0)\|_{\LL^{-s_0}_\gamma (B_{R_2})} \ge C  \|u(T_1)\|_{\LL^{-\tilde\varepsilon}_\gamma (B_{R_1})}$ with $s_0=(r^*/2)^{k+1}\tilde\varepsilon>1-m$\,, and then conclude by  H\"older's inequality, observing that
$\|u(T_1)\|_{\LL^{-\tilde\varepsilon}_\gamma (B_{R_1})}\ge \mu_\gamma(B_{R_1})^{\frac{1}{\varepsilon}-\frac{1}{\tilde\varepsilon}} \|u(T_1)\|_{\LL^{-\varepsilon}_\gamma (B_{R_1})}$\,, since $-\tilde\varepsilon> -\varepsilon$; finally we notice that $\frac{1}{\varepsilon}-\frac{1}{\tilde\varepsilon}=-\frac{1}{1-m}\big(\frac{r^*}{2}\big)^k\frac{r^*-2}{r^*+2}$.\qed

\subsection{Step 3. Reverse H\"older inequalities}\label{ssec.parab.JN}
In this subsection we prove the Step 3 of the proof of our positivity result. More precisely we prove $\LL^{-\varepsilon}\to\LL^{\varepsilon}$ estimates, in the form of reverse H\"older inequalities.
\begin{prop}[Reverse H\"older inequality for \ref{LARGE.DIRICHLET.PROBLEM}]\label{Lem.Step3}
Let $u_\delta$ be a solution to \ref{LARGE.DIRICHLET.PROBLEM} on $(0, \infty) \times B_{R_0}(x_0)$.  Let $u_0\chi_{B_R(x_0)} \in \LL^p_\gamma(B_{R_0}(x_0))$ with $ p > p_c $ when $m\in (0,m_c]$ and with $p \geq 1  $ when $m\in (m_c,1)$. Then for all   $t\ge t_0 > 0$  and all $0<R_1<R$, there exists $\nu_\delta=\nu_\delta(t_0,u_0)>0$  such that
\begin{equation}\label{Reverse.Holder.nu.delta}
\|u_\delta(t)\|_{\LL^{\varepsilon}_\gamma(B_{R_1}(x_0))} \le \ka_{7}^{2/ \varepsilon} \mu_\gamma(B_{R_1}(x_0))^{2/  \varepsilon } \|u_\delta(t)\|_{\LL^{-\varepsilon}_\gamma(B_{R_1}(x_0))}\qquad\mbox{for all $0 < \varepsilon <   \nu_\delta$,}
\end{equation}
where
\begin{equation*}\label{nu.delta.def}
 \nu_\delta :=  \frac{1}{\ka_{14} \, \ka_{6}} \left[ 1 + \frac{R^\sigma}{t_0}\left(\frac{|x_0|}{R}\vee 1\right)^{\beta-\gamma} M_p(u_0, \delta, t_0)^{1-m}\right]^{-1/2},
\end{equation*}
$M_p$ is given in \eqref{uniform.boundedness.approximate.solutions}, $\ka_{14}$ is as in Lemma \ref{Caccioppoli.JohnNirenberg.lemma.LDP}, and  $\ka_{6}, \ka_{7}$ are as in Corollary \ref{Reverse.Holder.Corollay}.
\end{prop}
\noindent\textbf{Remark. }As it happens in the elliptic case, the above reverse H\"older inequality plays a fundamental role in the proof of the lower bounds; the above Proposition can be considered the parabolic analog of the celebrated John-Nirenberg Lemma, cf \cite{FG-PAMS}. As far as we know, in the literature of parabolic equations, there are basically only two techniques that allow  one  to prove the above estimates: one is due to Moser \cite{Moser}, the other is due to Bombieri and Giusti  \cite{BGInv72}, see also \cite{MoserCpam71}. None of the  previous  techniques applies directly to our nonlinear setting: in order to ensure the validity of the reverse H\"older inequalities of  Corollary \ref{Reverse.Holder.Corollay},  we need to show that $\log u\in BMO_\gamma$ ($u$ is a solution to the MDP or to the \ref{LARGE.DIRICHLET.PROBLEM}). In order to obtain a quantitative control on the $BMO_\gamma$ norm, we will use the Caccioppoli inequalities \eqref{caccioppoli.inequality} of Lemma \ref{Lem.Caccioppoli}, combined with the weighted Poincar\'e inequality of Proposition \ref{zero.mean.Poincare} as follows.

\begin{lem} \label{Caccioppoli.JohnNirenberg.lemma.LDP}
Let $u_\delta$ be a non-negative solution to \ref{LARGE.DIRICHLET.PROBLEM} on $(0,\infty)\times\Omega$ and let $B_{2R}(x_0)\subset \Omega$. Let $u_0 \in \LL^p_{\gamma}(\Omega)$ for $p > p_c$ if $0 < m \leq m_c$ or $p \geq 1$ if $m_c < m < 1$. Then for any $t> 0 $ the function  $\log u_\delta (t) \in BMO_\gamma(B_R(x_0))$,  more precisely there exists a constant  $\ka_{14} = \ka_{14}(N, m \gamma, \beta) >0$  such that for any $t > 0$
\begin{equation}\label{BMO.norm.delta}
\| \log u_\delta(t)\|_{BMO_{\gamma}(B_R(x_0))} \le \ka_{14}\left[ 1 + \frac{R^\sigma}{t}\left(\frac{|x_0|}{R}\vee 1\right)^{\beta-\gamma} M_p(u_0, \delta, t)^{1-m}     \right]^{\frac{1}{2}}= \frac{1}{\ka_6 \nu_\delta},
\end{equation}
where $M_p$ is given in \eqref{uniform.boundedness.approximate.solutions} and $\ka_6$ is as in Corollary \ref{Reverse.Holder.Corollay}.
\end{lem}
\noindent\textbf{Proof of Lemma \ref{Caccioppoli.JohnNirenberg.lemma.LDP}. }We will write $u=u_\delta$, since no confusion arises here. Let $ R > \rho > 0$, $h >0 $ and   $\psi \in C^{\infty}_{c}(B_ {2\rho}(x_0))$ . Then by  Cacciopoli's inequality \eqref{caccioppoli.inequality} of  Lemma \ref{Lem.Caccioppoli}, we get
 \begin{equation*}
 \begin{split}
   &\frac{m^2(1-m)}{2}\int_{B_{2\rho}(x_0)}\frac{1}{h}\int_{\tau}^{\tau+h}\psi^2|\nabla \log u|^2 \dt |x|^{-\beta} \dx \\
   &  \leq 2(1-m) \int_{B_{2\rho}(x_0)}|\nabla \psi|^2 |x|^{-\beta} \dx
   + \int_{B_{2\rho}(x_0)}\frac{u^{1-m}(\tau+h, x) - u^{1-m}(\tau, x) }{h} \psi^2 |x|^{-\gamma} \dx.
\end{split}
\end{equation*}
By Lebesgue's Differentiation Theorem, the Steklov averages $\frac{1}{h}\int_{s}^{s+h}\psi^2|\nabla \log u|^2 \dt$ converge,  for almost every $t$, to $\psi^2|\nabla \log u|^2$ as $h \rightarrow 0$. Using the time monotonicity property of $u$, namely that
$u(\tau+h,x)^{1-m} \le u(\tau,x)^{1-m}\left(\frac{\tau+h}{\tau}\right)$, we get
 \begin{equation*}
 \int_{B_{2\rho}(x_0)}\frac{u^{1-m}(\tau+h, x) - u^{1-m}(\tau, x) }{h} \psi^2 |x|^{-\gamma} \dx \leq \frac{1}{\tau}\int_{B_{2\rho}(x_0)}u^{1-m}(\tau,x) \psi^2 |x|^{-\gamma} \dx.
 \end{equation*}
 Now we can take  $\psi = 1$ on $B_{\rho}(y)$ and $\psi= 0 $ outside $B_{2 \rho}(x_0)$,  such that $|\nabla \psi|\le c_N\rho^{-1}$  and let $v=\log u$,  letting $h \rightarrow 0$ we obtain
 \begin{equation}\label{upperbound.gradient.log}
   \begin{split}
        \int_{B_{2\rho(x_0)}}\psi^2|\nabla v(\tau)|^2 \frac{ \dx }{|x|^{\beta}}&\leq \frac{4}{m^2} \int_{B_{2\rho(x_0)}}|\nabla \psi|^2 \frac{ \dx }{|x|^{\beta}}  + \frac{2}{m^2(1-m)\tau}\int_{B_{2\rho(x_0)}}u^{1-m}(\tau) \psi^2 \frac{ \dx }{|x|^{\gamma}} \\
   & \leq\frac{4 c_N\mu_\beta (B_{2 \rho}(x_0))}{ \rho^2 m^2} + \frac{2 \mu_\gamma(B_{2 \rho}(x_0))}{\tau(1-m)m^2} M_p(u_0, \delta, \tau)^{1-m}.
   \end{split}
 \end{equation}
In order to estimate the $BMO_\gamma$ norm of $v=\log u$ on $B_R(x_0)$, we need to estimate the quantity $\mu_\gamma(B_\rho(y))^{-1}\int_{B_\rho(y)}|v-\overline{v}_{B_\rho(y)}| |x|^{-\gamma} \dx$ on any ball $B_\rho(y) \Subset B_R(x_0)$.  To this end, we use the weighted Poincar\'e inequality \eqref{Weighted.Poincare.ineq}, H\"{o}lder inequality and estimate \eqref{upperbound.gradient.log}, and we obtain the following:
 \begin{equation*}\label{estimate.bmo.1}\begin{split}
   \bigg(\frac{1}{\mu_\gamma(B_\rho(y))} \int_{B_\rho(y)}|v(\tau)   &-  \overline{v(\tau)}_{B_\rho(y)}|\, |x|^{-\gamma} \dx\bigg)^{2}  \le   \frac{P_{\gamma,\beta}^2\, \rho^2}{\mu_\beta(B_\rho(y))}\int_{B_{2\rho}(y)} \psi^2 |\nabla v(\tau)|^2 |x|^{-\beta}  \dx\, ,\\
   &\leq \frac{4 c_N P_{\gamma,\beta}^2 D_{\beta}}{  m^2} +  \frac{2 P_{\gamma, \beta}^2 D_\gamma}{m^2(1-m)} \frac{\rho^2 \mu_\gamma(B_{\rho}(y))}{\mu_\beta(B_\rho(y))} \frac{M_p(u_0, \delta, \tau)^{1-m}}{\tau},
 \end{split}\end{equation*}
where $D_\gamma (D_\beta)$ is the doubling constant of the measure $\mu_\gamma, (\mu_\beta)$ respectively, defined in \eqref{Doubling.Constant}. Finally, by Lemma \ref{technical.lemma.measures}
\begin{equation*}\label{lemma7.3more}
\frac{\rho^2 \mu_\gamma(B_{\rho}(y))}{\mu_\beta(B_\rho(y))} \le \ka_{16} \left(\int_{B_{\rho}(x_0)} |x|^{(\beta -\gamma)\frac{N}{2}} \dx \right)^{\frac{2}{N}}
\le \ka_{16} \left(\int_{B_R(x_0)} |x|^{(\beta -\gamma)\frac{N}{2}} \dx \right)^{\frac{2}{N}}\le c_1\left(\frac{|x_0|}{R}\vee 1\right)^{\beta-\gamma}\,R^\sigma\,,
\end{equation*}
 where $\ka_{16},c_1>0$ only depend on $N,\gamma,\beta$ and $\ka_{16}$ given in \ref{rho.pseudo}. This concludes the proof. \qed

\noindent {\bf Proof of Proposition \ref{Lem.Step3}.~}We have shown  that  $\log u_\delta(t) \in BMO_\gamma(B_{R_1}(x_0))$ for any $t \geq t_0 > 0$ in Lemma \ref{Caccioppoli.JohnNirenberg.lemma.LDP}, more precisely, inequality \eqref{BMO.norm.delta} gives $\|\log u_\delta(t)\|_{BMO_{\gamma}(B_{R_1}(x_0))}\le 1/\ka_6 \nu_\delta$  for all $t\ge t_0$. We are now in the position to apply inequality \eqref{rev.hold.cor.intro} of Corollary \ref{Reverse.Holder.Corollay}, which gives inequality \eqref{Reverse.Holder.nu.delta}, namely
\begin{equation*}
\|u(t)\|_{\LL^{s}_\gamma(B_{R_1}(x_0))} \le \ka_{7}^{2/s} \mu_\gamma(B_{R_1}(x_0))^{2/s} \|u(t)\|_{\LL^{-s}_\gamma(B_{R_1}(x_0))}\,,
\end{equation*}
for all $t\ge t_0 > 0$, and all $0 < s <   \nu_\delta\le 1/ \big(\ka_{6} \|\log u(t)\|_{BMO_{\gamma}(B_{R_1}(x_0))}\big)$. Letting $s=\varepsilon$ concludes the proof.\qed

The above results extend to the case of the \ref{MINIMAL.DIRICHLET.PROBLEM} as in the following lemma.

\begin{prop}[Reverse H\"older inequality for \ref{MINIMAL.DIRICHLET.PROBLEM}] \label{Caccioppoli_JohnNirenberg_lemma}
Let $u$ be a solution to \ref{MINIMAL.DIRICHLET.PROBLEM} on $(0,\infty)\times B_{R_0}(x_0)$, corresponding to the initial datum $u_0\chi_{B_R(x_0)} \in \LL^p_\gamma(B_{R_0}(x_0))$   with $ p > p_c $ when $m\in (0,m_c]$ and with $p \geq 1  $ when $m\in (m_c,1)$, with $4R\le R_0$ and assume that $B_{R_0}(x_0)$ satisfies either (1), (2) or (3); let $T=T(u_0)$ be its extinction time.
Then, estimate \eqref{BMO.norm.delta} holds for $u$, with $\delta=0$.
Moreover, let $H_p$ be as in \eqref{def.HP}, then for every $\tau_*\in (0,1]$  we define
\begin{equation}\label{def.nu0}
 \nu_0 := \frac{\tau_*^{\sigma p\vartheta_p}}{\ka_6 \ka_{15}}\left[1+\left(\frac{|x_0|}{R}\vee 1\right)^{\beta-\gamma}
 H_p(u_0, x_0, R)^{1-m} \right]^{-\frac{1}{2}}>0\,,
\end{equation}
so that for every $t\in[\tau_* t_*,t_*]\subset(0,T)$ with $t_*=t_*(u_0, x_0, R)$ is as in  \eqref{choice_of_t}. We have
\begin{equation}\label{BMO.norm.estimate2}
\| \log u(t)\|_{BMO_{\gamma}(B_R(x_0))} \le (\ka_6 \nu_0)^{-1}\,.
\end{equation}
Finally, for all $t\in[\tau_* t_*,t_*]\subset(0,T)$ and all $0<R_1<R$, we have
\begin{equation}\label{Reverse.Holder}
\|u(t)\|_{\LL^{\varepsilon}_\gamma(B_{R_1}(x_0))} \le \ka_{7}^{2/ \varepsilon } \mu_\gamma(B_{R_1}(x_0))^{2/ \varepsilon } \|u(t)\|_{\LL^{-\varepsilon}_\gamma(B_{R_1}(x_0))}\qquad\mbox{for all $0 < \varepsilon <   \nu_0$.}
\end{equation}
The constant $\ka_{15}=\ka'_{15} [m(1-m)]^{-1}$  depends on $N, m, \gamma, \beta$ and $\ka_6,\ka_7>0$ are as in Corollary \ref{Reverse.Holder.Corollay}.
\end{prop}
\noindent\textbf{Proof.} Inequality \eqref{BMO.norm.estimate2} follows by letting $\delta\to 0$ in inequality \eqref{BMO.norm.delta},  exploiting the lower semicontinuity of the $BMO_{\gamma}(B_R(x_0))$-norm,  then substituting $t=\tau_*t_{\ast}<T$, with $t_*$ given in \eqref{choice_of_t}, and finally noticing that
\begin{equation*}
\frac{1}{\tau_*t_*}M_p(u_0, \delta, \tau_*t_*)^{1-m}\xrightarrow[]{\delta\to 0^+} \frac{c_3}{R^\sigma\tau_*^{\sigma p\vartheta_p}}  H_p(u_0, x_0, R)^{1-m}\,,
\end{equation*}
where $c_3>0$ depends on $N,\gamma,\beta$ and $m$. Inequality \eqref{Reverse.Holder} then follows as in the proof of Proposition \ref{Lem.Step3}.\qed

\medskip

\begin{rem}\label{remark.nu.0}\rm When we are in the good fast diffusion range, i.e $m \in (m_c, 1)$, we can choose $\nu_0$ independent of $u_0$\,, indeed, by letting $p=1$ in \eqref{def.nu0} and recalling that  $H_1(u_0,x_0, R)= \mu_\gamma(B_R(x_0))^{\sigma \vartheta_1} R^{-\sigma(N-\gamma)\vartheta_1}$, we have that
\begin{equation*}\label{def.nu0.p1}
 \nu_0 := \frac{\tau_*^{\sigma \vartheta_1}}{\ka_6 \ka_{15}}\left[1+\left(\frac{|x_0|}{R}\vee 1\right)^{\beta-\gamma}\left(
 \frac{\mu_\gamma(B_R(x_0))^{\sigma \vartheta_p}}{R^{\sigma(N-\gamma)\vartheta_p}}\right)^{1-m} \right]^{-\frac{1}{2}}\,.
\end{equation*}
This will have important consequences, but in particular we immediately obtain  an \textit{absolute bound of the BMO norm of $\log u$ }on intrinsic cylinders, namely $\| \log u(t)\|_{BMO_{\gamma}(B_R(x_0))}\le (\ka_6 \nu_0)^{-1}$, for all $t\in[\tau_* t_*,t_*]\subset(0,T)$.
Unfortunately the dependence on $u_0$ cannot be dropped in the very fast diffusion range, i.e. when $m\in (0,m_c]$.
\end{rem}

\subsection{End of Step 3. Proof of  Proposition \ref{lower.smoothing.final.form} and $\LL^{-\infty}-\LL^\varepsilon$ estimates for MDP}
We now sum up all the results of the first three Steps to prove the $\LL^{-\infty}-\LL^\varepsilon$ estimates for the \ref{LARGE.DIRICHLET.PROBLEM}. Next we prove analogous estimates for the MDP by letting $\delta\to0$. \\
\smallskip
\noindent {\bf Proof of  Proposition \ref{lower.smoothing.final.form} .~}Let us first fix $\varepsilon>0$ as in \eqref{nu.delta.def.Final3}, namely such that\vspace{-2mm}
\begin{equation*}
0 < \varepsilon <   \nu_\delta\wedge (1-m)\qquad\mbox{with}\qquad\nu_\delta = \frac{1}{\ka_{14} \, \ka_{6}} \left[ 1 + \frac{R_1^\sigma}{t_0}\left( \frac{|x_0|}{R_1}\vee 1\right)^{\beta-\gamma} M_p(u_0, \delta, t_0)^{1-m}\right]^{-1/2}.
\end{equation*}
where $M_p$ is given in \eqref{uniform.boundedness.approximate.solutions}, $\ka_{14}$ is as in Lemma \ref{Caccioppoli.JohnNirenberg.lemma.LDP}, so that the Reverse H\"older inequality \eqref{Reverse.Holder.nu.delta} holds. Then we are in the position to use the $L^{-s_\varepsilon}-\LL^{-\varepsilon}$ smoothing effect of Proposition \ref{prop.finite.iter.negative.norms} with $\varepsilon\in (0,1-m)$ as above\,, together with the $\LL^{-\infty}-\LL^{-s_\varepsilon}$ lower bounds of Proposition \ref{lower.spacetime.smoothing}\,; combining all the above results we obtain, choosing $\overline{R}=\left(R_1+R_2\right)/2$ and times as in  the statement:
\begin{equation}\label{Prop34.proof.1}
\begin{split}
  &\inf_{(T_1, T_2]\times B_{R_2}(x_0)}{u}
  \geq \kb_3 \left[\left(1+M_{p}(u_0,\delta, T_0)^{1-m}\right)
  \left(\frac{h_{\sigma}\left( \overline{R}, R_2, x_0 \right)}{\left(\overline{R}-R_2\right)^{\sigma}} + \frac{1}{T_{1}-T_0}\right) \right]^{\eta_{\varepsilon}}\nonumber \\[-1mm]
  & \times \left[\left(\frac{T_0}{T_2}\right)^{\frac{1}{1-m}} \left(\int_{B_{\overline{R}}(x_0)}u^{-s}(T_2, x) |x|^{-\gamma}\dx\right)^{-\frac{1}{s_\varepsilon}}\right]^{\frac{s_\varepsilon}{s_\varepsilon+m-1}}\left(T_2-T_0\right)^{-\frac{1}{s_\varepsilon+m-1}}\\[-1mm]
  &\geq \kb_3\kb_4^{\frac{s_\varepsilon}{s_\varepsilon+m-1}}\bigg[2^{(2\vee\sigma)}\left(1+M_{p}(u_0,\delta, T_0)^{1-m}\right)
  \left(\underline{\mu}\vee\left(\frac{h_\sigma(R_1, R_2, x_0)}{(R_1-R_2)^\sigma}+ \frac{1}{T_1-T_0}\right)\right) \bigg]^{\eta_{\varepsilon}+\frac{s_\varepsilon}{s_\varepsilon+m-1}\zeta_\varepsilon }\nonumber \\[-1mm]
  & \times \left[\left(\frac{T_0}{T_2}\right)^{\frac{1}{1-m}} \left(\frac{T_2}{T_3}\right)^{\frac{1}{1-m}}
  \left( \int_{B_{R_1}(x_0)}u_\delta(T_3, x)^{-\varepsilon}|x|^{-\gamma} \dx \right)^{-\frac{1}{\varepsilon}}\right]^{\frac{s_\varepsilon}{s_\varepsilon+m-1}}\left(T_2-T_0\right)^{-\frac{1}{s_\varepsilon+m-1}}\\[-1mm]
  &\geq \kb_3\kb_4^{\frac{s_\varepsilon}{s_\varepsilon+m-1}}\bigg[2^{(2\vee\sigma)
  }\left(1+M_{p}(u_0,\delta, T_0)^{1-m}\right)
  \left(\underline{\mu}\vee\left(\frac{h_\sigma(R_1, R_2, x_0)}{(R_1-R_2)^\sigma}+ \frac{1}{T_1-T_0}\right)\right) \bigg]^{\eta_{\varepsilon}+\frac{s_\varepsilon}{s_\varepsilon+m-1}\zeta_\varepsilon }\nonumber \\[-1mm]
  & \times \left[\left(\frac{T_0}{T_3}\right)^{\frac{1}{1-m}}
  \ka_{7}^{-\frac{2}{s_\varepsilon}} \mu_\gamma(B_{R_1}(x_0))^{-\frac{2}{s_\varepsilon}}
  \left( \int_{B_{R_1}(x_0)}u_\delta(T_3, x)^{\varepsilon}|x|^{-\gamma} \dx \right)^{\frac{1}{\varepsilon}}\right]^{\frac{s_\varepsilon}{s_\varepsilon+m-1}}\left(T_2-T_0\right)^{-\frac{1}{s_\varepsilon+m-1}}\\[-1mm]
\end{split}
\end{equation}
where we have put $\underline{\mu}:= \mu_\gamma(B_{R_2}(x_0))^{-\frac{\sigma}{N-\gamma}}$. In the first step  we have used \eqref{formula.inf} with $s_\varepsilon>1-m$ as in \eqref{const.finite.iter.negative.norms}, and $0 < R_2 < \overline{R} < R_1$ and $T_1\in(T_0, T_2)$. Notice that $\eta_{\varepsilon}:={-\frac{1}{s_\varepsilon+m-1}\left(\frac{N-\gamma}{2+\beta-\gamma}+1 \right)}$, $\zeta_\varepsilon = {-\frac{1}{\varepsilon}\frac{1-\left(\frac{2}{r^*}\right)^{\ke}}{1-\frac{2}{r^*}}}$; $h_\sigma$ and $M_p$ are defined in \eqref{function_h} and \eqref{uniform.boundedness.approximate.solutions} respectively, and $\kb_3>0$ is as in \eqref{formula.inf} only  depending  on $s, \widetilde{p}, N, \gamma, \beta, m$.
In the second step we have used \eqref{ineq.finite.iter.negative.norms} with $\kb_4$ as in \eqref{const.finite.iter.negative.norms}, noticing that $M_p(u_0, \delta, T_2)\le M_p(u_0, \delta, T_0)$ and that $h_{\sigma}\left( \overline{R}, R_2, x_0 \right)\asymp h_{\sigma}\left( R_1, \overline{R}, x_0 \right)\le  2^{(2-\sigma)_{+}}h_{\sigma}\left( R_1, R_2, x_0 \right) $, and where $s_\varepsilon>1-m$ and $\kb_4>0$ are given in \eqref{const.finite.iter.negative.norms}.
In the third step we have used the Reverse H\"older inequality \eqref{Reverse.Holder.nu.delta}  with $\varepsilon$ and $\nu_\delta$ as above, with $\ka_{7}$ are as in Corollary \ref{Reverse.Holder.Corollay}. \qed
\medskip

\noindent\textbf{Taking the limit $\delta \rightarrow 0$. $\LL^{-\infty}-\LL^\varepsilon$ interior estimates for MDP. }
Consider the solution $u$ of the \ref{MINIMAL.DIRICHLET.PROBLEM} with initial data $u_0$. Then the solutions of the ``lifted problem'' \ref{LARGE.DIRICHLET.PROBLEM} $u_\delta$ are ordered with respect to $\delta$: more precisely, for $\delta > \delta'$, for any $x \in B_R(x_0)$ and for any $t \in (0, \infty)$
\begin{equation*}
u_\delta(t,x) \ge u_{\delta'}(t,x).
\end{equation*}
In particular, for any $x \in B_R(x_0)$ and for any $t \in (0, \infty)$, the limit as $\delta \rightarrow 0$ exists and it equal to $u(t,x)$. See \cite[Section B.3]{BFR} for more details.
 Note  that the constants in the inequality $\eqref{formula.inf.finalBBBB}$ remain stable as $\delta\to 0^+$ (see also Proposition \ref{Caccioppoli_JohnNirenberg_lemma}). As an immediate consequence of Proposition \ref{lower.smoothing.final.form} we get

\begin{cor}[$\LL^{-\infty}-\LL^{\varepsilon}$ estimates for MDP]\label{weak.positivity.MDP}
Let $u$ be a strong (super)solution to \ref{MINIMAL.DIRICHLET.PROBLEM} on $(0,\infty)\times B_{R_0}(x_0)$, corresponding to the initial datum $u_0\chi_{B_R(x_0)} \in \LL^p_\gamma(B_{R_0}(x_0))$ with $ p > p_c $ when $m\in (0,m_c]$ and with $p \geq 1  $ when $m\in (m_c,1)$, moreover assume that $8R=R_0$  and that $B_{R_0}(x_0)$ satisfies either (1), (2) or (3); let $T=T(u_0)$ be the extinction time.
  Let $\widetilde{H}_p$ be as in \eqref{def.HP}, and define, for any fixed $\tau_*\in (0,\frac{1}{3})$,
\begin{equation}\label{def.nu0bbbb}
 \nu_0 := \frac{m(1-m) \tau_*^{\sigma p\vartheta_p}}{\ka'_{15} \widetilde{H}_p(u_0,x_0,R)^{\frac{1}{2}}}  \in (0,1-m)\,.
\end{equation}
Then, $[2\tau_*  \ta_{\ast},(1-\tau_*) \ta_{\ast}]\subset[\tau_*  \ta_{\ast}, \ta_{\ast}] \subset(0,T)$ with $ \ta_{\ast}= \ta_{\ast}(u_0, x_0, R)$ is as in \eqref{choice_of_t.1}.
Moreover, for any $\varepsilon\in (0,\nu_0)$ there exist $s_\varepsilon>1-m$ and $\kb_\varepsilon>0$ such that
\begin{align}\label{formula.inf.final}
  \inf_{[2\tau_* \ta_{\ast},(1-\tau_*)\ta_{\ast}]\times B_{2R}(x_0)}u
  &\geq \kb_\varepsilon
  \left( \int_{B_{4R}(x_0)}u(\ta_{\ast}, x)^{\varepsilon}|x|^{-\gamma} \dx \right)^{\frac{s_\varepsilon}{\varepsilon(s_\varepsilon+m-1)}}
\end{align}
with
\begin{align}\label{kappa-epsilon.MDP}
\kb_\varepsilon &:= \kb_2\,\tau_*^{\theta_\varepsilon}\mu_\gamma(B_{4R}(x_0))^{-\frac{2}{s_\varepsilon+m-1}} \left((1-2\tau_*)\ta_{\ast}\right)^{-\frac{1}{s_\varepsilon+m-1}}\\ \nonumber
&\times\left[ \widetilde{H}_p(u_0,x_0,R)\left( 1\vee \frac{\ta_{\ast}}{R^{\sigma}}\right) \left(\mu_\gamma(B_{2R}(x_0))^{-\frac{\sigma}{N-\gamma}}\vee\left(\frac{h_{\sigma}( 4R, 2R, x_0)}{\left(2R\right)^{\sigma}} + \frac{1}{\tau_*\ta_{\ast}}\right)\right)\right]^{\eta_{\varepsilon}+\frac{s_\varepsilon\, \zeta_\varepsilon }{s_\varepsilon+m-1}}
\end{align}
 where $s_\varepsilon>0, \eta_{\varepsilon}, \zeta_\varepsilon<0$ are as in \eqref{exponents.lower.smoothing.Final3},
$\theta_\varepsilon=(1-p\sigma \vartheta_p)\left(\eta_{\varepsilon}+\frac{s_\varepsilon\, \zeta_\varepsilon}{s_\varepsilon+m-1}\right)+\frac{1}{1-m}\frac{s_\varepsilon}{s_\varepsilon+m-1}$;
$h_\sigma$, $M_p$ and $\vartheta_p$ are defined in \eqref{function_h}, \eqref{uniform.boundedness.approximate.solutions} and \eqref{definition.sigma.thethap} respectively;
finally, $\kb_2>0$ depends only on  $N,m,p,\beta,\gamma, \varepsilon$ , through $\kb_3,\kb_4$ defined in \eqref{formula.inf}, \eqref{const.finite.iter.negative.norms}, and through $\ka_6,\ka_7 $, which are defined in Corollary \ref{Reverse.Holder.Corollay};
$\ka'_{15}$ is the same as in Proposition \ref{Caccioppoli_JohnNirenberg_lemma}.
\end{cor}

\subsection{Step 4. Reverse $\LL^1-\LL^\varepsilon$ smoothing effects and interior lower bounds for MDP. }\label{ssec.step4}
Next  we obtain  a useful Lemma about quantitative positivity of local $\LL^1_\gamma$ norms and a local reverse $\LL^\varepsilon-\LL^1$ smoothing effects for solutions to the MDP.
\begin{lem}\label{positivity.MDP.thm}
Let $u$ be the solution to \ref{MINIMAL.DIRICHLET.PROBLEM} on $(0,\infty)\times B_{R_0}(x_0)$, corresponding to the initial datum $u_0\chi_{B_R(x_0)} \in \LL^p_\gamma(B_{R_0}(x_0))$  with $ p > p_c $ when $m\in (0,m_c]$ and with $p \geq 1  $ when $m\in (m_c,1)$, with $4R\le R_0$  and assume that $B_{R_0}(x_0)$ satisfies either (1), (2) or (3); let $T=T(u_0)$ be its extinction time and define $0\le \ta_*\le T$ as
\begin{equation} \label{choice_of_t.1}
 \ta_{\ast}= \ta_{\ast}(u_0, x_0, R)= \frac{ \ka_{10}'^{-1} }{2^{2-m}} \frac{R^\sigma}{\mu_\gamma(B_R(x_0))^{1-m}} \|u_0\|_{\LL^1_{\gamma}(B_R(x_0))}^{1-m},
\end{equation}
where $ \ka_{10}'\ge 1 $ is the constant defined in \eqref{estimate.testfunction}
 depending  only on $N, \gamma, \beta$ and $m$.
Then there exists $\kb_0>0$  such that
\begin{equation}\label{positivity.MDP}
    \frac{ \kb_0}{\mu_\gamma(B_R(x_0))}\int_{B_R(x_0)}{u_0 |x|^{-\gamma} \dx} \leq \frac{1}{\mu_\gamma(B_{4R}(x_0))} \int_{B_{4R}(x_0)}{u\left(\ta_{\ast} , x \right) |x|^{-\gamma}  \dx},
\end{equation}
where $ \kb_0$ depends on $N, \gamma, \beta, m$. Moreover, for  any $ \varepsilon \in \left(0, 1\right)$  the following estimate  holds
\begin{equation}\label{positivity.general.H}
  \left(\frac{1}{\mu_\gamma(B_R(x_0))} \int_{B_R(x_0)}{u_0 |x|^{-\gamma}  \dx} \right)^\varepsilon \leq \kb_0^{-1}\kb_9^{1-\varepsilon} \frac{H_p\left( u_0, x_0, R \right)^{1-\varepsilon} }{\mu_\gamma(B_{4R}(x_0))} \int_{B_{4R}(x_0)}{u^\varepsilon( \ta_{\ast}, x) |x|^{-\gamma} \dx},
\end{equation}
where $ H_p \left( u_0, x_0, R \right) $ is defined in \eqref{def.HP} , $\kb_9=\ka_{12} \left(\ka_{10}' 2^{2-m}\right)^{(N-\gamma)\vartheta_p} \omega_{N, \gamma}^{\sigma\vartheta_p}$, with $\ka_{12}>0$ depending only on $N, m, \gamma, \beta, p$ (given in \eqref{SmoothingEffectInequality}) and $\omega_{N, \gamma}$ being such that $ \omega_{N, \gamma}R^{N-\gamma}=\mu_\gamma(B_R(0))$.
\end{lem}
\noindent\textbf{Proof.} Let $u(t,x)$ be a solution to \ref{MINIMAL.DIRICHLET.PROBLEM} over the cylinder $B_{R_0}(x_0) \times (0, T)$.
Applying inequality \eqref{Herrero.Pierre.123} with times $t=0, \tau=\ta_{*}$ and radii  $R$ and $2R$  we obtain
\begin{equation}\label{positivity.through.HerreroPierre}\begin{split}
  \frac{1}{\mu_\gamma({B_R(x_0)})}&\int\limits_{{B_R(x_0)}}u_0 \frac{\dx}{|x|^{\gamma}}
  \le \frac{2^{\frac{1}{1-m}}}{\mu_\gamma({B_R(x_0)})}
  \left[ \int\limits_{B_{4R}}u(\ta_{\ast}, x) \frac{\dx}{|x|^{\gamma}} + \frac{\ka_{10}'^{\frac{1}{1-m }}\left(\ta_{\ast} \right)^{\frac{1}{1-m}}}{R^{\frac{\sigma}{1-m}}} \mu_\gamma({B_R(x_0)}) \right]  \\
  & \le \frac{2^{\frac{1}{1-m}} D_\gamma^2}{\mu_\gamma(B_{4R}(x_0))}\int\limits_{B_{4R}(x_0)}u(\ta_{\ast}, x) \frac{\dx}{|x|^{\gamma}} + \frac{1}{2}\frac{1}{\mu_\gamma(B_{R})}\int\limits_{{B_R(x_0)}}u_0 \frac{\dx}{|x|^{\gamma}},
\end{split} \end{equation}
where we used the fact that $u_0$ is supported in $B_R$, the doubling property of the measure $\mu_\gamma$ and the fact that $u(t,x)>0$. Inequality  \eqref{positivity.MDP} is then deduced from \eqref{positivity.through.HerreroPierre} with constant $\kb_0 = \frac{1}{2} 2^{-\frac{1}{1-m}}D_\gamma^{-2}$ .
We now turn our attention to inequality \eqref{positivity.general.H}, which will be deduced from \eqref{positivity.MDP}.
 Let $\varepsilon \in (0,1)$: using the smoothing-effect  inequality \eqref{SmoothingEffectInequality}, namely
$\|u(t)\|_{\LL^{\infty}(B_{R_0}(x_0))}\leq \ka_{12} \|u_0\|^{\sigma p\vartheta_p}_{\LL^p_{\gamma}(B_R(x_0))} t^{-(N-\gamma)\vartheta_p}$
(recall that $\ka_{12}$ does not depend on $B_R(x_0)$) applied at $t=\ta_{\ast}$ we get
\begin{equation*}\begin{split}
   \|u(\ta_{\ast})\|_{\LL^{\infty}(B_{R_0})}
    & \le  \kb_9 H_p(u_0,x_0, R) \frac{\|u_0\|_{\LL^1_\gamma(B_{R}(x_0))}}{\mu_\gamma(B_R(x_0))}
\end{split}\end{equation*}
 where in the last step we have used inequality \eqref{positivity.MDP} and the equality $ \omega_{N, \gamma}R^{N-\gamma}=\mu_\gamma(B_R(0))$; the constant $\kb_9=\ka_{12} \left(\ka_{10}' 2^{2-m}\right)^{(N-\gamma)\theta_p} \omega_{N, \gamma}^{\sigma\vartheta_p} $.  We finally combine the above inequality with  $\|u(\ta_*)\|_{\LL^1_\gamma(B_{4R}(x_0))}\le \|u(\ta_{\ast})\|_{\LL^{\infty}(B_{4R}(x_0))}^{1-\varepsilon}\|u(\ta_{\ast})\|^\varepsilon_{\LL^\varepsilon_\gamma(B_{4R}(x_0))}$, and we obtain inequality \eqref{positivity.general.H}.\qed

\noindent Putting together all the results of the 4 Steps, we obtain the following
\begin{cor}[$\LL^{-\infty}-\LL^{1}$ estimates for MDP]\label{weak.positivity.MDP.L1}
Let $u$ be a strong (super)solution to \ref{MINIMAL.DIRICHLET.PROBLEM} on $(0, \infty) \times B_{R_0}(x_0)$, corresponding to $u_0\chi_{B_R(x_0)} \in \LL^p_{\gamma}(B_{R_0}(x_0))$ with  $p > p_c$ if $m\in(0,m_c]$ or $p \geq 1$ if $m\in (m_c,1)$,  moreover assume that $4R=R_0$ and that $B_{R_0}(x_0)$ satisfies either (1), (2) or (3); let $T=T(u_0)$ be its extinction time and define $ t_*\in [0,T]$  by
\begin{equation}\label{t.star.MDP.L1}
 t_*= t_*(u_0, x_0, R)=  \kappa_* R^\sigma  \frac{\|u_0\|_{\LL^1_{\gamma}(B_R(x_0))}^{1-m}}{\mu_\gamma(B_R(x_0))^{1-m}}
\end{equation}
where  $\kappa_*=5^{-1}2^{m}\ka_{10}'^{-1}$,   $ \ka_{10}'\ge 1 $  depends  only on $N, \gamma, \beta$ and $m$ as in \eqref{estimate.testfunction}.
Then, there exists $\kb>0$ such that
\begin{align}\label{ineq.weak.positivity.MDP.L1}
  &\inf_{[\frac{t_*}{2},t_*]\times B_{2R}(x_0)}u
    \ge \kb(\widetilde{H}_p, R)\,  t_*^{\frac{r^*}{(1-m)(r^*-2)} },
\end{align}
where $\kb$ has an explicit form given in \eqref{const.lower.cor.L1},  in particular for $R$ bounded  and $\widetilde{H}_p>>1$ we have
\begin{equation}\label{Hp.tilde.def}
\kb\asymp \left(\frac{R^{c_1}}{\widetilde{H}_p^{c_2}}\right)^{\frac{\widetilde{H}_p^{1/2}}{m(1-m)}} \quad\mbox{with}\quad\widetilde{H}_p(u_0, x_0, R):= 1+\left(\frac{|x_0|}{R}\vee 1\right)^{\beta-\gamma} H_p(u_0, x_0, R)^{1-m} \ge 1\,,
\end{equation}
and with $c_1,c_2>0$ only  depending  on $N,m,p,\beta,\gamma$.
\end{cor}
\noindent {\bf Proof.~}Let $H_p$ be as in \eqref{def.HP} and fix $\overline{\tau}_*=5^{-1}$,   $\kappa_*=(1-\overline{\tau}_*)2^{m-2}\ka_{10}'^{-1}$ , and $t_*:=(1-\overline{\tau}_*)\overline{t}_*$, where  $\overline{t}_*$ is defined  in \eqref{choice_of_t.1}; recall that $D_\gamma\ge 1$ is the doubling constant of the measure $\mu_\gamma$, and $ \ka_{10}'\ge 1  $ is as in \eqref{estimate.testfunction} (since $\ka_{10}'$ is a constant in an upper bound, hence without loss of generality we can take it bigger than 1). Next we fix $\nu_0$ (that depends on $ \tau_*$)  as in \eqref{def.nu0bbbb} with $\tau_*=\overline{\tau}_*$. Finally, we choose $\varepsilon=\varepsilon_0:= (2/r^*)^{k_0}(1-m)$, where $k_0$ is the smallest integer such that   $\varepsilon_0<\nu_0$. Note that  we have  $k_0\ge \log(\nu_0)/\log(2/r^*)-\log(1-m)/\log(2/r^*)$.
With these choices, we know that $k_\varepsilon=k_0+1$ and the exponents $s_\varepsilon>1-m, \eta_{\varepsilon}, \zeta_\varepsilon<0$ given in \eqref{exponents.lower.smoothing.Final3} become
\begin{align*}
& s_0:=s_{\varepsilon_0} =\left(\frac{r^*}{2}\right)^{\ke}{\varepsilon_0} =(1-m)\frac{\sr}{2}>1-m  \,, \quad
     \frac{s_0}{s_0+m-1}=\frac{r^*}{r^*-2} >0  \\
&  \eta_{s_0}:=\eta_{\varepsilon_0}=-2\frac{\frac{N-\gamma}{2+\beta-\gamma}+1}{(1-m)(r^*-2)} <0\,,
 \quad
\zeta_0:=\zeta_{\varepsilon_0}={-\frac{1}{\varepsilon_0}\frac{r^*(1-m)- 2\varepsilon_0 }{(1-m)(r^*-2)}}<0\, .\nonumber
\end{align*}
 Note  that, even if $\varepsilon_0$ depends on $\nu_0$ (through $k_0$), the exponents $s_\varepsilon$ and $\eta_{\varepsilon}$ now only depend on $N,\gamma,\beta, m$.
We are now in the position to combine inequalities \eqref{formula.inf.final} and \eqref{positivity.general.H} (with $R_0\ge 4R$)
as follows: for any $0< R_1<4R\le R_0$ we have
\begin{align*}
  &\inf_{[ 2\overline{\tau}_* \overline{t}_*,  (1-\overline{\tau}_*)\overline{t}_*]\times B_{2R}(x_0)}u
  \geq \kb_{\varepsilon_0}
    \left( \int_{B_{4R}(x_0)}u(\overline{t}_*, x)^{{\varepsilon_0}}|x|^{-\gamma} \dx \right)^{\frac{s_{0}}{{\varepsilon_0}(s_{0}+m-1)}}\\
  &\ge \kb_{\varepsilon_0} \left(\frac{\mu_\gamma(B_{4R}(x_0))}{\kb_9
    H_p\left( u_0, x_0, R\right)^{1-{\varepsilon_0}}}\right)^{\frac{s_{0}}{{\varepsilon_0}(s_{0}+m-1)}}
    \left(\frac{1}{\mu_\gamma(B_R(x_0))}\int_{B_R(x_0)}{u_0 |x|^{-\gamma}  \dx}\right)^{\frac{s_{0}}{s_{0}+m-1} } \nonumber\\
    &\ge \kb(\widetilde{H}_p,R) \,  t_*^{\frac{s_{0}}{(1-m)(s_{0}+m-1)} } =\kb(\widetilde{H}_p,R)\,  t_*^{\frac{r^*}{(1-m)(r^*-2)} }\nonumber
\end{align*}
where in the second step we have used inequality \eqref{positivity.general.H} and  $\kb_\varepsilon$ is as in \eqref{kappa-epsilon.MDP} and $\kb_9$ is as in \eqref{positivity.general.H}. Finally we have used the expression of $t_*$ given in \eqref{t.star.MDP.L1}.  We estimate the  $\kb(\widetilde{H}_p,R)$  as follows
\begin{align}
 &\kb_{\varepsilon_0} \left(\frac{\mu_\gamma(B_{4R}(x_0))}{\kb_9
    H_p\left( u_0, x_0, R\right)^{1-{\varepsilon_0}}}\right)^{\frac{s_{0}}{{\varepsilon_0}(s_{0}+m-1)}} \\ \nonumber
   &\geq  \kb_2\,\overline{\tau}_*^{\theta_\varepsilon}\mu_\gamma(B_{4R}(x_0))^{\frac{s_0-2\varepsilon_0}{\varepsilon_0(s_0+m-1)}} \left((1-2\overline{\tau}_*)\ta_{\ast}\right)^{-\frac{1}{s_0+m-1}}\widetilde{H}_p\left( u_0, x_0, R\right)^{\frac{(\varepsilon_0-1)s_{0}}{{\varepsilon_0}(s_{0}+m-1)}}\kb_9^{-\frac{s_{0}}{{\varepsilon_0}(s_{0}+m-1)}}\\ \nonumber
&\times\left[ \widetilde{H}_p(u_0,x_0,R)\left( 1\vee \frac{\ta_{\ast}}{R^{\sigma}}\right) \left(\mu_\gamma(B_{2R}(x_0))^{-\frac{\sigma}{N-\gamma}}\vee\left(\frac{h_{\sigma}( 4R, 2R, x_0)}{\left(2R\right)^{\sigma}} + \frac{1}{\overline{\tau}_*\ta_{\ast}}\right)\right)\right]^{\eta_{s_0}+\frac{s_0\, \zeta_0 }{s_0+m-1}},
\end{align}
where we just used the expression of $\widetilde{H}_p=\widetilde{H}_p(u_0, x_0, R)$ given in \eqref{Hp.tilde.def} and rewritten in the constant appearing in \eqref{kappa-epsilon.MDP}. We then estimate
\begin{align}
&\left[ \widetilde{H}_p(u_0,x_0,R)\left( 1\vee \frac{\ta_{\ast}}{R^{\sigma}}\right) \left(\mu_\gamma(B_{2R}(x_0))^{-\frac{\sigma}{N-\gamma}}\vee\left(\frac{h_{\sigma}( 4R, 2R, x_0)}{\left(2R\right)^{\sigma}} + \frac{1}{\overline{\tau}_*\ta_{\ast}}\right)\right)\right]^{\eta_{s_0}+\frac{s_0\, \zeta_0 }{s_0+m-1}}  \\ \nonumber
&\geq \widetilde{H}_p(u_0,x_0,R)^{\eta_{s_0}+\frac{s_0\, \zeta_0 }{s_0+m-1}}\left[ \left(1 \vee \frac{t_*}{(1-\overline{\tau}_*)R^{\sigma}} \right)\right. \times  \\ \nonumber &\hspace{5mm}\times \left.\left(\frac{h_{\sigma}( 4R, 2R, x_0)}{(2R)^{\sigma}} \vee \mu_\gamma(B_{2R}(x_0))^{-\frac{\sigma}{N-\gamma}}\right)\left(1 + \frac{(2R)^{\sigma}(1-\overline{\tau}_*)}{\overline{\tau}_*t_*}\right) \right]^{\eta_{s_0}+\frac{s_0\zeta_0}{s_0+m-1}},
\end{align}
where we just used the expression of $\overline{t}_*$ given in the beginning of the proof. Recall that the expression of $\kb_2$, given in Proposition \ref{lower.smoothing.final.form} is:
\[
\kb_2=\kb_3\kb_4^{\frac{s_0}{s_0+m-1}}\ka_{7}^{-\frac{2}{s_\varepsilon+m-1}}2^{(2\vee\sigma)(\eta_{s_0}+\frac{s_0}{s_0+m-1}\zeta_0)},
\]
where $\kb_4=\kb_4'\mu_\gamma(B_{4R})^{-\frac{1}{\varepsilon_0}\frac{r^*-2}{r^*+2}}$. All the above estimates finally give the expression of $\kb$
\begin{align}\label{const.lower.cor.L1}
\kb:=&\kb_3\kb_4'^{\frac{s_0}{s_0+m-1}}\ka_{7}^{-\frac{2}{s_\varepsilon+m-1}}2^{(2\vee\sigma)(\eta_{s_0}+\frac{s_0}{s_0+m-1}\zeta_0)}\kb_9^{-\frac{s_{0}}{{\varepsilon_0}(s_{0}+m-1)}} \\ \nonumber \times &  \overline{\tau}_*^{\theta_\varepsilon} \mu_\gamma(B_{4R}(x_0))^{-\frac{2-\frac{s_0}{\varepsilon_0}(\frac{2\sr}{\sr+2})}{s_0+m-1}}\, \widetilde{H}_p^{\eta_{s_0}+\frac{s_0}{s_0+m-1}(\zeta_0 + \frac{\varepsilon_0-1}{\varepsilon_0})} \left(\frac{1-2\overline{\tau}_*}{1-\overline{\tau}_*}t_{\ast}\right)^{-\frac{1}{s_0+m-1}} \\ \nonumber
\times & \left[ \left(1 \vee \frac{t_*}{(1-\overline{\tau}_*)R^{\sigma}} \right)\left(\frac{h_{\sigma}( 4R, 2R, x_0)}{(2R)^{\sigma}} \vee \mu_\gamma(B_{2R}(x_0))^{-\frac{\sigma}{N-\gamma}}\right)\left(1 + \frac{(2R)^{\sigma}(1-\overline{\tau}_*)}{\overline{\tau}_*t_*}\right) \right]^{\eta_{s_0}+\frac{s_0\zeta_0}{s_0+m-1}}.
\end{align}
 We also recall that as in  Lemma \ref{positivity.MDP.thm},  $\kb_9=\ka_{12} \left(\ka_{10}' 2^{2-m}\right)^{(N-\gamma)\vartheta_p} \omega_{N, \gamma}^{\sigma\vartheta_p}$,  with $\ka_{12}>0$ that  depends only on $N, m, \gamma, \beta, p$ and is given in \eqref{SmoothingEffectInequality} and $\omega_{N, \gamma}$ being such that $ \omega_{N, \gamma}R^{N-\gamma}=\mu_\gamma(B_R(0))$.
Notice that for sufficiently small $\varepsilon_0<\nu_0\sim m(1-m)/\widetilde{H}_p^{1/2}$, we have that $\zeta_0 + \frac{\varepsilon_0-1}{\varepsilon_0} \sim  \frac{c}{\varepsilon_0}$. Finally, notice that   when $R$ is bounded  and $H_p $ is large enough, we have that $\kb\sim \big(R^{c_4}/\widetilde{H}_p^{c_5}\big)^{\widetilde{H}_p^{1/2}/m(1-m)} $, where $c_i>0$ only depend on $N,m,p,\beta,\gamma$.\qed

\subsection{Positivity for local solutions. End of the proof of Theorems \ref{LOCAL.LOWER.BOUNDS} and \ref{positivity.MDP.thm.general}. }\label{sec.final.proofs.positivity}

We are now in the position to conclude the proof of the main  results  of this Part, Theorems \ref{LOCAL.LOWER.BOUNDS} and \ref{positivity.MDP.thm.general}.

\noindent\textbf{End of the proof Theorem \ref{positivity.MDP.thm.general}. } Let $u(t, x)$ be a solution to the \ref{MINIMAL.DIRICHLET.PROBLEM} on the cylinder $B_{4R}(x_0) \times (0, T)$, where $T=T(u_0)$ is the extinction time. Recall that $0\le u_0 \in \LL^p_{\gamma}(B_R(x_0))$ for $p > p_c$ if $0 < m \leq m_c$ or $p \geq 1$ if $m_c < m < 1$.
Let $M=\int_{B_R(x_0)}u_0 |x|^{-\gamma} \dx > 0$ and define the rescaled solution $\hat u$ as follows
\begin{equation*}
 \hat{u}\left(\hat t, \hat x \right)= \frac{R^{N-\gamma}}{M} u\left(\tau \hat t, R\hat x \right), \qquad \tau= R^{\sigma - (N-\gamma)(1-m)} M^{1-m}.
\end{equation*}
The rescaled solution $\hat{u}$ solves the \ref{MINIMAL.DIRICHLET.PROBLEM} on the cylinder $B_{4}(R^{-1} x_0)\times (0, \hat{T}) $ with mass $1$ and extinction time $\hat{T}$. We are in the position to apply Corollary \ref{weak.positivity.MDP.L1} to get: (recall that $\hat x= R^{-1}x$ and $\hat t= \tau^{-1}t$)
\begin{align}\label{ineq.weak.positivity.MDP.L1.resc}
  &\inf_{x\in B_{2}(\hat x_0)}\hat u(\hat t_*,\hat x)
    \ge \kb(\widetilde{H}_p, 1)\,  \hat t_*^{\frac{r^*}{(1-m)(r^*-2)} }\,, \qquad\mbox{where}\qquad \hat t_*=
    \kappa_* \frac{ 1}{\mu_\gamma(B_1(\hat x_0))^{1-m}}
\end{align}
where the value of $\kb(\widetilde{H}_p, 1)$ is given in \eqref{const.lower.cor.L1}, while $\widetilde{H}_p(\hat u_0, \hat x_0, 1)$ and $\kappa_*$ are given in Corollary \ref{weak.positivity.MDP.L1}.
Note that the quantity $\widetilde{H}_p$ is actually scaling invariant, namely
\begin{equation}\label{HP.scaling.inv}
\widetilde{H}_p(\hat u_0, \hat x_0, 1)=\widetilde{H}_p( u_0,  x_0, R).
\end{equation}
Note also  that $\hat t_*$ only depends on $N, m, \gamma, \beta$, but not on $u,u_0$ nor  $R, x_0$ (here is where we use either assumption (1), (2) or (3)); indeed, using $\mu_\gamma (B_1( \rho^{-1} x_0)) =  \rho^{\gamma - N} \mu_\gamma (B_{ \rho}( x_0))$ it is straightforward to check that $\hat t \asymp \k_*$.
Recalling that $\hat t\mapsto \hat t^{-\frac{1}{1-m}}\hat u(\hat t,\hat x)$ is non-increasing in time for almost every $\hat x \in B_{1}(\hat x_0)$\,, we get as a consequence of \eqref{ineq.weak.positivity.MDP.L1.resc}, for all $0\le \hat t\le \hat t_*$
\begin{align}\label{ineq.weak.positivity.MDP.L1.resc.2}
  &\inf_{\hat x\in B_{2}(\hat x_0)}\hat u(\hat t,\hat x)\ge  \left(\frac{\hat t}{\hat t_*}\right)^{\frac{1}{1-m}}\inf_{\hat x\in B_{2}(\hat x_0)}\hat u(\hat t_*,\hat x)
    \ge \kb(\widetilde{H}_p, 1)\,  \hat t_*^{\frac{2}{(1-m)(r^*-2)} } \hat t^{\frac{1}{1-m}}:= \hat \kb(\widetilde{H}_p, 1)\,   \hat t^{\frac{1}{1-m}}\,.
\end{align}
We have used a scaling argument to obtain a cleaner constant $\kb$ in the final lower bound \eqref{LOCAL.LOWER.INEQUALITY.MDP}, in this way, $\kb=\kb(\widetilde{H}_p, 1)$ shall depend on $R,x_0$ only through $H_p$. This is a consequence of our assumptions (1),(2) or (3) and the explicit expression of $\hat \kb(\widetilde{H}_p, 1)$ given in \eqref{const.lower.cor.L1} (recall that $\overline{\tau}_*=1/5$):
\begin{align}
\kb &:= \kb_3\kb_4'^{\frac{s_0}{s_0+m-1}}\ka_{7}^{-\frac{2}{s_\varepsilon+m-1}}2^{(2\vee\sigma)(\eta_{s_0}+\frac{s_0}{s_0+m-1}\zeta_0)}\kb_9^{-\frac{s_{0}}{{\varepsilon_0}(s_{0}+m-1)}} 5^{-\theta_{\varepsilon_0}} \left(3 \hat{t}_*\right)^{-\frac{1}{s_0+m-1}}\\ \nonumber
&C_{N, \gamma}^{-\frac{2-\frac{s_0}{\varepsilon_0}(\frac{2\sr}{\sr+2})}{s_0+m-1}}\, \widetilde{H}_p^{\eta_{s_0}+\frac{s_0}{s_0+m-1}(\zeta_0 + \frac{\varepsilon_0-1}{\varepsilon_0})}\left[ \left(1 \vee \frac{5}{4}\hat t_* \right) C_{\gamma, \beta} \left(1 + \frac{4}{  \hat t_*}\right) \right]^{\eta_{s_0}+\frac{s_0\zeta_0}{s_0+m-1}} \nonumber
\end{align}

In the computation of the above constant, we have used systematically the identity $\mu_\gamma (B_1( \rho^{-1} x_0)=  \rho^{\gamma - N} \mu_\gamma (B_{ \rho}( x_0))$ which holds under our assumptions: as a consequence all the constants in the right-hand side of formula \eqref{const.lower.cor.L1} will depend only on $N,m,\gamma,\beta$, and some of them on $H_p$. More precisely, $\kb_9 \asymp \kb_0^{-1}\ka_{12}$ as well as $\mu_\gamma(B_4(R^{-1}x_0)) \asymp C_{N, \gamma}$ and  $ h_{\sigma}\left( 4, 2, R^{-1}x_0 \right) \asymp C_{\gamma, \beta} $, where $C_{N, \gamma}, C_{\gamma, \beta}$ only depend on $N, \gamma$ and $\gamma, \beta$ respectively. Recall also that  $R^\sigma\asymp R^2 \mu_\gamma(B_R(x_0))/ \mu_\beta(B_R(x_0))$. Finally, we observe that when $\widetilde{H}_p$ is large, we have $\hat\kb(\widetilde{H}_p)\asymp  \widetilde{H}_p^{-\frac{c_2\widetilde{H}_p^{1/2}}{m(1-m)}} $, where $c_2$ only depends on $N, m, \gamma, \beta$.  Undoing the rescaling we obtain the lower bound \eqref{LOCAL.LOWER.INEQUALITY.MDP} and the proof of Theorem \ref{positivity.MDP.thm.general} is concluded.\qed

\medskip

\noindent\textbf{Proof of Theorem \ref{LOCAL.LOWER.BOUNDS}. }Once the positivity result is proven for solutions to the \ref{MINIMAL.DIRICHLET.PROBLEM}, namely once Theorem \ref{positivity.MDP.thm.general} is established, then by a standard comparison argument, the positivity result can be extended to any nonnegative local (super)solution. For strong (super) solutions the result is immediate, while for more general concepts of solutions, such as weak energy or very weak, the proof follows by a long but straightforward limiting process; see \cite{HerreroPierre} and \cite{JLVPorousMedioum} for more details about the non-weighted case; the case with weights follows along similar lines.\qed


\section{Part III. Harnack inequalities and H\"older continuity}  In this third part of the paper we study  regularity estimates for nonnegative solutions to both linear and nonlinear equations.

\subsection{The linear case}\label{sec.linear.regularity}
We are going to prove Harnack inequalities and local space-time H\"older continuity for nonnegative local solutions to the linear equation with Caffarelli-Kohn-Nirenberg weights. The equation
\begin{equation}\label{WHE.DW.GEN}
 v_t=w_\gamma\sum_{i,j=1}^N\partial_i  \left( A_{i,j}(t,x)\,  \partial_j v  \right),
\end{equation}
is posed on the cylinder $Q:=(0,T)\times\Omega$, where $A_{i,j}=A_{j,i}$ and for some $\gamma,\beta < N$ satisfying \eqref{paramt.range}, i.e. $\gamma -2 < \beta \leq \left(\frac{N-2}{N}\right) \gamma$,  as well we suppose that there exist constants $0 < \lambda_0 < \lambda_1 < +\infty$ such that
\begin{equation}\label{WHE.DW.GEN.paramt}
w_\gamma\asymp |x|^{\gamma} \qquad\mbox{and}\qquad 0<\lambda_0|x|^{-\beta} |\xi|^2\le \sum_{i,j=1}^N A_{i,j}(t,x)\xi_i \xi_j\le \lambda_1 |x|^{-\beta}|\xi|^2\,.
\end{equation}
The regularity estimates that we present in this section  are not present in the literature in the full range of parameters that we consider here, but several results have been obtained in different settings, see \cite{CS-AA87,CS-RSMUP85,CS-AMPA84,CS-CPDE84,FGS86, GW,GW2,Moser,MoserCpam71, SY}. We will only sketch the proofs, since they are minor modifications  of those obtained by Chiarenza-Serapioni and Gutierrez-Wheeden, \cite{CS-RSMUP85, GW, GW2} combined with the original proof of Moser \cite{MoserCpam71}. We shall keep  track of  the dependence of the Harnack constant by $\lambda_0,\lambda_1$ in a quantitative way as in \cite{MoserCpam71}, since in the nonlinear case this  will have  remarkable consequences.

In this weighted setting the Harnack inequality holds on suitable cylinders which take into account the geometry of the  problem;  recall that under assumptions (1), (2) or (3) we have
\[
 \rho_{x_0}^{\gamma,\beta}(R) := \left( \int_{B_R(x_0)} |x|^{(\beta - \gamma) \frac{N}{2}} \dx \right)^{\frac{2}{N}}
 \asymp\frac{\mu_\gamma(B_R(x_0))}{\mu_\beta(B_R(x_0))} R^2 \asymp R^{2+\beta-\gamma}.
\]
The following cylinders generalize the standard parabolic ones:
\begin{align}\label{Parabolic.cylinders.lin}
Q_R(t_0,x_0)&:=\left\{(t,x) \in \RR^+ \times \RR^d :\, t_0- \rho_{x_0}^{\gamma,\beta}(R)  <t\le t_0\,, \  |x-x_0|<2\,R  \right\}\,,\nonumber\\
Q_R^+(t_0,x_0)&:=\left\{(t,x) \in \RR^+ \times \RR^d :\,t_0-\frac{1}{4}\,\rho_{x_0}^{\gamma,\beta}(R)<t\le t_0\,, \ |x-x_0|<\frac{1}{2}\,R \right\}\,,\\
Q_R^-(t_0,x_0)&:=\left\{(t,x) \in \RR^+ \times \RR^d :\,t_0-\frac{7}{8}\,\rho_{x_0}^{\gamma,\beta}(R) <t\le t_0-\frac{5}{8}\,\rho_{x_0}^{\gamma,\beta}(R)\,, \ |x-x_0|<\frac{1}{2}\,R \right\}\,.\nonumber
\end{align}

It is convenient to introduce a suitable parabolic quasi-metric which carries on the information of the weights. Let $(t, x), (s,y) \in (0, \infty)\times\RR^N$ we define
\begin{equation}\label{parabolic.quasi.metric}
d_{\gamma, \beta}\left((t,x),(s,y)\right):= |x-y|\vee  \left(\rho^{\gamma,\beta}_{\overline{xy}}\right)^{-1}(|t-s|),
\end{equation}
where $\overline{xy}:=(x+y)/2$; the behaviour of $\left(\rho^{\gamma,\beta}_{\overline{xy}}\right)^{-1}$ is analyzed in Lemma \ref{inverse.technical.inequality.measures.lemma}.  Similar quantities have already been introduced in \cite{GN-CPDE88} in order to prove H\"{o}lder continuity of the solutions to weighted parabolic equation similar to \eqref{WHE.DW.GEN.paramt}, but with different weights.
In \cite{CS-CPDE84} it has been observed that although they imply continuity, they do not always imply H\"older continuity. For general classes of weights it is not possible to deduce any uniform modulus of continuity with respect to a standard parabolic quasi-distance.
However, for our class of weights, we still manage to deduce H\"older continuity from Harnack inequalities. Indeed, the quasi-metric $d_{\gamma, \beta}$ is controlled (on bounded space-time domains) by the following, more standard, parabolic quasi-distance (see Lemma \ref{inverse.technical.inequality.measures.lemma}):
\begin{equation}\label{standard.parabolic.quasi.metric}
\tilde{d}_\sigma\left((t,x),(s,y)\right):=\left\{
\begin{array}{lll}
|x-y|+|t-s|^{\frac{1}{\sigma}}&\,\mbox{if }  \sigma=2+\beta-\gamma \ge 2 ,\\
|x-y|+|t-s|^{\frac{1}{2}}&\,\mbox{if }  0<\sigma < 2\,.\\
\end{array}\right.
\end{equation}

The following first result generalizes the Harnack inequality of Moser \cite{MoserCpam71}, in the spirit of Chiarenza-Serapioni \cite{CS-RSMUP85,CS-AMPA84,CS-CPDE84} and Gutierrez-Wheeden \cite{GW,GW2}:
\begin{thm}[Parabolic Harnack inequality in the linear case]\label{thm.harnack.lin}
Let $v$ be a nonnegative bounded local weak solution to equation \eqref{WHE.DW.GEN} on $Q:=(0,T)\times\Omega$,  under assumption  \eqref{WHE.DW.GEN.paramt}. Then, for all $Q_{R}(t_0,x_0)\subset Q$, there exists $\ka_\ell>0$ such that
\begin{equation}\label{thm.harnack.lin.ineq}
\sup_{Q^-_R(t_0,x_0)} v\le \ka_\ell^{\,\lambda_0^{-1}+\lambda_1} \inf_{Q^+_R(t_0,x_0)} v.
\end{equation}
The constant $\ka_\ell>0$ depends on $N,\gamma, \beta$, but not on $v$ nor on $\lambda_0, \lambda_1$.
\end{thm}
\noindent\textbf{Remark. }As remarked before, although this result has been proven before at least in some range of parameters, the dependence of the Harnack constant on the ellipticity constants $\lambda_0,\lambda_1$ was not clear nor explicit; such dependence is needed in the proof of H\"older continuity for nonlinear equations , as we will show at the end of this section; this was pointed out by Moser in \cite{MoserCpam71}, where a complete proof of \eqref{thm.harnack.lin.ineq} in the unweighted case $\beta=\gamma=0$ can be found. The (nontrivial) fact that $\ka_\ell$ only depends on $N,\gamma, \beta$ is also pointed out by Gutierrez and Wheeden in \cite{GW2} after the statement of their Harnack inequalities, Theorem A; indeed we sketch here an adaptation of their proof to our case.

\medskip

As it often happens for linear parabolic equations, H\"older continuity follows by Harnack inequalities using a nowadays standard argument, cf. \cite{Moser}, that we sketch in the proof of Proposition \ref{Prop.Lin.Harn.HoCont}. We will assume in what follows, without loss of generality, that $\ka_\ell\ge 2$ and $0<\lambda_0\le  \lambda_1$,  and we define
\begin{equation}\label{Holder.exponent.lin}
\alpha:= \log_A\frac{\ka_\ell^{\,\lambda_0^{-1}+\lambda_1}}{\ka_\ell^{\,\lambda_0^{-1}+\lambda_1}-1} \in (0,1)\,,
\end{equation}
where $A>4$ which depends on $\gamma, \beta, N$ and is given in \eqref{expression.factor}. As well we introduce the notion of distance between sets of the form $Q=(0, T)\times \Omega \subset (0, \infty) \times \RR^N$. Let $Q'=(T_1, T_2) \times \Omega' \subset Q$, we define
\begin{equation}\label{parabolic.distance.sets}
d_{\gamma, \beta}(Q, Q'):= \inf\limits_{\substack{(t,x)\in \{[0,T]\times\partial\Omega\}\cup\{\{0\}\times\Omega\},\\ (s,y) \in Q'}}\, |x-y|\vee \left( \rho^{\gamma, \beta}_y\right)^{-1}(|t-s|)\,.
\end{equation}
We observe that if the $d_{\gamma, \beta}(Q, Q')=2D$ then for any $(t,x)\in Q'$ the parabolic cylinder $Q_{D}(t,x)\subset Q$.
\begin{prop}[H\"older Continuity in the linear case]\label{Prop.Lin.Harn.HoCont}
Let $v$ be a nonnegative bounded local weak solution to equation \eqref{WHE.DW.GEN} on $Q:=(0,T)\times\Omega$, under the assumption \eqref{WHE.DW.GEN.paramt}. Let $Q':=(T_1, T_2) \times \Omega' \subset Q$ and let $2D=d_{\gamma, \beta}(Q, Q')$. Then there exist $\alpha\in (0,1) $ as in \eqref{Holder.exponent.lin} and $\ka_\alpha >0$\,, such that for all $(t,x), (s,y) \in Q' $
\begin{equation}\label{Prop.Lin.Harn.HoCont.ineq}
\sup_{(t,x), (\tau,y)\in Q'}\frac{|v(t,x)-v(\tau,y)|}{d_{\gamma, \beta}\left((t,x), (s,y)\right)^\alpha}\le\frac{\ka_\alpha}{D^\alpha} \|v\|_{\LL^\infty(Q)},
\end{equation}
where $\ka_\alpha>0$ depends only on $N, \gamma, \beta, \lambda_0, \lambda_1$.
\end{prop}
The following corollary is immediate, and shows how the above estimates imply a more uniform modulus of continuity.
\begin{cor}\label{Cor.Lin.Harn.HoCont}
Under the assumptions of Proposition \ref{Prop.Lin.Harn.HoCont}, there exist $\alpha\in (0,1) $ as in \eqref{Holder.exponent.lin} and $\ka_\alpha' >0$\,, such that for all $(t,x), (s,y) \in Q' $
\begin{equation}\label{Prop.Lin.Harn.HoCont.ineq.111}
\sup_{(t,x), (\tau,y)\in Q'}\frac{|v(t,x)-v(\tau,y)|}{(|x-y|+|t-s|^{\frac{1}{2\vee \sigma}})^\alpha}\le\frac{\ka_\alpha'}{D^\alpha} \|v\|_{\LL^\infty(Q)},
\end{equation}
where $\ka_\alpha'>0$ is given by
\begin{equation*}
\ka_\alpha'= \ka_\alpha \ka_{19}^{\alpha}\left\{
\begin{array}{lll}
1\,,&\,\mbox{if }  \sigma \ge 2 ,\\
\left(T^{\frac{1}{\sigma}}\vee \sup\limits_{x_0\in\Omega}|x_0|\right)^{\frac{\gamma-\beta}{2}}\,,&\,\mbox{if }  0<\sigma < 2\,.\\
\end{array}\right.
\end{equation*}
Where $\ka_{19}>0$ depends only on $N, \gamma, \beta$ and is given in \eqref{inverse.technical.inequality.measures}.
\end{cor}

The proof of the above results relies on the following upper and lower bounds.
\begin{prop}\label{spacetime.smoothing.theorem.linear}
    Let $ u\in \LL^p_{\gamma, \rm loc}(\left(0, T \right)\times B_R(x_0))$ with $p>0$ be a nonnegative local strong (sub)solution to \eqref{WHE.DW.GEN}  and let $x_0 \in \RR^N$, $ 0 < R_1 < R_0 < R $ such that $0 \not \in \overline{B_{R_0}(x_0)\setminus B_{R_1}(x_0)}$ and let $0 \leq T_0 < T_1 < T$. Then there exists a constant $\ka_{lin}>0$ depending only on $\gamma, \beta, N, p$ such that the following inequality holds\vspace{-1mm}
    \begin{equation}\label{upper.bounds.lin} \begin{split}
        \sup\limits_{(\tau,y)\in (T_1, T] \times B_{R_1}(x_0)} u(\tau,y) &\leq  \ka_{lin}
        \left[\frac{h_{\sigma}\left( R_0, R_{1}, x_0 \right)}{\left(R_0 - R_{1}\right)^{\sigma}} + \frac{1}{T_{1} - T_0} \right]^{\frac{N-\gamma +\sigma}{\sigma p}}
\left[ \int_{T_0}^T\int_{B_{R_0}(x_0)}{u^{p} \frac{\dx\dt}{|x|^{\gamma}}  }\right]^{\frac{1}{p}}
    \end{split}\end{equation}
where $\sigma$ is defined in \eqref{definition.sigma.thethap},
$h_\sigma(R_0, R_1, x_0)$ is defined in \eqref{function_h} and $\ka_{lin}$ is a computable constant  such that
$\ka_{lin}\lesssim  S_{\gamma, \beta}^{\frac{2(N-\gamma)}{p\sigma}} \left(\lambda_0^{-1}\lambda_1\right)^{\frac{(N-\gamma + \sigma)}{\sigma p}}$,
with $S_{\gamma, \beta}$  as in Proposition \ref{prop.Sobolev.Balls}.
\end{prop}
\begin{prop}\label{cor.lin.lower}
Let $ u $ be a nonnegative local strong (super)solution to \eqref{WHE.DW.GEN} on $(0, T)\times B_R(x_0)$, with $ 0 < R_1 < R_0 < R $ such that  $0 \not \in \overline{B_{R_0}(x_0)\setminus B_{R_1}(x_0)}$ and let $0 \leq T_0 < T_1 < T$. Then for any $p>0$ there exists a constant $\kb_{lin}>0$ depending only on $\gamma, \beta, N, p$ such that the following inequality holds\vspace{-1mm}
    \begin{equation}\label{lower.bounds.lin}  \begin{split}
        \inf\limits_{(\tau,y)\in (T_1, T] \times B_{R_1}(x_0)} u(\tau,y) &\geq  \kb_{lin}
        \left[\frac{h_{\sigma}\left( R_0, R_{1}, x_0 \right)}{\left(R_0 - R_{1}\right)^{\sigma}} + \frac{1}{T_{1} - T_0} \right]^{-\frac{N-\gamma +\sigma}{\sigma p}}
\left[ \int_{T_0}^T\int_{B_{R_0}(x_0)}{u^{-p} \frac{\dx\dt}{|x|^{\gamma}}  }\right]^{-\frac{1}{p}}
    \end{split}\end{equation}
where $\sigma$ is defined in \eqref{definition.sigma.thethap},
$h_\sigma$ is defined in \eqref{function_h} and $\kb_{lin}$ is a computable constant  such that
$\kb_{lin}\gtrsim  S_{\gamma, \beta}^{-\frac{2(N-\gamma)}{p\sigma}} \left(\lambda_0^{-1}\lambda_1\right)^{-\frac{(N-\gamma + \sigma)}{\sigma p}}$
where $S_{\gamma, \beta}$ is as in Proposition \ref{prop.Sobolev.Balls}.
\end{prop}
\noindent\textbf{Remark. }The above estimates have been previously obtained by several authors in different settings, we just mention here the closest results: Lemma 3.17 of \cite{GW2} (in the case of general $\mathcal{A}_2$ weights), Lemma 2.1 of \cite{CS-RSMUP85}, and Lemma 1 of \cite{MoserCpam71} when there are no weights. The proof follows Moser's idea: using weighted Sobolev inequalities and upper (resp. lower)  iterations , to obtain upper (resp. lower) space-time smoothing effects;  indeed, the space-time upper bounds \eqref{upper.bounds.lin} can be obtained also by taking $m=1$ in the proof of Theorem \ref{spacetime.smoothing.theorem}: note that the two factors $u^{p+m-1}$ and $u^p$ in the energy estimates \eqref{sup.energy.inequality.upper} are now the same, hence the proof can be done directly with $u$, and we do not need to use the subsolution $v= u \vee 1$; as a consequence, the factor $+1$ in the integral in the right-hand side of formula \eqref{spacetime.smoothing} disappears.  Analogously, the space-time lower bounds \eqref{lower.bounds.lin} follow by a minor modification of the proof of Proposition \ref{lower.spacetime.smoothing} with $m=1$, more precisely we just repeat the Steps 1, 2 and 3 of the proof and we obtain the analogous of formula \eqref{lower.first.inequality.9}, which can be rewritten in the form \eqref{lower.bounds.lin} .    Note that these proofs are considerably simpler than in the nonlinear setting, $m\in (0,1)$.

\medskip

\noindent {\bf Proof of Harnack inequalities, Theorem \ref{thm.harnack.lin}.~}The proof  follows the lines of the original Moser proof in \cite{MoserCpam71}. Once obtained local upper and lower bounds, \eqref{upper.bounds.lin} and \eqref{lower.bounds.lin}, the   hardest  part of the proof consists in obtaining a reverse H\"older inequality that allows one to join them and deduce the Harnack inequalities   \eqref{thm.harnack.lin.ineq}.  To our knowledge only two techniques are known to perform this task: one  originally due to Moser \cite{Moser} that gives a suitable reverse space-time H\"older inequality on shifted cylinders, and another due to Bombieri and Giusti \cite{BGInv72}, see also \cite{MoserCpam71}, which shows how estimates \eqref{upper.bounds.lin} and \eqref{lower.bounds.lin} imply (local) absolute upper and lower bounds that allow to obtain the Harnack inequality \eqref{thm.harnack.lin.ineq}. We will follow the latter strategy, and just sketch the proof, which is essentially the same as Section 4 and 5 of \cite{GW2}, see also \cite{MoserCpam71}; we shall focus on the points where we some non-straightforward changes are needed. Last, we just remark that it is enough to work in a cube of size 1, then the result will follow by rescaling.

\noindent$\bullet~$\textsc{Step 1. }\textit{Bombieri-Giusti Lemma. }We are going to use a weighed version of Bombieri-Giusti Lemma as it has been done in Section 5 of \cite{GW2} (see also Lemma 3 of \cite{MoserCpam71}). Notice that the following Lemma applies to measurable functions $f$, not necessarily solutions to a PDE.\\
\textsl{Claim}. Let $A,B, \overline{p}, \varrho,\delta$ be positive constants,   and $Q_1, Q_{R_0}, Q_{R_1}, Q_{ \varrho}$ as in \eqref{Parabolic.cylinders.lin}. Also,  we assume that the positive measurable function $f$ defined on $Q_1$,  and the doubling measure $\nu$ on $\RR^{N+1}$ satisfy the inequalities\vspace{-2mm}
\begin{equation}\label{BoGiu.hyp}
\sup_{Q_{R_1}}f^p\le \frac{A}{(R_0-R_1)^B}\iint_{Q_{R_0}}f^p\nu(t,x)\dx\dt
\qquad\mbox{and}\qquad
\nu\left\{(t,x)\in Q_1\,:\, \log f > s\right\}\le \left(\frac{1}{s\,\overline{p}}\right)^\delta\nu(Q_1)\vspace{-2mm}
\end{equation}
for all $s>0$, $\frac{1}{2}\le \varrho \le R_1<R_0<1$, all $p\in (0, \overline{p})$. Then there exists $c_0=c_0(A,B,\delta)>0$ such that\vspace{-2mm}
\begin{equation}\label{BoGiu}
\log\sup_{Q_\varrho }f\le \frac{c_0}{\overline{p}(1-\varrho)^{2B}}.\vspace{-2mm}
\end{equation}
The proof of the above claim   is a minor modification of the proof of Lemma 5.1 of \cite{GW2}, see Section 5 of \cite{GW2} for more details. Indeed, for some range of parameters, for instance when our weights fall in the Muckenhoupt class $\mathcal{A}_2$, the proof is exactly the same. The only point where we can not directly adapt those proofs, is when a suitable ``localized'' weighted Poincar\'e inequality is used: in our context, such inequality reads\vspace{-2mm}
\begin{equation}\label{localized.poincare}
\int_{B_R(y_0)}|f(x)-\overline{f}|^2\frac{\varphi(x)}{|x|^{\gamma}}\dx \le c_\varphi \frac{\mu_\gamma(B_R(y_0))}{\mu_\beta(B_R(y_0))} R^2  \int_{B_R(y_0)} |\nabla f|^2\frac{\varphi(x)}{|x|^{\beta}}\dx\,,\vspace{-2mm}
\end{equation}
where $\overline{f}=\left(\int_{B_R(y_0)} \frac{\varphi(x)}{|x|^{\gamma}}\dx\right)^{-1}\left(\int_{B_R(y_0)}f(x) \frac{\varphi(x)}{|x|^{\gamma}}\dx\right)$, for on any ball $B_R(y_0)\subset\RR^N$ and for an extra ``weight''    $\varphi\in C_0(B_R(y_0))$, $0\le \varphi\le 1$ with convex super-level sets, where $c_\varphi=c_{N,\gamma,\beta}\left(|B_R(y_0)|/\int_{B_R(y_0)}\varphi\dx\right)^2$. This inequality is proven in Lemma 4.1 of both \cite{GW2,GW}, and the relies on results of \cite{SW} (or, in the non-weighted case see Lemma 3 of \cite{Moser}); this is the point where the restriction on the class of weights appears. We recall that the results of \cite{GW2,GW} do not cover all the range of parameters $\gamma, \beta$ that we consider here:  they hold  for weights which satisfy  the $\mathcal{A}_2$ property (or generalizations of it), and this is not always the case in our setting.   A closer inspection of the proof of Lemma 3 of \cite{Moser} reveals that it is enough to prove \eqref{localized.poincare} just for one function $\varphi$  with the properties that $0\le \varphi\le 1$ on $B_R(y_0)$\,, and for some $\delta\in (0,1)$ and some $R_\delta\in (0,R)$ we also have that $\varphi\ge \delta$ on $B_{R_\delta}(y_0)\subset B_R(y_0)$ and $\varphi=0$ on $\partial B_R(y_0)$. We are going to show that inequality \eqref{localized.poincare} is indeed a consequence of the so-called Intrinsic Poincar\'e inequality\vspace{-2mm}
\begin{equation}\label{intrinsic.Poincare}
(\lambda_2-\lambda_1)\int_{B_R(y_0)}|f(x)-\overline{f}|^2\frac{\varphi^2_1(x)}{|x|^{\gamma}}\dx \le C_{\gamma,\beta}  \int_{B_R(y_0)} |\nabla f|^2\frac{\varphi^2_1(x)}{|x|^{\beta}}\dx,\vspace{-2mm}
\end{equation}
where $\varphi_1$ is the first eigenfunction of the operator $\mathcal{L}_{\gamma, \beta}$ (with Dirichlet boundary  conditions and with unitary $\LL^2_\gamma$ norm) on $B_R(y_0)$,
$\overline{f}=\left(\int_{B_R(y_0)}f(x) \frac{\varphi_1^2(x)}{|x|^{\gamma}}\dx\right)$, and $\lambda_1, \lambda_2$ are respectively the first and the second eigenvalue of the  $\mathcal{L}_{\gamma, \beta}$ on $B_R(y_0)$. The proof of inequality \eqref{intrinsic.Poincare} is quite standard: this inequality is indeed equivalent to the second Poincar\'e inequality\vspace{-2mm}
\[
\lambda_2 \|g\|^2_{\LL^2_\gamma(B)}\le \|\nabla g\|^2_{\LL^2_\beta(B)}\qquad\mbox{for all $g\in \D_{\gamma, \beta}\left( B \right)$ such that}\qquad
\int_\Omega g\varphi_1\frac{\dx}{|x|^{\gamma}}=0\,.\vspace{-2mm}
\]
The above inequality is true on balls as a consequence of the compactness of the embedding $\D_{\gamma, \beta}\left( B \right)\subset \LL^2_\gamma(B)$, where $\D_{\gamma, \beta}$ is defined in Subsection \ref{result.organization}; finally, inequality \eqref{intrinsic.Poincare} follows by letting $g=(u-\overline{u})\varphi_1$, see for instance Lemma 3.1 of \cite{BGV-JMPA}. The claim is proven.\\

\noindent$\bullet~$\textsc{Step 2. }\textit{Proof of the Harnack inequality. }We are going to apply twice the result of the previous step to get local upper and lower bounds that will finally combine in the Harnack estimates \eqref{thm.harnack.lin.ineq}. In this case the proof is an adaptation of Section 6 of \cite{GW2}, see also Section 3 of \cite{MoserCpam71} for the non-weighted case, we will just emphasize the essential changes.

Let us assume first that we are in the position to apply the Bombieri-Giusti inequality with $R_0=3/4$ and $ R_1=R=2/3$  and $\varrho=1/2$,  to both $f\sim v$ and $f\sim v^{-1}$ on $Q^-_R(t_0,x_0)$ and $Q^+_R(t_0,x_0)$ respectively; we will briefly explain at the end of this step how to proceed to ensure that the assumptions \eqref{BoGiu.hyp} are satisfied by both $v$ and $v^{-1}$. Using inequality \eqref{BoGiu} for $f=\ee^{-M_2+V}v$ on $Q^-_R(t_0,x_0)$, where $V$ and $M_2$ are chosen as in \eqref{log.est}, we obtain the desired absolute local upper bounds:\vspace{-2mm}
\begin{equation}\label{BoGiu.upper}
\sup_{Q^-_{1/2}(t_0,x_0)}v=\ee^{M_2-V}\sup_{Q^-_{1/2}(t_0,x_0)}f\le \ee^{M_2-V}\exp\left(\frac{c_0}{\overline{p}(1/2)^{2B}}\right)\le \ka_{\ell,1}^{\,\lambda_0^{-1}+\lambda_1}\ee^{-V}.\vspace{-2mm}
\end{equation}
Proceeding analogously for $f=\ee^{-M_2-V}v^{-1}$ on $Q^+_R(t_0,x_0)$, we obtain the desired absolute local lower bounds:\vspace{-1mm}
\begin{equation}\label{BoGiu.lower}
\inf_{Q^+_{1/2}(t_0,x_0)}v=\ee^{-M_2-V}\left(\sup_{Q^+_{1/2}(t_0,x_0)}f\right)^{-1}
\ge \ee^{-M_2-V}\exp\left(-\frac{c_0}{\overline{p}(1/2)^{2B}}\right)
\ge\ka_{\ell,2}^{-\,\lambda_0^{-1}-\lambda_1}\ee^{-V}.\vspace{-1mm}
\end{equation}
Notice that the last inequalities in \eqref{BoGiu.upper} and \eqref{BoGiu.lower} follow by the choice of $V, M_2$ as in \eqref{log.est} and can be proven by following exactly Section 6 of \cite{GW2}, hence we omit the details. Finally, the Harnack inequality \eqref{thm.harnack.lin.ineq} follows by combining inequalities \eqref{BoGiu.upper} and \eqref{BoGiu.lower} and $\ka_{\ell}=\ka_{\ell,1}\cdot \ka_{\ell,2}$.

It only remain to show that we can actually use inequality \eqref{BoGiu} for $f=u$ and $f=u^{-1}$, hence we need to ensure the validity of hypotheses \eqref{BoGiu.hyp} in both cases. This is done by proving the so-called logarithmic estimates, see for instance  Lemma 4.9 of \cite{GW2}. The proof of that Lemma can be repeated also in our setting, and shows that: for any nonnegative bounded solution $u$ defined on $(a,b)\times B_{3/2}$, bounded below by a positive constant in $(a,b)\times B_{1}$, then there exist $c_1, M_2,\delta$ and $V$ such that, for any $s>0$\vspace{-1mm}
\begin{equation}\label{log.est}\begin{split}\vspace{-1mm}
\mu_\gamma\left\{(t,x)\in (t_0,b)\times B_1(x_0)\,:\, \log u<-s-M_2(b-t_0)-V\right\}
    &\le c_1\left[\frac{1}{s}\frac{\mu_\gamma(B_1(x_0))}{\mu_\beta(B_1(x_0))}\frac{1}{b-t_0}\right]^\delta(b-t_0),\\
\mu_\gamma\left\{(t,x)\in (a,t_0)\times B_1(x_0)\,:\, \log u>s-M_2(a-t_0)-V\right\}
&\le c_1\left[\frac{1}{s}\frac{\mu_\gamma(B_1(x_0))}{\mu_\beta(B_1(x_0))}\frac{1}{t_0-a}\right]^\delta(t_0-a),
\end{split}
\end{equation}
where the constants $c_1,\delta>0$ only depend on $N,\beta,\gamma$, $M_2\sim \mu_\beta(B_1(x_0)) /\mu_\gamma(B_1(x_0))$, and $V$ depends on $v$, but it is the same in both cases, as explained carefully in Section 6 of of \cite{GW2}. Details about the proof of the above estimates can be found in Section 4 of of \cite{GW2}, which in turn extend ideas of Moser (Section 2 and 3 of \cite{MoserCpam71}) to the weighted case. The latter estimates, together with the local upper and lower bounds, \eqref{upper.bounds.lin} and \eqref{lower.bounds.lin}, allow to apply the Bombieri-Giusti result in both cases. Hence \eqref{BoGiu.upper} and \eqref{BoGiu.lower} hold and the proof is concluded. The general statement follows by a scaling argument.

Finally we recall that,  as Moser first noticed in \cite{MoserCpam71}, with the present method it is possible to keep track of the dependence on the ``ellipticity'' constants $0<\lambda_0 \le  \lambda_1<+\infty$ throughout the proof: also in the present weighted setting we were able to keep track of such dependence in the constants.\qed

\medskip
\noindent {\bf Proof of H\"older continuity, Proposition \ref{Prop.Lin.Harn.HoCont}.~}We adapt the Moser's proof to our weighted case, see \cite{Moser,MoserCpam71}. We first prove how the oscillation of the solution decrease geometrically on parabolic cylinders.  We recall that if the $d_{\gamma, \beta}(Q, Q')=2D$ then for any $(t,x)\in Q'$ the following inclusion holds $Q_{D}(t,x)\subset Q$. Fix $r \in (0,D/2)$ and denote for simplicity $Q_{r}:=Q_r(t_0,x_0)$ and $Q_r^{\pm}:=Q_r^{\pm}(t_0,x_0)$.  Let us introduce the following quantities:
\[
\M_r :=\sup_{Q_r} v\,, \quad \M_r^{\pm}:=\sup_{Q^{\pm}_r} v\,, \quad \m_r:=\inf_{Q_r} v\,, \quad \m_r^{\pm}:=\inf_{Q^{\pm}_r}v\,.
\]
We are in the position to apply the Harnack inequality \eqref{thm.harnack.lin.ineq} to the nonnegative solution $\M_{2r}-v$ to obtain
\[
\M_{2r}-\m_{r}^-=\sup_{Q_{r}^-}(\M_{2r}-v) \le H\,\inf_{Q_{r}^+}(\M_{2r}-v)=H\,(\M_{2r}-\M^+_{r})\,.
\]
Notice that without loss of generality we can set $H:=\ka_\ell^{\,\lambda_0^{-1}+\lambda_1}\ge 4$\,, with $\ka_\ell$ as in \eqref{thm.harnack.lin.ineq}.
Similarly, using $v-\m_{2r}$ we obtain  $\M_{r}^- - \m_{2r}\le H\,(\m^+_{r}-\m_{2r})$ which, summed up with the previous inequality, gives
\[
H\,(\M^+_{r}-\m^+_{r}) + \M_{r}^- - \m_{r}^-\le (H-1)\,(\M_{2r}-\m_{2r})\,.
\]
Using $Q_{r/A}\subset Q^+_r$ (see Lemma \ref{multiplicative.factor}, formula \ref{inclusion.parabolic.cylinders}), we conclude that
\[
\osc_{Q_{r/A}}v\le \osc_{Q^+_r}v=\M^+_{r}-\m^+_{r} \le \frac{H-1}{H}\,(\M_{2r}-\m_{2r})=\frac{H-1}{H}\,\osc_{Q_{2r}}v\,.
\]
Recall that without loss of generality  we have assumed that  $H/(H-1)\le 4<A$\,,   see also \eqref{Holder.exponent.lin}; a well-known iteration technique (see, \emph{e.g.},~\cite[Lemma~6.1]{Giusti}) then shows that
\begin{equation}\label{oscillation.reduction}
\osc_{Q_{r}}v\le A^{\alpha}\,\frac{r^{\alpha}}{D^\alpha}\,\osc_{Q_{R}}v\qquad\qquad\mbox{for all }\,r\in(0,D]\,,
\end{equation}
with $\alpha:=\log({H}/(H-1))/\log A\in(0,1)$, as in \eqref{Holder.exponent.lin}, and $H>1$ as above.

Now, we fix $(t,x),(s,y) \in Q'$, and we first assume that $\Omega'$ is convex. The first case that we analyze corresponds to $d_{\gamma, \beta}((t,x),(s,y))\leq D$. Hence, there exists an integer $k \geq 1$ such that
\[
\frac{D}{A^{k+1}}\leq d_{\gamma, \beta}((t,x),(s,y))<\frac{D}{A^{k}},
\]
from which follows that $(t,x),(s,y) \in Q_{\frac{D}{A^{k}}}(t\vee s,\overline{xy})\subseteq Q_{D}(t\vee s,\overline{xy}) \subset Q$, where $\overline{xy}=(x+y)/2$. Using \ref{oscillation.reduction} we get the following estimate
\begin{equation}
|v(t,s)-v(s,y)|\leq \osc_{Q_{D/A^{k}(\overline{xy}, t\vee s)}}v  \leq \left(\frac{A\, D}{A^{k}D }\right)^\alpha \|v\|_{\LL^{\infty}(Q)}\leq \frac{A^{2\alpha}}{D^\alpha}d_{\gamma, \beta}((t,x),(s,y))^\alpha\,\|v\|_{\LL^{\infty}(Q)}\,.
\end{equation}
The second case corresponds to $d_{\gamma, \beta}((t,x),(s,y))> D$, and we proceed as follows
\[
|v(t,s)-v(s,y)|\leq 2 \|v\|_{\LL^{\infty}(Q)} \leq \frac{2 \|v\|_{\LL^{\infty}(Q)}}{D^\alpha}d_{\gamma, \beta}((t,x),(s,y))^\alpha\,.
\]
The constant $\ka_\alpha>0 $ is given by
\begin{equation}\label{constant.alha}
\ka_\alpha:= 1 \vee A^\alpha\,,
\end{equation}
where $\alpha$ is as in \eqref{Holder.exponent.lin} and $A$ as in \eqref{multiplicative.factor}; it depends on $N, \gamma, \beta$ and $\lambda_0,\lambda_1$. In the case when $\Omega'$ is not convex, the result follows by a standard covering argument, however for the purposes of the present work, we only need quantitative information on balls. The proof is now concluded.\qed
\medskip
\noindent {\bf Proof of Corollary \ref{Cor.Lin.Harn.HoCont}.~}As a consequence of inequality \eqref{bound.inverse.rho} of Lemma \ref{inverse.technical.inequality.measures.lemma}, we know that there exist a constants $\ka'>0$ such that for any $(t,x), (s,y) \in Q$ we have
\begin{equation*}
d_{\gamma, \beta}((t,x),(s,y))\leq \ka'\left( |x-y|+|t-s|^{\frac{1}{2\vee\sigma}}\right),
\end{equation*}
which proves the Corollary.\qed

\subsection{The nonlinear case}
This Subsection essentially contains the proofs of Harnack inequalities and H\"older continuity for WFDE, Theorems \ref{Harnack.2} and \ref{thm.holder.nonlin} respectively.
\begin{thm}[Alternative form of Harnack inequality]
Under the assumptions of Theorem \ref{Harnack.2}, for any $t_0 > 0$ there exist constants $\ka_8, \ka'_9,\kappa_*>0$ such that
  \begin{equation*}
  \sup_{x \in B_R(x_0)} u(t, x) \le \ka_8 \frac{\|u(t_0)\|_{\LL^p_{\gamma}(B_{2R}(x_0))}^{p\sigma \vartheta_p}}{t_0^{(N-\gamma)\vartheta_p}} + \ka'_9 \inf_{x \in B_R(x_0)} u(t \pm \theta, x)
  \end{equation*}
for any
  \begin{equation*}\label{t.star.haranck.alt}
    t, t\pm \theta \in (t_0, t_0+t_{\ast}(t_0)) \cap (0, T), \qquad \mbox{and} \qquad  t_\ast(t_0)=\kappa_* R^{\sigma} \mu_\gamma(B_{R}(x_0))^{m-1} \|u(t_0)\|_{\LL^1_{B_{R}(x_0)}}^{1-m}.
  \end{equation*}
The constants $\ka_8, \ka'_9,\kappa_*>0$ depend on $N,m,\gamma,\beta$; $\kappa_*>0$ is given in the proof of Corollary \ref{weak.positivity.MDP.L1},
 and $\ka'_9=\ka_9 \kb^{-1}$  where $\ka_8,\ka_9>0$ are as in \eqref{local_upper_bounds}, while $\kb>0$ has an (almost) explicit expression is given in \eqref{const.lower.cor.L1};   note  that $\kb, \ka_8$ and $\ka_9$   depend on $R$ and $x_0$ and, when  $0<m\le m_c$, $\kb$ depends also on $ H_p(u_0, x_0, R)$ defined in \eqref{def.HP}.
\end{thm}
\noindent\textbf{Proof. }It follows immediately by combining inequalities \eqref{upper.provisional.2} and \eqref{LOCAL.LOWER.INEQUALITY}. \qed

 \medskip

\noindent\textbf{Proof of the Harnack inequalities of Theorem \ref{Harnack.2}. }Due to the time translation invariance of the equation it suffices to prove the result for $t_0=0$. Assume $t \in (\varepsilon t_*, t_*)$, for $\varepsilon\in (0,1)$ fixed. Recall that $R^{\sigma}\asymp R^2 \mu_\gamma(B_R(x_0))\mu_\beta(B_R(x_0))^{-1}$.  Using the upper bound \eqref{upper.bounds.rho}, inequality \eqref{r.sigma} and formula \eqref{def.HP} we get\vspace{-2mm}
\begin{equation}\label{harnack.estimate3bbb}\begin{split}
&\sup_{x \in B_R(x_0)} u(t, x) \le \ka_1 \frac{\|u_0\|_{\LL^p_\gamma(B_{2R}(x_0))}^{p\sigma\vartheta_p}}{t^{(N-\gamma)\vartheta_p}}+\ka_2\left[\frac{t}{\ka_{17} R^\sigma} \right]^{\frac{1}{1-m}}
\leq\left[\frac{\ka_1\,\|u_0\|_{\LL^p_\gamma(B_{2R}(x_0))}^{p\sigma\vartheta_p}R^{\frac{\sigma}{1-m}}}{(\varepsilon t_{\ast})^{(N-\gamma)\vartheta_p+\frac{1}{1-m}} } + \frac{\ka_2}{\ka_{17}^{\frac{1}{1-m}}} \right]\left[\frac{t}{R^\sigma}\right]^{\frac{1}{1-m}}\\
&\le \left[\frac{\ka_1 \omega_\gamma^{\sigma \vartheta_p}}{2^{\frac{\sigma}{1-m}}}\frac{H_p(u_0, x_0, 2R)}{\varepsilon^{\frac{p\sigma\vartheta_p}{1-m}}} + \frac{\ka_2}{\ka_{17}^{\frac{1}{1-m}}} \right]\left[\frac{t}{R^\sigma}\right]^{\frac{1}{1-m}} \leq \left[ \frac{\ka_1 \omega_\gamma^{\sigma \vartheta_p}}{2^{\frac{\sigma}{1-m}}}+ \frac{\ka_2}{\ka_{17}^{\frac{1}{1-m}}} \right]\left[1 \vee\frac{H_p(u_0, x_0, 2R)}{\varepsilon^{\frac{\sigma p \vartheta_p}{1-m}}}\right]  \left[\frac{t}{R^\sigma}\right]^{\frac{1}{1-m}},\vspace{-2mm}
\end{split}\end{equation}
where $\omega_\gamma=B_1(0)$ and $\ka_{17}$ as in inequality \eqref{r.sigma}.  We recall next the lower bound \eqref{LOCAL.LOWER.INEQUALITY}\,, that in this case reads\vspace{-2mm}
\begin{equation*}\label{harnack.estimate4}
    \inf_{x \in B_{R}(x_0)} u(t,x) \ge\inf_{x \in B_{2R}(x_0)} u(t,x)\ge \kb \left[\frac{t}{(2R)^{\sigma}}\right]^{\frac{1}{1-m}}\qquad\mbox{for any $t \in[0, t_*] \cap (0, T)$.}\vspace{-2mm}
\end{equation*}
 By  combining the two above inequalities we get for any $t \in (\varepsilon t_*, t_*)$:\vspace{-2mm}
\begin{equation*}\label{harnack.estimate5}\begin{split}
\sup_{x \in B_R(x_0)} u(t, x)
\le  \left[ \frac{\ka_1 \omega_\gamma^{\sigma \vartheta_p}}{\kb}+ \frac{\ka_22^{\frac{\sigma}{1-m}}}{\kb\,\ka_{17}^{\frac{1}{1-m}}} \right]\left[1 \vee\frac{H_p(u_0, x_0, 2R)}{\varepsilon^{\frac{\sigma p \vartheta_p}{1-m}}}\right]
\inf_{x \in B_{R}(x_0)} u(t,x):= \ka_3 \inf_{x \in B_{R}(x_0)} u(t,x)\vspace{-2mm}
\end{split}\end{equation*}
The constants $ \ka_1, \ka_2>0$ depend on $N,m,\gamma,\beta$ and are given in  \eqref{upper.bounds.rho};  $\kb>0$ is given in \eqref{const.lower.cor.L1}: notice that,  when $0<m\le m_c$,  $\kb$ depends also on $ H_p(u_0, x_0, 2R)$ defined in \eqref{def.HP}. We finally recall that\vspace{-2mm}
\begin{equation}\begin{split}\label{const.harnack}
\ka_3&:=  \left[ \frac{\ka_1 \omega_\gamma^{\sigma \vartheta_p}}{\kb}+ \frac{\ka_22^{\frac{\sigma}{1-m}}}{\kb\,\ka_{17}^{\frac{1}{1-m}}} \right]\left[1 \vee\frac{H_p(u_0, x_0, 2R)}{\varepsilon^{\frac{\sigma p\vartheta_p}{1-m}}}\right] \\
 & \asymp  \left[\ka_1 \omega_\gamma^{\sigma \vartheta_p}+ \frac{\ka_22^{\frac{\sigma}{1-m}}}{\ka_{17}^{\frac{1}{1-m}}} \right] \frac{H_p(u_0, x_0, 2R)}{\varepsilon^{\frac{\sigma p \vartheta_p}{1-m}}} \left(\frac{\widetilde{H}_p^{c_2}}{R^{c_1}}\right)^{\frac{\widetilde{H}_p^{1/2}}{m(1-m)}}\,\,\,\, \mbox{when}\,\,\,\,\, \widetilde{H}_p\gg 1\,, \vspace{-2mm}
\end{split}\end{equation}
where $\widetilde{H}_p(u_0, x_0, 2R):= 1+\left(\frac{|x_0|}{2R}\vee 1\right)^{\beta-\gamma} H_p(u_0, x_0, 2R)^{1-m} \ge 1$
and  $c_1,c_2>0$ only depend on $N,m,p,\beta,\gamma$. See also Corollary \ref{weak.positivity.MDP.L1} for a more detailed the expression of $\kb, c_1,c_2$. This concludes the proof. \qed
\medskip

We can prove an analogous continuity result for local solutions to the (WFDE), using the upper and lower bounds of Theorems \ref{local_upper_bounds} and \ref{LOCAL.LOWER.BOUNDS}, and the  the linear results of the previous subsection.

\medskip

\noindent {\bf Proof of the H\"older continuity estimate of Theorem \ref{thm.holder.nonlin}.~}We split the proof in two steps.

\noindent$\bullet~$\textsc{Step 1. }\textit{Intrinsic rescaling. }We begin by considering a local solution $u$ on the cylinder $Q:=(0,T]\times \Omega$. Fix $t_0, R_0>0$ such that $Q^*_{4R_0}(t_0, x_0):=\left(t_0, T\wedge(t_0+ t_*) \right]\times B_{4R_0}(x_0)\subset Q$.
We define the rescaled solution $\hat{u}$ as follows:
\[
\hat{u}(\hat{t},\hat{x}):=M_0^{-1}u(t,x)\qquad\mbox{with}\qquad t=R_0^\sigma M_0^{1-m}\hat{t}\,,\; x=R_0\hat{x},
\]
 where $M_0$ is any positive real number such that $M_0\ge \|u\|_{\LL^\infty(Q_{4R_0}(t_0,x_0))}$.  It is easy to check that if $u$ is a local solution on $Q^*_{4R_0}(t_0, x_0)$, then $\hat{u}$ is a local solution to the same equation on $Q^*_{4}(\hat{t}_0,\hat{x}_0):=[\hat{t}_0,\hat T\wedge(\hat{t}_0+\hat{t}_*)]\times B_{4}(\hat{x}_0)$, where
\begin{equation*}\label{tstar.hat}
     \hat{t}_*= \hat{t}_*(\hat{u}(\hat{t}_0), \hat{x}_0,4)= \kappa_*\,4^\sigma \frac{ \|\hat{u}(\hat{t}_0)\|_{\LL^1_{\gamma}(B_4(\hat{x}_0))}^{1-m}}{\mu_\gamma(B_4(\hat{x}_0))^{1-m}}
= \kappa_*\,4^\sigma  \left[\frac{ \|u(\hat{t}_0)\|_{\LL^1_{\gamma}(B_4(\hat{x}_0))}}{M_0\,\mu_\gamma(B_4(\hat{x}_0))}\right]^{1-m} \,,
\end{equation*}
Moreover, $\|\hat{u}\|_{\LL^\infty(Q^*_{4}(\hat{t}_0,\hat{x}_0))}\le 1$\,, since by assumption we have $M_0\ge \|u\|_{\LL^\infty(Q_{4R_0}(t_0,x_0))}$.\\
We are now in the position to apply the lower bounds of Theorem \ref{LOCAL.LOWER.BOUNDS} to $\hat{u}$ on $Q^*_{4}(\hat{t}_0,\hat{x}_0)$:
\begin{equation}\label{LOCAL.LOWER.INEQUALITY.Calpha}
       \inf_{x \in B_{2}(\hat{x}_0)} \hat{u}(t,x) \ge \kb \left[\frac{\hat{t}-\hat{t}_0}{2^{\sigma}}\right]^{\frac{1}{1-m}}\qquad\mbox{for any $\hat{t} \in[\hat{t}_0+\frac{1}{4}\hat{t}_*, \hat{t}_0+\hat{t}_*]\cap(0,\hat T)$},
    \end{equation}
 where $2$ is the radius of the ball.  Notice that $\kb$ has an (almost) explicit expression is given in \eqref{const.lower.cor.L1},  and  (in the very fast diffusion range, i.e. when $m<m_c$ and $p>1$) depends on $ H_p(\hat{u}(\hat{t}_0),\hat{x}_0, 4)$ defined in \eqref{def.HP}.
Clearly inequality \eqref{LOCAL.LOWER.INEQUALITY.Calpha} implies
\begin{equation*}\label{LOCAL.LOWER.INEQUALITY.Calpha.1}
  \inf_{(t,x) \in Q^*_{2}(\hat{t}_0,\hat{x}_0)} \hat{u}(t,x)
  =\inf_{(t,x) \in [\hat{t}_0+\hat{t}_*/2, \hat{t}_0+\hat{t}_* ]\cap(0,\hat T)\times B_{2}(\hat{x}_0)} \hat{u}(t,x)
  \ge 4^{-\frac{\sigma}{1-m}}\kb\, \hat{t}_*^{\frac{1}{1-m}}\,.
    \end{equation*}

\noindent$\bullet~$\textsc{Step 2. }\textit{Application of the linear result. }$\hat{u}$ can be considered a solution to the linear equation \eqref{WHE.DW.GEN} with $a(t,x)=m\hat{u}^{m-1}(t,x)$; we are now in the position to apply the result of Corollary \ref{Cor.Lin.Harn.HoCont}
inside the cylinder $Q^*_2:=[\hat{t}_0 + \hat{t}_*/2,  \hat{t}_0+\hat{t}_* ]\cap(0,\hat T)\times B_{2}(\hat{x}_0)$, since in $Q^*_2$ we have
\begin{equation}\label{lambda01bbb}
\lambda_0:= m\le m\hat{u}^{m-1}=a(t,x)\le 4^\sigma\kb^{m-1}\hat{t}_*^{-1}=:\lambda_1\,.
\end{equation}
 Then, on the cylinder $Q^*_1:=[\hat{t}_0 + (5/8)\hat{t}_*,  \hat{t}_0+(7/8)\hat{t}_* ]\cap(0,\hat T)\times B_{1}(\hat{x}_0)$, the result of Corollary \ref{Cor.Lin.Harn.HoCont} implies that there exist $\alpha\in (0,1)$ given in \eqref{Holder.exponent.lin}, and $\ka_\alpha' >0$\,, depending on $N,\gamma,\beta, \lambda_0, \lambda_1$ such that
\begin{equation}\label{c.alpha.regularity.before.change}
\sup_{(\hat{t},\hat{x}), (\hat{\tau},\hat{y})\in Q^*_1}\frac{|\hat{u}(\hat{t},\hat{x})-\hat{u}(\hat{\tau},\hat{y})|}{\left(|\hat{x}-\hat{y}| + |\hat{t}-\hat{\tau}|^{\frac{1}{2\vee\sigma}} \right)^\alpha}
\le\frac{\ka_\alpha'}{D^\alpha} \|\hat{u}\|_{\LL^\infty(Q_{4\hat{R}})}\le \frac{\ka_\alpha'}{D^\alpha}\,,
\end{equation}
where we have used that $\|\hat{u}\|_{\LL^\infty(Q_{4\hat{R}})}\le \|\hat{u}\|_{\LL^\infty(Q^*_{4}(\hat{t}_0,\hat{x}_0))}\le 1$, with  $D=d_{\gamma,\beta}(Q^*_2, Q^*_1)$ defined in \eqref{parabolic.distance.sets}. Due to the particular form of the cylinders $Q^*_1$ and  $Q^*_2$\,, we have that
\[
D=1\wedge \inf\limits_{\hat{y} \in B_1(\hat{x})}  \left(\rho^{\gamma,\beta}_{\hat{y}}\right)^{-1}(\hat{T} \wedge  \hat{t}_*/8)
\ge 1\wedge \ka_{19}^{-2}(\hat{T} \wedge \hat{t}_*/8)^{1/\sigma}  \wedge \ka_{19}^{-2}\left(\rho^{\gamma,\beta}_{\hat{x_0}}\right)^{-1}(\hat{T} \wedge \hat{t}_*/8):=D_0
\]
The latter inequality follows by assumptions (1), (2) or (3).
Undoing the intrinsic change of variables, \eqref{c.alpha.regularity.before.change} transforms into \eqref{Prop.HoCont.ineq} and the proof is concluded.

Finally,  note  that when $H_p $ is large enough, by Corollary \ref{weak.positivity.MDP.L1}  we know that $\kb\sim \big(R^{c_4}/\widetilde{H}_p^{c_5}\big)^{\widetilde{H}_p^{1/2}/m(1-m)} $,  hence
$ \alpha= \log_A\frac{\ka_\ell^{\,\lambda_0^{-1}+\lambda_1}}{\ka_\ell^{\,\lambda_0^{-1}+\lambda_1}-1}$ given in \eqref{Holder.exponent.lin}, with $\lambda_0,\lambda_1$ given in \eqref{lambda01bbb}, behaves like $\alpha\sim {\rm exp}\left(- \frac{c_6}{t_*}\widetilde{H}_p^{\frac{c_7}{m}\widetilde{H}_p^{1/2} }\right)$, recalling that $c_i>0$ only depend on $N,m,p,\beta,\gamma$.\qed \vspace{-4mm}

\section{Appendices}
We collect in this Appendix several technical facts and proofs, used in the rest of the paper.
\subsection{Appendix-A} \label{sec:appendix.A}
This Appendix is devoted to the proof of the upper and lower energy estimates of Lemma \ref{theorem_energy_inequality} and of the Caccioppoli estimate of Lemma \ref{Lem.Caccioppoli}.

\noindent\textbf{Approximation (via truncation) of powers of strong solutions. }The proof of the energy estimates relies on the idea of using $u^{p-1}\psi$ as a test function, where $\psi$ is a suitable smooth cutoff function and $u$ is a solution to the WFDE. As the reader may guess, this is not an admissible test function,  hence we need to proceed by a careful approximation. An additional difficulty is represented by the presence of singular/degenerate weights: under our assumptions, $u$ is merely a function in  $C_{\rm loc}((T_0, T); \LL^{p}_{\gamma, \rm loc}(\Omega))$ such that $u^m \in \LL^2_{\rm loc}((T_0, T); H^{1}_{\gamma, \beta, \rm loc}(\Omega))$; as already observed, $u$ need  not be  a function in $\LL^1_{\rm loc}\left( \Omega \right)$ and its gradient $\nabla u$ needs not to be the distributional one,  see \cite{FKS, HKM}.
The goal of the next Lemma is to show that a suitable truncation of a strong solution to \ref{WFDE} belongs to the class of admissible test functions, hence an approximate energy identity holds. Here we follow the approach used in \cite{AS-ARMA67} and in \cite{Moser}. Let $p > 1$ and $1< l< k$: we define following auxiliary functions for $u>0$ \vspace{-2mm}
\begin{equation*}
J_p(u):=\left\{
\begin{array}{lll}
((u\, \wedge \,k)^{\frac{p-1}{m}}-l^{-1})_{+}&\,\mbox{if } 1 < p \leq 1+m ,\\
(u\, \wedge \,l)^{\frac{p-1}{m}-1}(u\, \wedge \,k)&\,\mbox{if }  p > 1+m,\\
\end{array}\right.\qquad\;
G_p(u):=\int_0^u J_p(s^m)\ds\,.\vspace{-2mm}
\end{equation*}
Note that $J_p$ is a bounded Lipschitz function, for all choices of $k>l>1$ and $p>1$. Recall that $u^{p-1}$ is not an admissible test function, hence we will use a truncation of it, in the precise form of $J_p(u^m)$.

\begin{lem}\label{approximate.solutions}
Let $u$ be a non-negative strong local solution to \ref{WFDE} in $(T_0,T)\times B_{R}(x_0)$. For every $p > 1$ and for any $[t_1, t_2] \subset (T_0, T)$ the following equality holds
\begin{equation}\label{approximate.solutions.formula}
    \int_{t_1}^{t_2}\int_{B_{R}(x_0)}u_t J_p(u^m) \psi \dx \dt + \int_{t_1}^{t_2}\int_{B_{R}(x_0)} \nabla u^m \cdot \nabla (J_p(u^m)\psi) \dx \dt = 0,
\end{equation}
for any  $\psi \in C^2((T_0,T); C^2_c(B_{R}(x_0)))$.   A local strong sub (resp. super) solution satisfies  \eqref{approximate.solutions.formula} with $\le$ (resp. $\ge$) for any nonnegative test function in the same class.
\end{lem}
\noindent\textbf{Proof. }By definition $u^m \in \LL^2_{\rm loc}(T_0, T; H^{1}_{\gamma, \beta, \rm loc}(B_{R}(x_0)))$, hence there is a sequence $ \phi_n\in C^{\infty}_c((T_0, T) \times B_{R}(x_0))$ which converges strongly to $u^m$ in $\LL^2_{\rm loc}(T_0, T; H^{1}_{\gamma, \beta, \rm loc}(B_{R}(x_0)))$.  Since $J_p(\cdot)$ is a Lipschitz function, the family $\{\psi J_p(\phi)\}$ is a subset of $W^{1, 2}_{\rm loc}(T_0, T; \LL^{2}_{\gamma}(K)) \cap \LL^2_{\rm loc}(T_0, T; \D_{\gamma, \beta}(B_{R}(x_0)))$, hence an admissible test function in the sense of Definition \ref{def.local.weak.sol}, so that
\begin{equation*}\begin{split}
 \int_{\Omega}\left[u(t_2, x)\right. & \left.\psi(t_2, x)J_p(\phi_n)(t_2,x)- u(t_1, x) \psi(t_1, x)J_p(\phi_n)(t_1,x) \right] |x|^{-\gamma} \dx \\
 &=   \int_{t_1}^{t_2}\int_{\Omega}{u (\psi J_p(\phi_n))_t \ |x|^{-\gamma} \dx } \dt - \int_{t_1}^{t_2}\int_{\Omega}{ \nabla{u^m}\cdot\nabla{(\psi J_p(\phi_n))} \ |x|^{-\beta} \dx } \dt.
\end{split}\end{equation*}
Integrating by parts in time the first integral on the right-hand side, we obtain
\begin{equation}\label{approximate.solutions.formula.1}
    \int_{t_1}^{t_2}\int_{B_{R}(x_0)}u_t J_p(\phi_n) \psi \dx \dt + \int_{t_1}^{t_2}\int_{B_{R}(x_0)} \nabla u^m \cdot \nabla (J_p(\phi_n)\psi) \dx \dt = 0;
\end{equation}
the reader may observe that this integration by parts makes sense since $u$ is assumed to be a strong solution, i.e. $u_t \in \LL^1((T_0,T)\times B_{R}(x_0))$. Taking the limit as $n \rightarrow \infty$ in  \eqref{approximate.solutions.formula.1} gives \eqref{approximate.solutions.formula}. \qed

\subsubsection{Proof of the energy estimates of Lemma \ref{theorem_energy_inequality}}\label{sec.511}

We split the proof of Lemma \ref{theorem_energy_inequality} in several parts: first we prove the upper estimate, then the lower.

\medskip

\noindent\textbf{Proof of the upper energy inequality \eqref{sup.energy.inequality.upper}}
Let us fix $x_0\in \RR^N$\,, and simply denote $B_R=B_R(x_0)$ when there is no ambiguity.

\noindent$\bullet~$\textsc{Step 1. }\textit{Reduction. }The upper energy inequality \eqref{sup.energy.inequality.upper} follows by a slightly different inequality:
\begin{align}\label{upper_energy_inequality}
   \int_{B_{R_1}} & u(T, x)^p \, |x|^{-\gamma} \dx  +
   \int_{T_1}^T\int_{B_{R_1}}{\left|\nabla u^{\frac{p+m-1}{2}} \right|^2 \ |x|^{-\beta} \dx  \dt} \nonumber \\
  & \leq C\left(m, p\right) \left[\frac{ h_{\sigma}\left( R, R_1, x_0 \right)}{\left(R-R_1 \right)^{\sigma}} + \frac{1}{T_1-T_0} \right] \int_{T_0}^T\int_{B_R}{\left( u^{p+m-1} + u^{p} \right) \ |x|^{-\gamma} \dx \dt}\,.
\end{align}
Indeed, inequality \eqref{sup.energy.inequality.upper} follows from \eqref{upper_energy_inequality} by letting $ T=\tau $ and taking the supremum in $ \tau \in \left[T_1, T \right]$.

\noindent$\bullet~$\textsc{Step 2. }\textit{First energy inequality. }In this step we want to prove the following inequality:
\begin{equation}\label{energy.inequality.1}\begin{split}
\frac{p-1}{p}\int_{B_R}\left[u(T, x)^p \right. &\,\psi^2\left(T, x\right) - \left.u(T_0, x)^p \psi^2\left(T_0, x\right) \right] \ |x|^{-\gamma} \dx \\
 &\,+ \frac{2 m \left(p-1 \right)^2}{\left(p+m-1 \right)^2} \int_{T_0}^{T}\int_{B_R} \psi^2  \left|\nabla u^{\frac{p+m-1}{2}} \right|^2 \ |x|^{-\beta} \dx  \dt  \\
  &\,\leq 2 \left[\int_{T_0}^{T}\int_{B_R}{u^p \ \psi \ |\psi_t| \ |x|^{-\gamma} \dx  \dt} +  m \int_{T_0}^{T}\int_{B_R}{u^{p+m-1} \ |\nabla \psi|^2\ |x|^{-\beta} \dx  \dt} \right].
\end{split}\end{equation}
Following Moser's approach \cite{Moser}, we would like to test the equation with $u^{p-1}\psi^2$, but unfortunately this is not an admissible test function: we shall proceed by approximation, using $J_p(u^m)\psi^2$ as in Lemma \ref{approximate.solutions}; this approximation  extends  to our weighted setting some ideas from Aronson and Serrin \cite{AS-ARMA67}:
\begin{equation}\label{approx.inequality.0}
  \int_{T_0}^{T}\int_{B_{R}}u_t J_p(u^m) \psi^2 |x|^{-\gamma}\dx \dt + \int_{T_0}^{T}\int_{B_{R}} \nabla u^m \cdot \nabla (J_p(u^m)\psi^2) |x|^{-\beta}\dx \dt = 0.
\end{equation}
 Recalling  that $\partial_t G_p(u)= J_p(u^m) u_t$, an integration by parts (in time) in the left-hand side of \eqref{approx.inequality.0} gives
\begin{equation}\label{approx.inequality.1}\begin{split}
  \int_{B_R}&\left[G_p(u)(T,x)\,\psi^2\left(T, x\right) - G_p(u)(T_0, x)^p \psi^2\left(T_0, x\right) \right] \ |x|^{-\gamma} \dx \\
   & + \int_{T_0}^{T}\int_{B_{R}} \nabla u^m \cdot \nabla (J_p(u^m)\psi^2) |x|^{-\beta}\dx \dt
   \le 2\int_{T_0}^{T}\int_{B_{R}} |\psi| |\psi_t| G_p(u) |x|^{-\gamma} \dx\dt.
\end{split}\end{equation}
 Note that $J_p(u^m)\to \tilde{J_p}(u^m) $ and $G_p(u) \to \tilde{G}_p(u)$ as $k \rightarrow \infty$ where
\begin{equation*}
\tilde{J_p}(u):=\left\{
\begin{array}{lll}
(u^{\frac{p-1}{m}}-l^{-1})_{+}&\,\mbox{if } 1 < p \leq 1+m ,\\
(u \wedge l)^{\frac{p-1}{m}-1}\,u&\,\mbox{if }  p > 1+m,\\
\end{array}\right.
\tilde{G_p}(u):=\int_0^u \tilde{J_p}(s^m)\ds\,.
\end{equation*}

Since $J_p$ is Lipschitz, taking limits as $k\to\infty$ in inequality \eqref{approx.inequality.1} we get  (by dominated convergence)  :
\begin{equation}\label{approx.inequality.2}\begin{split}
  \int\limits_{B_R}\left[\tilde{G_p}(u)(T,x) \right. &\,\psi^2\left(T, x\right) - \left.\tilde{G_p}(u)(T_0, x)^p \psi^2\left(T_0, x\right) \right] \frac{\dx}{|x|^{\gamma}} + \int\limits_{T_0}^{T}\int\limits_{B_{R}} \psi^2 \tilde{J_p}'(u^m) \left| \nabla u^{m} \right|^2 \frac{\dx \dt}{|x|^{\beta}}\\
   & \le 2\int_{T_0}^{T}\int_{B_{R}} |\psi| |\psi_t| \tilde{G_p}(u) |x|^{-\gamma} \dx\dt - 2 \int_{T_0}^{T}\int_{B_{R}} \psi \tilde{J_p}(u^m) \nabla u^m \cdot \nabla \psi\, |x|^{-\beta}\dx\dt.
\end{split}\end{equation}
We combine now the following numerical inequality $ \tilde{J}_p^2(u^m) \leq \left(\frac{m}{p-1}\right) \tilde{J}_p'(u^m) u^{p+m-1}\left[1+ \left(\frac{p-1}{m}-1\right)\one\{u^m>l\}\right]$ with Young's inequality $|v \cdot w| \leq |v|^2 / 4 + |w|^2$  to obtain
\begin{equation}\label{approx.inequality.2bb}
2 |\psi \tilde{J_p}(u^m) \nabla u^m \cdot \nabla \psi |
\leq \frac{1}{2} \psi^2 \tilde{J_p}'(u^m) |\nabla u^m|^2 + \frac{2m}{p-1}u^{p+m-1}|\nabla \psi|^2 f(l,u),
\end{equation}
where $f(l,u)=\left[1+ \left(\frac{p-1}{m}-1\right)\one\{u^m>l\}\right]$. Combining \eqref{approx.inequality.2} and \eqref{approx.inequality.2bb} we get
\begin{equation*}
\begin{split}
\int_{B_R}\left[\tilde{G_p}(u)(T,x) \right. &\,\psi^2\left(T, x\right) - \left.\tilde{G_p}(u)(T_0, x)^p \psi^2\left(T_0, x\right) \right] |x|^{-\gamma} \dx + \frac{1}{2}\int_{T_0}^{T}\int_{B_{R}} \psi^2 \tilde{J_p}'(u^m) \left| \nabla u^{m} \right|^2 |x|^{-\beta}  \dx \dt\\
   & \le 2\int_{T_0}^{T}\int_{B_{R}} |\psi| |\psi_t| \tilde{G_p}(u) |x|^{-\gamma} \dx\dt + \frac{2m}{p-1}  \int_{T_0}^{T}\int_{B_{R}} u^{p+m-1} |\nabla \psi|^2 f(l,u)|x|^{-\beta} \dx \dt.
\end{split}
\end{equation*}
Finally, we obtain \eqref{energy.inequality.1} by taking the limit as $l \rightarrow \infty$ in the above inequality:  we notice that in such limit $\tilde{J_p}'(u^m) \left| \nabla u^{m} \right|^2\to \tilde{c}_{p,m}\left|\nabla u^{\frac{p+m-1}{2}} \right|^2$ in the appropriate integral sense, where $\tilde{c}_{p,m}$ is a suitable multiplicative constant, as well $\tilde{G}_p \to  u^p/p$ and $f(l,u) \to 1$ by dominated convergence.

\medskip

\noindent$\bullet~$\textsc{Step 3. }\textit{Choice of the test function $\psi$. }By a suitable choice of test function, we can show that inequality \eqref{energy.inequality.1} implies \eqref{upper_energy_inequality}.  It is always possible to choose a smooth $0\le \psi\le 1$ supported in $[T_0,T]\times B_R$, such that $\psi\equiv 1$ on $[T_1,T]\times B_{R_1}$, $\psi(T_0,x)=0$ for all $x\in B_R$ and $\psi(t,x)=0$ for all $(t, x) \in [T_0,T]\times \partial B_R$\,, so that there exists $K_\psi>0$ (depending only on $N$) such that $| \nabla \psi(t,x)|^2 \leq K_{\psi} \left(R -R_1 \right)^{-2}$ and $|\psi_t(t,x) | \leq K_{\psi}(T_1-T_0)^{-1}$ for all $(t,x)\in  \left(T_0, T\right]\times B_R \setminus B_{R_1}$. With this test function, we estimate the two sides of \eqref{energy.inequality.1} separately.

\noindent\textit{Estimating the right-hand side of \eqref{energy.inequality.1}. }We will show that

\begin{align}\label{energy.inequality.2b}
 &  \left[\int_{T_0}^{T}\int_{B_R}{u^p \ \psi \ |\psi_t| \ |x|^{-\gamma} \dx  \dt} + \int_{T_0}^{T}\int_{B_R}{u^{p+m-1} \ |\nabla \psi|^2\ |x|^{-\beta} \dx  \dt} \right]  \nonumber \\
 &  \leq 2 K_{\psi}  \left[\frac{1}{T_1-T_0} +  \frac{h_{\sigma}\left( R, R_1, x_0 \right)}{\left(R - R_1 \right)^{\sigma}} \right] \left[\int\limits_{T_0}^{T}\int\limits_{B_R}{u^p \frac{\dx \dt}{|x|^{\gamma} }} + \int_{T_0}^{T}\int\limits_{B_R}{u^{p+m-1} \frac{\dx \dt}{|x|^{\gamma} }} \right],
\end{align}
where the function $ h_{\sigma}\left( R, R_1, x_0 \right) $ is defined in \eqref{function_h}. Indeed,
plugging the above chosen $\psi$ in the right-hand side of inequality \eqref{energy.inequality.1} we get
\begin{align}\label{energy.inequality.2}
  c_{m, p} & \left[\int_{B_R}{u(T, \cdot)^p \ \psi^2\left(T, \cdot \right) \ |x|^{-\gamma} \dx} + \int_{T_0}^{T}\int_{B_R}{  \psi^2 \ |\nabla u^{\frac{p+m-1}{2}} |^2 \ |x|^{-\beta} \dx  \dt}  \right] \nonumber \\
 &  \leq 2 \left[\int_{T_0}^{T}\int_{B_R}{u^p \ \psi \ |\psi_t| \ |x|^{-\gamma} \dx  \dt} + \int_{T_0}^{T}\int_{B_R}{u^{p+m-1} \ |\nabla \psi|^2\ |x|^{-\beta} \dx  \dt} \right]
 \end{align}
where the constant $ c_{m, p}  $ is given by
\begin{equation*}\label{constant.sup.energy.inequality.upper}
 c_{m, p} =  \frac{p-1}{p} \wedge \frac{2m \left(p-1 \right)^2}{\left(p+m-1 \right)^2}.
\end{equation*}
 We just have to estimate the quotient $|x|^{-\beta}/|x|^{-\gamma}$ in the right-hand side of \eqref{energy.inequality.2} in terms of $h_{\sigma}\left( R, R_1, x_0 \right) $ and $R-R_1$ to get \eqref{energy.inequality.2b}.  First, recall that  $ \nabla \psi\left(t, \cdot \right) $ is supported in $ B_R\left( x_0 \right) \setminus B_{R_1}\left( x_0 \right) $ for all $ t \in \left(T_0, T_1\right] $. Next we split two cases, namely $\sigma<2$ and $\sigma\ge 2$.\\
\noindent\textit{-- Case $0 < \sigma < 2 $.} In $\overline{B_{R}(x_0)\setminus B_{R_1}(x_0)}$ we have $ |x|^{-\beta} = |x|^{\gamma-\beta} |x|^{-\gamma} \leq \left( |x_0| + R \right)^{\gamma-\beta} |x|^{-\gamma}$, hence
\begin{equation}\label{1111}
 \frac{\left( |x_0|+R \right)^{\gamma-\beta}}{\left(R - R_1 \right)^2} = \left(\frac{ |x_0|+R }{R-R_1}\right)^{\gamma-\beta} \frac{1}{\left(R-R_1 \right)^{\sigma}}=\dfrac{h_\sigma\left( R, R_1, x_0 \right)}{\left( R-R_1 \right)^\sigma},
\end{equation}
since we recall that $ 0 <\sigma =2+\beta-\gamma< 2 $ means $ \gamma > \beta $ and $ \frac{R+|x_0|}{R-R_1}\geq\frac{R}{R-R_1} > 1$.

\noindent\textit{-- Case $ \sigma \geq 2 $. }Recall that $ \sigma \geq 2 $ means that $ \gamma \leq \beta $. We now consider two sub-cases.
Recall that we consider balls $B_{R_1}(x_0)\subset B_R(x_0)$ such that $0 \not \in \overline{B_{R}(x_0)\setminus B_{R_1}(x_0)}$\,.\\
If $ 0 \not \in B_R(x_0)$, i.e. $|x_0|> R$, then in $\overline{B_{R}(x_0)\setminus B_{R_1}(x_0)}$ we have   $|x|^{-\beta} = |x|^{-\gamma}/ |x|^{\beta - \gamma}  \leq  |x|^{-\gamma}/\left( |x_0| - R \right)^{\beta-\gamma}$, hence
\begin{equation}\label{2222}
 \frac{\left( |x_0| - R \right)^{\gamma-\beta}}{\left(R - R_1 \right)^2} = \left(\frac{R - R_1}{ |x_0| - R }\right)^{\beta-\gamma} \frac{1}{\left(R-R_1 \right)^{\sigma}}\leq \frac{h_{\sigma}\left( R, R_1, x_0 \right)}{\left( R-R_1 \right)^{\sigma}}.
\end{equation}
 If $0 \in B_{R_1}(x_0)$, i.e. $R_1 > |x_0|$ then in $\overline{B_{R}(x_0)\setminus B_{R_1}(x_0)}$ we have $|x|^{-\beta} = |x|^{-\gamma}/|x|^{\beta - \gamma} \leq  |x|^{-\gamma}/(R_1-|x_0|)^{\beta-\gamma}$, hence
\begin{equation}\label{3333}
 \frac{(R_1-|x_0|)^{\gamma-\beta}}{\left(R - R_1 \right)^2} = \left(\frac{R - R_1}{R_1-|x_0|}\right)^{\beta-\gamma} \frac{1}{\left(R-R_1 \right)^{\sigma}} \leq \frac{h_{\sigma}\left( R, R_1, x_0 \right)}{\left( R-R_1 \right)^\sigma}.
\end{equation}

\noindent\textit{Estimating the r left-hand side and the right-hand side of \eqref{energy.inequality.2}. }
We just notice that since $ \psi = 1 $ in $ [T_1, T] \times B_{R_1} $ we get
\begin{equation*}
 \int_{B_{R_1}}{u(T,x)^p \ |x|^{-\gamma} \dx} \leq  \int_{B_R}{u(T,x)^p \ \psi^2\left(T, x\right) \ |x|^{-\gamma} \dx},
\end{equation*}
and
\begin{equation*}
 \int_{T_1}^{T}\int_{B_{R_1}}{  |\nabla u^{\frac{p+m-1}{2}} |^2 \ |x|^{-\beta} \dx  \dt} \leq \int_{T_0}^{T}\int_{B_R}{  \psi^2 \ |\nabla u^{\frac{p+m-1}{2}} |^2 \ |x|^{-\beta} \dx  \dt}.
\end{equation*}
Summing up, inequality \eqref{energy.inequality.2} becomes
\begin{align*}
   c_{m,p}\left[\int_{B_{R_1}}\right. &\,u(T, x)^p |x|^{-\gamma} \dx  + \left.\int_{T_1}^{T}\int_{B_{R_1}}{|\nabla u^{\frac{p+m-1}{2}} |^2 \ |x|^{-\beta} \dx  \dt}\right] \nonumber \\
  & \leq 2K_{\psi} \left[\frac{h_{\sigma}\left( R, R_1, x_0 \right)}{\left(R-R_1 \right)^{\sigma}} + \frac{1}{T_1-T_0} \right] \int_{T_0}^{T}\int_{B_R}{\left( u^{p+m-1} + u^{p} \right) \ |x|^{-\gamma} \dx \dt}\,,
\end{align*}
where $c_1\equiv C\left(m, p\right):= 2K_{\psi} \ c_{m,p}^{-1}$ and $K_\psi>0$ depends on $N$. The proof of the inequality \eqref{sup.energy.inequality.upper} is concluded.\qed

\subsubsection{Proof of the lower energy inequalities   \eqref{energy.inequality.last} and \eqref{sup.energy.inequality.lower}}
We will perform before a common step, used in the proof of both inequalities. Let us fix $x_0\in \RR^N$\,, and simply denote $B_R=B_R(x_0)$ when no confusion arises. We always consider $p\in \RR\setminus\{0\}$.

\noindent$\bullet~$\textsc{Step 1. }\textit{First energy inequality. }In this step prove the following inequality for $-p<1-m$\,, $p\ne 0$:
\begin{equation}\label{sup.energy.inequality.lower.last.1}\begin{split}
 \frac{p+1 }{p} & \left[ \int_{B_R}{\left(u(T, x)^{-p} \ \psi^2\left(T, x \right) -  u(T_0, x)^{-p} \ \psi^2\left(T_0,x\right) \right)  |x|^{-\gamma} \dx}  \right]  \\
 &+ \frac{ 2  m \left(p+1 \right)^2}{\left(m-p-1 \right)^2} \int_{T_0}^{T}\int_{B_R}{    \left|\nabla u^{\frac{m-p-1}{2}} \right|^2  \psi^2 |x|^{-\beta} \dx  \dt}  \\
 & \leq 2m\int_{T_0}^{T}\int_{B_R}{u^{m-p-1}\ |\nabla \psi|^2 \ |x|^{-\beta} \dx  \dt } + 2\frac{p+1 }{p} \int_{T_0}^{T}\int_{B_R}{u^{-p} \ \psi \ |\psi_t|  \ |x|^{-\gamma} \dx \dt }.
\end{split}\end{equation}
We just sketch the proof, since it is similar - but simpler - to Step 2 of the proof of the upper energy inequality \eqref{sup.energy.inequality.upper}: we approximate $u^{-p-1}\psi^2$ with admissible test functions, we use the weak formulation of the equation \eqref{weak.formulation} and after a double limiting process, we obtain:
\begin{equation*}\label{sup.energy.inequality.lower.last.2}\begin{split}
   -\left(p+1 \right)\int_{T_0}^{T}\int_{B_R}u_t \, u^{-\left(p+1 \right)}  \psi^2  |x|^{-\gamma} \dx \dt + \frac{4 m \left(p+1 \right)^2}{\left(m-p-1 \right)^2} \int_{T_0}^{T}\int_{B_R}{    \left|\nabla u^{\frac{m-p-1}{2}} \right|^2  \psi^2 |x|^{-\beta} \dx  \dt} &  \\
  \leq  \frac{4m\left(p+1 \right)}{\left(m-p-1\right)} \int_{T_0}^{T}\int_{B_R}{\psi  \nabla u^{\frac{m-p-1}{2}} \cdot u^{\frac{m-p-1}{2}} \nabla \psi \ |x|^{-\beta} \dx \dt}&.
\end{split}\end{equation*}
Using Young's inequality, i.e., $|v \cdot w| \leq \frac{|v|^2}{2\varepsilon} + \frac{\varepsilon |w|^2}{2}$, with $\varepsilon= -\frac{p+1 }{m-p-1 } > 0 $ we get the following inequality
\begin{equation*}\label{sup.energy.inequality.lower.last.3}\begin{split}
  \frac{p+1}{p}  \int_{T_0}^{T}\int_{B_R}{\partial_t(u^{-p}) \psi^2   \ |x|^{-\gamma} \dx \dt} + \frac{2 m \left(p+1 \right)^2}{\left(m-p-1 \right)^2} \int_{T_0}^{T}\int_{B_R}{   \left|\nabla u^{\frac{m-p-1}{2}} \right|^2  \psi^2 |x|^{-\beta} \dx  \dt}&  \\
     \leq  2m\int_{T_0}^{T}\int_{B_R}{u^{m-p-1}\ |\nabla \psi|^2 \ |x|^{-\beta} \dx  \dt }&.
\end{split}\end{equation*}
Integrating by parts in time the first term of the above inequality, we obtain \eqref{sup.energy.inequality.lower.last.1}.

\medskip

In the following steps we show that \eqref{sup.energy.inequality.lower.last.1} implies both  \eqref{energy.inequality.last} and \eqref{sup.energy.inequality.lower}. We will just sketch the proof, since it is very similar to the proof of inequality \eqref{sup.energy.inequality.upper}.

\noindent$\bullet~$\textsc{Step 2. }\textit{Proof of inequality \eqref{energy.inequality.last}}.   In order to keep the same notation in the proofs, we will change the sign of the exponent $p$ with respect to the statement of inequality \eqref{energy.inequality.last}, namely we will switch $p$ to  $-p$.   So $m-1 < p < 0$, hence we have $ \frac{p+1}{p} <  0 $.  We follow the Steps 3 of the proof of inequality \eqref{sup.energy.inequality.upper}: we choose a smooth $0\le \psi\le 1$ supported in $[T_0,T]\times B_R$, such that $\psi\equiv 1$ on $[T_0,T_1]\times B_{R_1}$, $\psi(T,x)=0$ for all $x\in B_R$ and $\psi(t,x)=0$ for all $(t, x) \in [T_0,T]\times \partial B_R$\,, so that there exists $K_\psi>0$ (depending only on $N$) such that $| \nabla \psi(t,x)|^2 \leq K_{\psi} \left(R -R_1 \right)^{-2}$ and $|\psi_t(t,x) | \leq K_{\psi}(T-T_1)^{-1}$ for all $(t,x)\in  \left(T_0, T\right]\times B_R \setminus B_{R_1}$.
With this choice of $\psi$ we obtain the following inequality:
\begin{equation*}\label{sup.energy.inequality.last.5}\begin{split}
  \int_{B_{R_1}(x_0)}& u(T_0,x)^{-p}  |x|^{-\gamma} \dx+ \int_{T_0}^{T_1}\int_{B_{R_1}} |\nabla u^{\frac{m-p-1}{2}}|^2  |x|^{-\beta} \dx  \dt  \\
& \leq c_4\left[ \int_{T_0}^{T}\int_{B_R}u^{m-p-1}|\nabla \psi|^2  |x|^{-\beta} \dx  \dt + \int_{T_0}^{T}\int_{B_R}u^{-p} |\psi_t| |x|^{-\gamma} \dx \dt \right].
\end{split}\end{equation*}
Proceeding as in Step 3 of the proof of inequality \eqref{sup.energy.inequality.upper} we obtain inequality \eqref{energy.inequality.last}, with
\begin{equation*}
c_2 =  K_{\psi} \frac{2\left(m \vee \frac{|p+1|}{|p|}\right)}{|p+1|\left(\frac{1}{|p|} \wedge \frac{2m|p+1|}{(m-p-1)^2}\right)}.
\end{equation*}

\noindent$\bullet~$\textsc{Step 3. }\textit{Proof of inequality \eqref{sup.energy.inequality.lower}. }We choose $\psi$ as in Step 3 of the proof of inequality \eqref{sup.energy.inequality.upper} and repeating the same estimates used there, we can estimate \eqref{sup.energy.inequality.lower.last.1} to get
\begin{equation*}\label{sup.energy.inequality.lower.last.4}\begin{split}
  & \int_{B_{R_1}(x_0)}{u(T, x)^{-p} \ |x|^{-\gamma} \dx} + \int_{T_1}^{T}\int_{B_R}{\left|\nabla u^{\frac{-p+m-1}{2}} \right|^2 \ |x|^{-\beta} \dx  \dt}  \\
  & \leq c_3 \left[\frac{h_\sigma(R, R_1, x_0)  }{\left(R-R_1 \right)^\sigma} + \frac{1}{T_1-T_0} \right] \int_{T_0}^{T}\int_{B_R}{\left( u^{-p+m-1} + u^{-p} \right) \ |x|^{-\gamma} \dx \dt},
\end{split}\end{equation*}
Finally, inequality \eqref{sup.energy.inequality.lower} follows by letting $ T=\tau $ and taking the supremum in $ \tau \in \left[T_1, T \right]$ in the above inequality. The constant $c_3>0$ becomes
\begin{equation*}\label{formula.c2}
c_3=\frac{4 K_{\psi} }{p+1}\,\frac{  m \vee \frac{p+1}{p}   }{ \frac{1}{p}\wedge \frac{4 m  (p+1)}{ (m-p-1)^2} }>0,\qquad\mbox{since }  p>0.
\end{equation*}
The proof of Lemma \ref{theorem_energy_inequality}  is  now concluded.\qed

\subsubsection{Proof of the Caccioppoli estimates of Lemma \ref{Lem.Caccioppoli}. }
We just sketch the proof. We use the test function $\psi^2 u^{-m}$, assuming first $0<\delta\le u\le M$, and we approximate it as in Step 2 of the proof of the upper energy inequality \eqref{sup.energy.inequality.upper} so that we obtain
\begin{equation}\label{caccioppoli.inequality.1}\begin{split}
 - \iint_{Q}{\psi^2 \partial_t(u^{1-m})  |x|^{-\gamma} \dx \dt} &+m^2 \left(1-m \right)\iint_{Q}{ \psi^2  |\nabla \log u|^2 |x|^{-\beta} \dx \dt} \\
 &\leq 2m\left(1-m \right) \iint_{Q}{ \psi \nabla \log u \cdot \nabla \psi \ |x|^{-\beta} \dx \dt},
\end{split}\end{equation}
where $Q=\left(\tau, t\right) \times B_R(x_0)$. Inequality \eqref{caccioppoli.inequality} follows by using Young's inequality $ab\le \varepsilon a^2+b^2/4\varepsilon$, with  $\varepsilon= m/4$ , on the right-hand side of \eqref{caccioppoli.inequality.1} and integrating by parts in time the first term of inequality \eqref{caccioppoli.inequality.1}. Note that the assumption $u\in [\delta,M]$ can be removed by a lengthy but straightforward approximation, but we refrain from doing this here, since we apply \eqref{caccioppoli.inequality} only to solutions to a ``lifted'' Dirichlet problem \eqref{LARGE.DIRICHLET.PROBLEM}, which we already know to be positive and bounded.\qed

\subsection{Appendix-B} \label{sec:appendix.C}
The goal of this Appendix is to prove the weighted Caffarelli-Kohn-Nirenberg and Poincar\'{e} Inequalities of Propositions \ref{prop.Sobolev.Balls} and \ref{zero.mean.Poincare}  and to provide some useful quantitative information about the auxiliary function $ \rho^{\gamma, \beta}_{x_0}$ and its inverse.  We recall here the expression of $ \rho^{\gamma, \beta}_{x_0}$ defined in \eqref{rho.pseudo}:
\begin{equation*}
\rho_{x_0}^{\gamma,\beta}(R):=\left( \int_{B_R(x_0)} |x|^{(\beta - \gamma) \frac{N}{2}} \dx \right)^{\frac{2}{N}}.
\end{equation*}
We begin with a technical lemma on the behaviour of the function  $\rho^{\gamma, \beta}_{x_0}$.
\begin{lem}\label{technical.lemma.measures}
    Let $N \ge 3$, assume that $ \gamma, \beta \in \RR$ satisfy \eqref{paramt.range}. Then there exists $\ka_{16}, \ka_{18} >0 $ such that for any $y \in \RR^N$ and for any $R>0$ the following inequalities hold
    \begin{align}\label{technical.inequality.measures}
      \ka_{16}^{-1}\, \rho_{y}^{\gamma,\beta}(R) & \leq   R^2 \frac{\mu_\gamma(B_R(y))}{\mu_\beta(B_R(y))} \leq \ka_{16}  \rho_{y}^{\gamma,\beta}(R)\,,\\ \nonumber
      \ka_{18}^{-1} R^2 \left[R \vee |y|\right]^{\beta-\gamma} &\leq  \rho_{y}^{\gamma,\beta}(R)  \leq \ka_{18} R^2 \left[R \vee |y|\right]^{\beta-\gamma}
     \end{align}
The constants $\ka_{16},\ka_{18}>0$ depend only on $N, \gamma, \beta$.
\end{lem}
\noindent\textbf{Proof. }  We will only prove  the first inequality appearing in \eqref{technical.inequality.measures}, the second one will follow by the same estimates, noticing that Step 2 and 3 correspond to the cases $|y|\le 2R$ and $|y|>2R$.  The proof will be divided in different cases.\\
\noindent$\bullet~$\textsc{Case 1. } \textit{Assume $y=0$ and $R > 0$.} This case is done  by a direct calculation.\\
\noindent$\bullet~$\textsc{Case 2. } \textit{Assume that $0 < |y| \le 2R $}. The reader may observe that in this case the following inclusions  holds
$ B_r(y) \subset B_{4r}(0) \subset B_{8r}(y)$ .
Then, by  the  doubling property, we obtain the following inequalities (recall that $\sigma=2+\beta - \gamma$)
\begin{equation*}\begin{split}
\mu_{-(\sigma-2)\frac{N}{2}}(B_R(y))^{\frac{2}{N}} \mu_\beta(B_R(y))
& \le \mu_{-(\sigma-2)\frac{N}{2}}(B_{4R}(0))^{\frac{2}{N}} \mu_\beta(B_{4R}(0))
\le C_1 R^{2+\beta - \gamma} R^{N-\beta}  \\ &\le C_2 R^2 \mu_\gamma(B_{4R}(0))
 \le C_3 R^2 \mu_\gamma(B_{8R}(y))   \le C_4 R^2 \mu_\gamma(B_{R}(y)).
\end{split}\end{equation*}
The other inequality is obtained by similar techniques. \\
\noindent$\bullet~$\textsc{Case 3. } \textit{Assume that $0 < 2R < |y| $}. Assume that $z \in B_R(y)$, therefore $\frac{|y|}{2}\leq |z| \le \frac{3|y|}{2}$. In order to prove inequality \eqref{technical.inequality.measures} we will show that the quantity $I$ defined by
\begin{equation*}
I = \left( \frac{1}{R^N}\int_{B_R(y)} |z|^{(\sigma-2)\frac{N}{2}} \dz \right)^\frac{2}{N} \left( \frac{1}{R^N} \int_{B_R(y)} |z|^{-\beta} \dz  \right) \left( \frac{1}{R^N}\int_{B_R(y)} |z|^{-\gamma} \dz  \right)^{-1},
\end{equation*}
is bounded (above and below) by a constant independent of $y$ and $R$. For any $\alpha > -N $ we can estimate $\int_{B_R(y)} |z|^{\alpha} \dz$ as
\begin{equation*}
C_5 R^N |y|^\alpha  \le \int_{B_R(y)} |z|^{\alpha} \dz \le C_6 R^N |y|^{\alpha},
\end{equation*}
where the  constants $C_5$ and $C_6$ depend  only on the dimension $N$. Therefore the quantity $I$ is bounded (above) by
\begin{equation*}
I \le C_7 |y|^{\beta -\gamma} |y|^{-\beta} |y|^{\gamma} \leq C_7,
\end{equation*}
recall that $\sigma=2+\beta-\gamma$. The very same technique works also for the other bound. \qed

\medskip
The function $\rho^{\gamma, \beta}_{x_0}(r)$, is increasing in $r$ therefore it has an inverse which we denote by $\left(\rho^{ \gamma, \beta}_{x_0}\right)^{-1}$, whose behaviour we show in the next lemma.
\begin{lem}\label{inverse.technical.inequality.measures.lemma}
Let $ \gamma, \beta \in \RR$ satisfy \eqref{paramt.range} and $N\ge 3$. Then there exists $\ka_{19} >0 $ such that for any $x_0 \in \RR^N$   and for any $s>0$ the following inequalities hold
    \begin{equation}\label{inverse.technical.inequality.measures}\begin{split}
      \ka_{19}^{-1}\, s^{\frac{1}{2}} \left[s^{\frac{1}{\sigma}} \vee |x_0|\right]^{\frac{\gamma-\beta}{2}} \leq   \left(\rho^{ \gamma, \beta}_{x_0}\right)^{-1}(s) \leq  \ka_{19}\, s^{\frac{1}{2}} \left[s^{\frac{1}{\sigma}} \vee |x_0|\right]^{\frac{\gamma-\beta}{2}}\,,
    \end{split}\end{equation}
where constant $\ka_{19}>0$ depends only on $N, \gamma, \beta$. As a consequence, for any $x_0 \in \Omega\subset\RR^N$   and for any $s\in (0,T]$, we have:
\begin{equation}\label{bound.inverse.rho}
\left(\rho^{ \gamma, \beta}_{x_0}\right)^{-1}(s)\le \left\{
\begin{array}{lll}
\ka_{19}\,s^{\frac{1}{\sigma}}\,,&\,\mbox{if }  \sigma \ge 2 ,\\
\ka_{19}\,s^{\frac{1}{2}}\,\left(T^{\frac{1}{\sigma}}\vee \sup\limits_{x_0\in\Omega}|x_0|\right)^{\frac{\gamma-\beta}{2}}\,,&\,\mbox{if }  0<\sigma < 2\,.\\
\end{array}\right.
\end{equation}
\end{lem}
\noindent\textbf{Proof. } We first observe that inequality \eqref{bound.inverse.rho} easily follows by \eqref{inverse.technical.inequality.measures}, hence we only have to prove the latter.  \\
\noindent$\bullet~$\textsc{Case 1. } \textit{Assume $x_0=0$.}  Under this assumption we know that $ \rho^{\gamma, \beta}_0(r) \asymp r^\sigma $. Therefore $\left(\rho^{ \gamma, \beta}_{x_0}\right)^{-1}(s) \asymp s^\frac{1}{\sigma}$.
\noindent$\bullet~$\textsc{Case 2. } \textit{Assume $x_0\neq0$.} Here we deal with two different cases. First, we observe that if $0\leq r \leq |x_0|$ we have  $\rho^{\gamma, \beta}_{x_0}(r) \asymp r^2 |x_0|^{\beta-\gamma}$, therefore $\left(\rho^{ \gamma, \beta}_{x_0}\right)^{-1}(s) \asymp s^\frac{1}{2}|x_0|^{\frac{\gamma-\beta}{2}}$ and the estimate holds when $r\asymp s^\frac{1}{2} |x_0|^\frac{\gamma-\beta}{2} \leq |x_0|$, i.e. when $s^\frac{1}{\sigma}\leq |x_0|$. Next, when $0 \leq |x_0| \leq r$ we have $\rho^{\gamma, \beta}_{x_0}(r) \asymp r^\sigma$ and therefore $\left(\rho^{ \gamma, \beta}_{x_0}\right)^{-1}(s) \asymp s^\frac{1}{\sigma}$, the estimate holds when $s^\frac{1}{\sigma}\geq |x_0|$. The two estimates give \eqref{inverse.technical.inequality.measures}, and this concludes the proof.\qed
\noindent\textbf{Proof of the weighted Poincar\'e inequality of Proposition \ref{zero.mean.Poincare}. }The Poincar\'e inequality \eqref{Weighted.Poincare.ineq} will easily follow from H\"older's inequality and from
the following weighted Sobolev-Poincar\'e inequality proven in \cite[Theorem I]{FGW}
\begin{equation}\label{Franchi.sobolev.poincare}
\left(\int_{B_R(y)} |\phi - \overline{\phi}|^{\sr} |x|^{-\gamma} \dx \right)^{\frac{1}{\sr}} \leq C_1 R \frac{ \mu_\gamma(B_R(y))^{\frac{1}{\sr}}}{\mu_\beta(B_R(y))^{\frac{1}{2}}}\left(\int_{B_R(y)} |\nabla \phi|^2 |x|^{-\beta} \dx \right)^{\frac{1}{2}}.
\end{equation}
where $\overline{\phi}=\mu_{\gamma}\left( B_R(y)\right)^{-1} \int_{B_R(y)}{\phi |x|^{-\gamma}  \dx }$,  $B_R(y)$ is any ball  and $C_1>0$ depends only on $N, \gamma$ and $\beta$.\qed
\medskip

\noindent\textbf{Proof of Proposition \ref{prop.Sobolev.Balls}. }
 Inequality \eqref{WSIcomplete} follows from \eqref{Franchi.sobolev.poincare}, estimating the constant as in the above proof, then using $\|f-\overline{f}\|_{\LL^p_\gamma(B_R(x_0))} \ge \|f\|_{\LL^p_\gamma(B_R(x_0))} - \overline{f} \mu_\gamma(B_R(x_0))^{\frac{1}{p}}$ and  H\"{o}lder's inequality.\qed

The following technical lemma is needed in the proof of Proposition \ref{Prop.Lin.Harn.HoCont}.
\begin{lem}\label{multiplicative.factor}
    Let $N \ge 3$, assume that $ \gamma, \beta \in \RR$ satisfy \eqref{paramt.range}. For any positive real number $A$,
    \begin{equation}\label{expression.factor}
    A \geq 4 \vee 2\ka_{18} \vee \left(4\ka_{18}^2\right)^\frac{1}{\sigma},
    \end{equation}
and for any $r>0$, for any $x_0 \in \RR^N$ the following inclusion holds
\begin{equation}\label{inclusion.parabolic.cylinders}
Q_{R/A}(t_0,x_0) \subset Q_R^+(t_0,x_0),
\end{equation}
where $Q_R^+$ and $Q_R$  are defined in \eqref{Parabolic.cylinders.lin}, and $\ka_{18}>0$ is as in \eqref{technical.inequality.measures}.
\end{lem}
\noindent\textbf{Proof. }  We prove only in the case $0 < \sigma < 2$, namely $\gamma>\beta$, since the case $\sigma\ge 2$ is actually simpler and follows by the very same steps. In order to prove the inclusion \eqref{inclusion.parabolic.cylinders} we need to verify two conditions: $2R/A \leq R/2 $ and $ 4\rho^{\gamma, \beta}_{x_0}(R/A) \leq \rho^{\gamma, \beta}_{x_0}(R)$. The first condition is automatically verified by \eqref{expression.factor}, hence we only need to verify the latter, which easily follows by the following estimates:
\begin{equation*}\begin{split}
4   \rho_{x_0}^{\gamma,\beta}(R/A) \leq 4 \ka_{18} \frac{R^2}{A^2}\left[\frac{R}{A}\vee |x_0|\right]^{-\left(\gamma-\beta\right)}\leq \frac{4 \ka_{18}^2}{A^\sigma}\ka_{18}^{-1}R^2\left[\frac{R\vee |x_0|}{A}\right]^{-\left(\gamma-\beta\right)}\le \rho^{\gamma, \beta}_{x_0}(R)\,,
\end{split}\end{equation*}
which follow from $\frac{R\vee |x_0|}{A} \leq \frac{R}{A}\vee |x_0|$ together with the condition \eqref{expression.factor}. The proof is concluded.\qed

\subsubsection{Further estimates on test functions}
The operator $\mathcal{L}_{\gamma, \beta}f=|x|^\gamma\nabla\cdot\left(|x|^{-\beta}\nabla f\right)$ acts on smooth functions as follows:
\begin{equation} \label{operator_on_cutoff_functions}
        \mathcal{L}_{\gamma, \beta}\left( \phi \right) =  |x|^{\gamma-\beta} \left[ \Delta \phi -\beta \frac{x}{ |x|^2}\cdot \nabla \phi \right].
\end{equation}
In the proof of Proposition \ref{Herrero_Pierre_Thm}, we use the following technical Lemma.
\begin{lem}\label{test.function.HP}
For any $x_0 \in \RR^N$ and any $R > 0$ there exists $\phi \in C^2_c(\RR^N)$ such that $\supp(\phi) \subset B_{2R}(x_0)$, $\phi \equiv 1 $ on $B_R(x_0)$, $0\le \phi\le 1$ and the following estimate holds:
\begin{equation}\label{estimate.testfunction}
\phi^{\frac{-m}{1-m}}(x) |\mathcal{L}_{\gamma, \beta}(\phi)(x) |^{\frac{1}{1-m}} \le  \ka_{10}  \left(\rho^{\gamma,\beta}_{x_0}(R)\right)^{-\frac{1}{1-m}}\qquad\mbox{for all }x\in \RR^N\,,
\end{equation}
where $\ka_{10}>0$ depends only on $N, \gamma, \beta$ and $m$.
\end{lem}
\noindent\textbf{Proof. }We define a function $ \phi:= \psi\left(  |x-x_0|^{\sigma}  R^{-\sigma} \right)^b $ with $b>0$ to be chosen later; we choose the cutoff function $ \psi: \left[0, \infty \right) \rightarrow [0,1] $ to be smooth. A simple calculation shows that
\begin{equation}\label{derivatives_of_cutoff_functions}
    \begin{aligned}
&\nabla\phi=  b\sigma R^{-\sigma} \psi \left( |x-x_0|^\sigma  R^{-\sigma} \right)^{b-1}
        \psi'\left(  |x-x_0|^\sigma R^{-\sigma} \right)  |x-x_0|^{\sigma-2} (x-x_0) , \\ \nonumber
&\Delta\phi
        =  b \sigma R^{-\sigma} \psi^{b-2}\, |x-x_0|^{\sigma-2} \left[ \sigma |x-x_0|^\sigma R^{-\sigma} \left((b-1) |\psi'|^2 + \psi\psi''  \right)+\psi \psi' (N+\beta-\gamma) \right]. \\
    \end{aligned}
\end{equation}
Using the expression \eqref{operator_on_cutoff_functions} we get
 \begin{equation*}\label{what_inside_the_integral}
 \begin{split}
   \phi^{\frac{-m}{1-m}}& \left| \mathcal{L}_{\gamma, \beta}\left( \phi \right) \right|^{\frac{1}{1-m}} =   \psi^{ \frac{-bm}{1-m}} \left| \mathcal{L}_{\gamma, \beta}\left( \psi^b \right) \right|^{\frac{1}{1-m}}
    =   \psi^{\frac{-bm+\left( b-2 \right)}{1-m}}  \left[ b \sigma R^{-\sigma}|x|^{2-\sigma}|x-x_0|^{\sigma-2}\right]^{\frac{1}{1-m}} \\
    & \times
    \left| \sigma |x-x_0|^\sigma R^{-\sigma} \left((b-1) |\psi'|^2 + \psi\psi''  \right)+\psi \psi' (N+\beta-\gamma)-\psi \psi' \beta |x|^{-2}(x-x_0)\cdot x \right|^{\frac{1}{1-m}}.\\
 \end{split}
\end{equation*}
We need to split the proof in two cases, depending on the relation between $|x_0|$ and $R$. \\
\noindent$\bullet~$\textit{When $0 \le |x_0| \le \frac{3}{2}R$\,: }Choosing $\psi= \psi(|x-x_0|) $ to be a equal to $1$ on $B_{(7/4)R}(x_0)$ and zero outside $B_{2R}(x_0)$, we have $\supp(\mathcal{L}_{\gamma, \beta}(\psi^b))\subseteq B_{(7/4)R}(x_0)^c\cap B_{2R}(x_0) $; since $B_{(1/4)R}(0) \subset B_{(7/4)R}(x_0)$, it turns out that $\supp(\mathcal{L}_{\gamma, \beta}(\psi^b))\subseteq \{(1/4)R \le |x| \le 3R \} \cap \{ (7/4)R \le |x-x_0| \le 2R \}$. Taking $ b \geq \frac{2}{1-m} $ we obtain
 \begin{equation*}
 \begin{split}
    \phi^{\frac{-m}{1-m}} \left| \mathcal{L}_{\gamma, \beta}\left( \phi \right) \right|^{\frac{1}{1-m}}
    \le C_{\sigma} \left[ b \sigma  R^{-\sigma}\right]^{\frac{1}{1-m}} \left[\sigma4^\sigma(|b-1||\psi'|^2+|\psi''|) + |\psi'|(N+\beta-\gamma + 4|\beta| )\right]^{\frac{1}{1-m}}
 \end{split},
 \end{equation*}
where we have used that $0 \le \psi \le 1$; note that $C_\sigma>0$  depends only on $\sigma$ and,  by \eqref{technical.inequality.measures} in this case  $\rho^{\gamma, \beta}_{x_0}(R)\asymp R^\sigma$,  this proves \eqref{estimate.testfunction}.  \\
\noindent$\bullet~$\textit{When $|x_0| \ge \frac{3}{2}R$\,:} In this case we choose $\psi(|x-x_0|)$ equal to $1$ on $B_R(x_0)$ and equal to $0$ outside $B_{(5/4) R}(x_0)$.
In this way, $\supp(\mathcal{L}_{\gamma, \beta}(\psi^b))\subseteq\{ R \le |x-x_0| \le (5/4)R \}\subseteq\{(1/3)|x_0|\leq |x|\leq(11/6)|x_0|\}$; noticing that $\{ R \le |x-x_0| \le (5/4)R \} \subset \{|x|\geq R/4 \}$,  using \eqref{technical.inequality.measures} and proceeding as in the previous case, we conclude the proof of \eqref{estimate.testfunction}.\qed

\newpage

\noindent {\large \sc Acknowledgments. }This work was partially funded by Projects MTM2014-52240-P and MTM2017-85757-P (Spain), and by the E.U. H2020 MSCA programme, grant agreement 777822. N.~S. was partially funded by the FPI-grant BES-2015-072962, associated to the project MTM2014-52240-P (Spain). We would like to warmly thank the anonymous reviewer for an extraordinary work that helped us a lot to improve the paper to its current status.

\smallskip\noindent {\sl\small\copyright~2018 by the authors. This paper may be reproduced, in its entirety, for non-commercial purposes.}



\addcontentsline{toc}{section}{~~~References}
\begin{thebibliography}{00} \itemsep1pt \parskip2pt \parsep0pt

\small

\bibitem{AP} B. Abdellaoui, A. I. Peral. \textit{H\"older regularity and Harnack inequality for degenerate parabolic equations related to Caffarelli-Kohn-Nirenberg inequalities. }
Nonlinear Anal. \textbf{57} \rm (2004), 971--1003.




\bibitem{AB-JDE}D.G. Aronson, P. Besala, \textit{Parabolic equations with unbounded coefficients.}
J. Differential Equations \textbf{3} \rm (1967), 1--14.

\bibitem{AB-CM} D.G. Aronson, P. Besala, \textit{Uniqueness of positive solutions of parabolic equations with unbounded coefficients.}
Colloq. Math. \textbf{18} \rm (1967), 125--135.


\bibitem{AS-ARMA67} D. G. Aronson, J. Serrin, \textit{Local behavior of solutions of quasilinear parabolic equations},
{ \rm Arch. Rational Mech. Anal.}  \textbf{25} \rm (1967), 81--122.

\bibitem{BoGa}L. Boccardo, T. Gallou\"et. \textit{Nonlinear elliptic and parabolic equations involving measure data. }J. Funct. Anal. {\bf 87} \rm (1989), 149--169.

\bibitem{BGInv72} { \rm E. Bombieri, E. Giusti}, \textit{Harnack's inequality for elliptic differential equations on minimal surfaces},
Invent. Math. \textbf{15} \rm (1972), 24--46.

\bibitem{BS2000} S. Bonafede, I. Skrypnik, \textit{On Hölder continuity of solutions of doubly nonlinear parabolic equations with weight. (English, Ukrainian summary)}
Ukraïn. Mat. Zh. \textbf{51} \rm (1999), 890--903; translation in
Ukrainian Math. J. \textbf{51} \rm (1999), 996–1012.


 \bibitem{BDMN2016b} M.~Bonforte, J.~Dolbeault, M.~Muratori,  B.~Nazaret, \textit{Weighted fast diffusion equations ({P}art {II}): {S}harp asymptotic
  rates of convergence in relative error by entropy methods},  Kin. Rel. Mod. \textbf{10} \rm (2017), 61--91.


 \bibitem{BDMN2016a} M.~Bonforte, J.~Dolbeault, M.~Muratori,  B.~Nazaret, \textit{Weighted fast diffusion equations ({P}art {I}): {S}harp asymptotic
  rates without symmetry and symmetry breaking in {C}affarelli-{K}ohn-{N}irenberg inequalities}, Kin. Rel. Mod. \textbf{10} \rm (2017), 33--59.

\bibitem{BFR} M. Bonforte, A. Figalli, X. Ros-Oton, \textit{Infinite speed of propagation and regularity of solutions to the fractional porous medium equation in general domains}, Comm. Pure Appl. Math. \textbf{70}, \rm (2017), 1472--1508.

\bibitem{BFV-Parabolic} M. Bonforte, A. Figalli, J. L. V{\'a}zquez. \textit{Sharp global estimates for local and nonlocal  porous medium-type equations in bounded domains. } Anal. PDE \textbf{11} \rm (2018), 945--982.

\bibitem{BGV-JMPA} {\rm M. Bonforte, G.~Grillo, J. L. V{\'a}zquez},
    {\em Behaviour near extinction  for the Fast Diffusion Equation on bounded domains}, J. Math. Pures Appl. \textbf{97}, \rm (2012), 1--38.

\bibitem{BGV-JEE} {\rm M. Bonforte, G. Grillo, J.L V\'azquez.} {\it
Fast diffusion flow on manifolds of nonpositive curvature.}
\textrm{ J. Evol. Equ.} \bf 8 \rm(2008),  99--128.

\bibitem{BGV-ARMA} M. Bonforte, G. Grillo, J. L. V{\'a}zquez, \textit{Special fast diffusion
    with slow asymptotics. Entropy method and flow on a Riemannian manifold},
    {\rm Arch. Rat. Mech. Anal.}  \bf 196, \rm (2010), 631--680.

\bibitem{BGV} {\rm M. Bonforte, G. Grillo, J.L. V\'azquez}.  \textit{Quantitative local bounds for subcritical semilinear elliptic equations},
Milan J. Math. \bf 80 (\rm 2012), 65--118.


\bibitem{BV-ADV-plap}M. Bonforte, R. G. Iagar, J. L. V{\'a}zquez, \emph{Local smoothing effects, positivity, and Harnack inequalities for the fast $p$\,-Laplacian equation}, Advances in Math. \bf 224, \rm (2010),  2151--2215.

\bibitem{BV-ADV}  {\rm M.  Bonforte, J.~L. V\'azquez. } {\it Positivity, local smoothing, and Harnack inequalities for very fast diffusion equations}, \textrm{Advances in Math.} \bf 223 \rm (2010), 529--578.

\bibitem{BV-PPR1} M. Bonforte,  J. L. V\'azquez, {\it A Priori Estimates  for Fractional Nonlinear  Degenerate Diffusion Equations on bounded domains},  Arch. Rat. Mech. Anal. \textbf{218} \rm (2015),  317--362.

\bibitem{XC-XRS} {\rm X. Cabr\'e, X. Ros-Oton,} {\it Sobolev and isoperimetric inequalities with monomial weights},
J. Differential Equations \textbf{255} \rm (2013),  4312--4336.


  \bibitem{Caffarelli-Kohn-Nirenberg-84} L.~Caffarelli, R.~Kohn, L.~Nirenberg, \textit{First order interpolation inequalities with weights},
 Compositio Math. \textbf{53} \rm (1984), 259--275.


\bibitem{Catrina-Wang-01} F.~Catrina,  Z.-Q. Wang, \textit{On the {C}affarelli-{K}ohn-{N}irenberg inequalities: sharp constants,
  existence (and nonexistence), and symmetry of extremal functions}, Comm. Pure Appl. Math. \textbf{54} \rm (2001), 229--258.

\bibitem{CF-AA84} F. Chiarenza, M. Frasca,
\textit{Boundedness for the solutions of a degenerate parabolic equation},
Applicable Anal. \textbf{17} \rm (1984),  243--261.


\bibitem{CS-AA87} F. Chiarenza, R. Serapioni,
\textit{Pointwise estimates for degenerate parabolic equations},
Appl. Anal. \textbf{23} \rm (1987), 287--299.

\bibitem{CS-RSMUP85} F. Chiarenza, R. Serapioni,
\textit{A remark on a Harnack inequality for degenerate parabolic equations},
Rend. Sem. Mat. Univ. Padova \textbf{73} \rm (1985), 179--190.

\bibitem{CS-AMPA84} F. Chiarenza, R. Serapioni,
\textit{Degenerate parabolic equations and Harnack inequality},
Ann. Mat. Pura Appl. \textbf{137} \rm (1984), 139--162.

\bibitem{CS-CPDE84} F. Chiarenza, R. Serapioni,
\textit{A Harnack inequality for degenerate parabolic equations},
Comm. Partial Differential Equations \textbf{9} \rm (1984), 719--749.

\bibitem{CP-JFA82} M. Crandall, M. Pierre, \textit{Regularizing effects for $u_t+A \phi(u)=0$ in $\LL^1$},
J. Funct. Anal. \textbf{45} \rm (1982), 194--212.

\bibitem{DGP-SIAM} A. Dall'Aglio, D. Giachetti, I. Peral.
\textit{Results on parabolic equations related to some Caffarelli-Kohn-Nirenberg inequalities},
SIAM J. Math. Anal. \textbf{36} \rm (2004/05), 691--716.

\bibitem{DPDS-LOG} P. Daskalopoulos, M. del Pino.
\textit{On the Cauchy problem for $u_t=\Delta\log u$ in higher dimensions}, Math. Ann. \textbf{313} (2), 189--206.

\bibitem{DaskaBook}{\rm P. Daskalopoulos, C. E. Kenig}. \textsl{``Degenerate diffusions. Initial value problems and local regularity theory''. }EMS Tracts in Mathematics, \textbf{1}. European Mathematical Society (EMS), Zürich, \rm 2007. x+198 pp. ISBN: 978-3-03719-033-3.


\bibitem{DBbook} {\rm E. DiBenedetto}. {\sl ``Degenerate parabolic
equations'',} Universitext. Springer-Verlag, New York, 1993.

\bibitem{DiBenGianVesBook} {\rm E. DiBenedetto, U. Gianazza, V. Vespri.} \textsl{``Harnack's inequality for degenerate and singular parabolic equations''.} Springer Monographs in Mathematics. Springer, New York, \rm 2012. xiv+278 pp. ISBN: 978-1-4614-1583-1.

\bibitem{DGV} {\rm E. DiBenedetto, U. Gianazza, V. Vespri,}
\textit{Forward, Backward and Elliptic Harnack Inequalities for Non–Negative Solutions to Certain Singular Parabolic Partial
Differential Equations}, Ann. Sc. Norm. Super. Pisa Cl. Sci. (5) \textbf{9} \rm (2010), 385--422.


\bibitem{DiBenedettourbVesp}  {\rm E. DiBenedetto, M. Urbano, V.
Vespri,}  \textit{Current issues on singular and degenerate evolution
equations}, {\em  Evolutionary equations. Vol. I,}  169--286,
Handb. Differ. Equ., North-Holland, Amsterdam, 2004.


\bibitem{1406} J.~Dolbeault, M.~J. Esteban, S.~Filippas, A.~Tertikas, \textit{Rigidity results with applications to best constants and symmetry of
  {C}affarelli-{K}ohn-{N}irenberg and logarithmic {H}ardy inequalities}, Calc. Var. Partial Differential Equations, \textbf{54} \rm (2015),
  2465--2481.

\bibitem{DEL2015} J.~Dolbeault, M.~J. Esteban, M.~Loss, \textit{Rigidity versus symmetry breaking via nonlinear flows on cylinders
  and {E}uclidean spaces}, Invent. Math. \textbf{206} \rm (2016), 397--440.

\bibitem{DELM2015} J.~Dolbeault, M.~J. Esteban, M.~Loss, M.~Muratori, \textit{Symmetry for extremal functions in subcritical
  {C}affarelli-{K}ohn-{N}irenberg inequalities}, C. R. Math. Acad. Sci. Paris \textbf{355} \rm (2017),  133--154.

\bibitem{DETT} J.~Dolbeault, M.~J. Esteban, G.~Tarantello, A.~Tertikas, \textit{Radial symmetry and symmetry breaking for some interpolation
  inequalities}, Calc. Var. Partial Differential Equations, \textbf{42} \rm (2011), 461--485.

\bibitem{DMN2017} J.~Dolbeault, M.~Muratori, B.~Nazaret, \textit{Weighted interpolation inequalities: a perturbation approach.}
Math. Ann. \textbf{369} \rm (2017), 1237–1270.

\bibitem{Eidus} D. Eidus,
\textit{The Cauchy problem for the nonlinear filtration equation in an inhomogeneous medium},
J. Differential Equations, \textbf{84} \rm (1990), 309--318.

\bibitem{KE-PAMS} D. Eidus, S. Kamin,
\textit{The filtration equation in a class of functions decreasing at infinity},
Proc. Amer. Math. Soc. \textbf{120} \rm (1994), 825--830.



\bibitem{FKS} {\rm E. B. Fabes, C. E. Kenig, R. P. Serapioni}.  \textit{The local regularity of solutions of degenerate elliptic equations},
Comm. Partial Differential Equations  \textbf{7}  \rm (1982),  77--116.

\bibitem{FGS86} E. B. Fabes, N. Garofalo, S. Salsa. \textit{A backward Harnack inequality and Fatou theorem for nonnegative solutions of parabolic equations}, Illinois J. Math. \textbf{30}  \rm (1986), 536--565.

\bibitem{Felli-Schneider-03} V.~Felli, M.~Schneider, \textit{Perturbation results of critical elliptic equations of
  {C}affarelli-{K}ohn-{N}irenberg type}, J. Differential Equations, \textbf{191} \rm (2003), 121--142.


\bibitem{FGW}{ \rm B. Franchi, C. Guti\`{e}rrez, R. L. Wheeden}
\textit{Weighted Sobolev-Poincaré inequalities for Grushin type operators},
\textrm{Comm. Partial Differential Equations} \bf 19 \rm (1994), 523--604.

\bibitem{FG-PAMS} E. Fabes, N. Garofalo \textit{Parabolic B.M.O. and Harnack’s inequality}, Proc.Am. Math. Soc. \textbf{50},
\rm (1985) 63--69.


\bibitem{Giusti} E.~Giusti, \textsl{``Direct Methods in the Calculus of Variations''}, World Scientific Publishing Co., Inc., River Edge, NJ, 2003.

\bibitem{GM14} G.~Grillo, M.~Muratori, \textit{Radial fast diffusion on the hyperbolic space}, Proc. Lond. Math. Soc. (3), \textbf{109} \rm  (2014), 283--317.



\bibitem{GMPo13} G.~Grillo, M.~Muratori, M.~M. Porzio, \textit{Porous media equations with two weights: smoothing and decay
  properties of energy solutions via {P}oincar\'e inequalities}, Discrete Contin. Dyn. Syst., \textbf{33} \rm (2013) 3599--3640.


\bibitem{GMP15} G.~Grillo, M.~Muratori, F.~Punzo, \textit{On the asymptotic behaviour of solutions to the fractional porous
  medium equation with variable density}, Discrete Contin. Dyn. Syst., \textbf{35} \rm (2015), 5927--5962.


\bibitem{GMV-manifolds}G.~Grillo, M.~Muratori, J. L. V\'azquez, \textit{The porous medium equation on Riemannian manifolds with negative curvature. The large-time behaviour.} Advances in Mathematics \textbf{314} \rm (2017), 328--377.


\bibitem{GN-CPDE88} C. E. Guti\'errez, G. S. Nelson, \textit{Bounds for the fundamental solution of degenerate parabolic equations},
Comm. Partial Differential Equations \textbf{13} \rm (1988), 635--649.


\bibitem{GW2} C. Gutiérrez, R. L. Wheeden, \textit{Mean value and Harnack inequalities for degenerate parabolic equations. }Colloq. Math. \textbf{60/61} \rm     (1990), 157--194.

\bibitem{GW} C. Gutiérrez, R. L. Wheeden, \textit{Harnack's inequality for degenerate parabolic equations.} Comm. Partial Differential Equations \textbf{16}  \rm (1991), 745--770.


\bibitem{Hadamard}J. Hadamard, \textit{Extension {\'a} l'{\'e}quation de la chaleur d'un th{\'e}or{\'e}me de A. Harnack. }Rend.
Circ. Mat. Palermo, \textbf{3} \rm (1954), 337--346.


\bibitem{HKM} { \rm J. Heinonen,T. Kilpeläinen, O. Martio}. \textsl{``Nonlinear potential theory of degenerate elliptic equations''.}
Oxford Mathematical Monographs. Oxford Science Publications. The Clarendon Press, Oxford University Press, New York, 1993. vi+363 pp. ISBN: 0-19-853669-0

\bibitem{HK} {\rm J. Heinonen, P. Koskela}
\textit{Weighted Sobolev and Poincaré inequalities and quasiregular mappings of polynomial type},
\textrm{Math. Scand.}  \bf 77  \rm (1995), 251--271.

\bibitem{HerreroPierre} {\rm M.A. Herrero, M. Pierre}. \textit{The Cauchy problem for $u_t=\Delta u^m$ when $0<m<1$}, Trans. Amer. Math. Soc. \textbf{291} \rm \rm (1985), 145--158.

\bibitem{IS14} R.~G. Iagar, A.~S{\'a}nchez, \textit{Large time behavior for a porous medium equation in a nonhomogeneous
  medium with critical density}, Nonlinear Anal. , \textbf{102} \rm (2014), 226--241.



\bibitem{I-NM} K. Ishige, \textit{On the behavior of the solutions of degenerate parabolic equations.}
Nagoya Math. J. \textbf{155} \rm (1999), 1--26.


\bibitem{IM-MZ} K. Ishige, M. Murata, \textit{An intrinsic metric approach to uniqueness of the positive Cauchy problem for parabolic equations.}
Math. Z. \textbf{227} \rm (1998), 313--335.





\bibitem{IM-SNS} K. Ishige, M. Murata, \textit{Uniqueness of nonnegative solutions of the Cauchy problem for parabolic equations on manifolds or domains.}
Ann. Scuola Norm. Sup. Pisa Cl. Sci. (4) \textbf{30} \rm (2001), 171--223.




\bibitem{KKT} S. Kamin, R. Kersner, A. Tesei,
\textit{On the Cauchy problem for a class of parabolic equations with variable density},
Atti Accad. Naz. Lincei Cl. Sci. Fis. Mat. Natur. Rend. Lincei  Mat. Appl. \textbf{9} \rm (1998), 279--298.


\bibitem{KRV-DCDS}  S. Kamin, G. Reyes, J.L. Vázquez, \textit{Long time behavior for the inhomogeneous PME in a medium with rapidly decaying density},
Discrete Contin. Dyn. Syst. \textbf{26} \rm (2010), 521--549.

\bibitem{KR-CPAM81} S. Kamin, P. Rosenau,
\textit{Propagation of thermal waves in an inhomogeneous medium}.
Comm. Pure Appl. Math. \textbf{34} \rm (1981), 831--852.

\bibitem{KR-CPAM82} S. Kamin, P. Rosenau,
\textit{Nonlinear diffusion in a finite mass medium},
Comm. Pure Appl. Math. \textbf{35} \rm (1982), 113--127.

\bibitem{KR-JMPH} S. Kamin, p. Rosenau,
\textit{Nonlinear thermal evolution in an inhomogeneous medium},
J. Math. Phys. \textbf{23} \rm (1982), 1385--1390.


\bibitem{KK-MA} J. Kinnunen, T. Kuusi, \textit{Local behaviour of solutions to doubly nonlinear parabolic equations.}
Math. Ann. \textbf{337} \rm  (2007), 705--728.



\bibitem{KO} {\rm A. Kufner, B. Opic}, \textit{How to define reasonably weighted Sobolev spaces},
Comment. Math. Univ. Carolin. \textbf{25} \rm (1984),  537--554.


\bibitem{KOBook} {\rm A. Kufner, B. Opic,} {\sl ``Hardy-type inequalities".}
Pitman Research Notes in Mathematics Series, 219. Longman Scientific \& Technical, Harlow, 1990. xii+333 pp. ISBN: 0-582-05198-3.


\bibitem{K-SNS} T. Kuusi, \textit{Harnack estimates for weak supersolutions to nonlinear degenerate parabolic equations.}
Ann. Sc. Norm. Super. Pisa Cl. Sci. (5) \textbf{7} \rm (2008), 673--716.


\bibitem{MiKuu} T. Kuusi, G.  Mingione. {\it Guide to Nonlinear Potential Estimates. } Bulletin of Math. Sci. {\bf 4} (2014), 1--82.


\bibitem{MadernaSalsa} C. Maderna, S. Salsa, \textit{Sharp estimates of solutions to a certain type of singular elliptic boundary value problems in two dimensions}, Applicable Anal. \textbf{12} \rm (1981), 307--321.


 \bibitem{Moser}  J. Moser, \textit{A Harnack inequality for parabolic differential equations}, Comm. Pure and Appl. Math.  \textbf{17} (1964), 101--134. And \textit{Correction to: ``A Harnack inequality for parabolic differential equations''},
Comm. Pure Appl. Math. \textbf{20} \rm (1967), 231--236.

\bibitem{MoserCpam71} {\rm J.Moser}, \textit{On a pointwise estimate for parabolic differential equations},
Comm. Pure Appl. Math. \textbf{24}  \rm (1971), 727--740.


\bibitem{NR13} S.~Nieto, G.~Reyes, \textit{Asymptotic behavior of the solutions of the inhomogeneous porous
  medium equation with critical vanishing density}, Commun. Pure Appl. Anal. , \textbf{12} \rm (2013), 1123--1139.




\bibitem{NPS-JDE2015} K. Nyström, H. Persson, O. Sande, \textit{Boundary estimates for solutions to linear degenerate parabolic equations},
J. Differential Equations \textbf{259} \rm (2015), 3577--3614.

\bibitem{Pinchover} Y. Pinchover, \textit{On uniqueness and nonuniqueness of the positive Cauchy problem for parabolic equations with unbounded coefficient.}
Math. Z. \textbf{223} \rm (1996), 569--586.



\bibitem{Pini} B. Pini, \textit{Sulla soluzione generalizzata di Wiener per il primo problema di valori al contorno
nel caso parabolico. }Rend. Sem. Mat. Univ. Padova, \textbf{23} \rm (1954), 422–434.



\bibitem{RV-CPAA09} G. Reyes, J.L. Vázquez,
\textit{Long time behavior for the inhomogeneous PME in a medium with slowly decaying density},
Commun. Pure Appl. Anal. \textbf{8} \rm (2009), 493--508.


\bibitem{RV-CPAA08} G. Reyes, J.L. Vázquez,
\textit{The inhomogeneous PME in several space dimensions. Existence and uniqueness of finite energy solution},
Commun. Pure Appl. Anal. \textbf{7} \rm (2008), 1275--1294.

\bibitem{RV-JEMS06} G. Reyes, J.L. Vázquez,
\textit{A weighted symmetrization for nonlinear elliptic and parabolic equations in inhomogeneous media},
J. Eur. Math. Soc. (JEMS) \textbf{8} \rm (2006), 531--554.

\bibitem{RV-NHM} G. Reyes, J.L. Vázquez,
\textit{The Cauchy problem for the inhomogeneous porous medium equation},
Netw. Heterog. Media \textbf{1} \rm (2006), 337--351.


\bibitem{RV-JMPA} G. Reyes, J.L. Vázquez,
\textit{Asymptotic behaviour of a generalized Burgers' equation},
J. Math. Pures Appl. \textbf{9} \rm (1999),  633--666.

\bibitem{SY} M. V. Safonov, Y. Yuan, \textit{Doubling properties for second order parabolic equations}, Ann. of Math. (2) \textbf{150} \rm (1999), 313--327.

\bibitem{SW} E. Sawyer, R. L. Wheeden, \textit{Weighted inequalities for fractional integrals on Euclidean and homogeneous spaces. }Amer. J. Math. \bf 114 \rm (1992), 813--874.

\bibitem{Surnachev-JDE}{\rm M. Surnachev}
\textit{A Harnack inequality for weighted degenerate parabolic equations},
J. Differential Equations \bf 248 \rm (2010), 2092--2129.

\bibitem{Surnachev-TMM}{\rm M. Surnachev}
\textit{Regularity of solutions of parabolic equations with a double nonlinearity and a weight},
Trans. Moscow Math. Soc.  \textbf{75} \rm (2014), 259--280.


\bibitem{T-CPAM68} S. N. Trudinger, \textit{Pointwise estimates and quasilinear parabolic equations},
\textrm{Comm. Pure Appl. Math.} \textbf{21}  \rm (1968), 205--226.



\bibitem{JLVSmoothing}{\rm J.~L. V\'azquez}. {\sl ``Smoothing and Decay Estimates for Nonlinear Diffusion Equations. Equations of Porous Medium Type''.} Oxford Lecture Series in Mathematics and its Applications, 33. Oxford University Press, Oxford, 2006.



\bibitem{JLVPorousMedioum}{\rm J.L.  V\'azquez}. \textsl{``The porous medium equation. Mathematical theory.''}
Oxford Mathematical Monographs. The Clarendon Press, Oxford University Press, Oxford, \rm 2007. xxii+624 pp. ISBN: 978-0-19-856903-9; 0-19-856903-3.

\bibitem{JLV-Hyp}{\rm J.L.  V\'azquez.} \textit{Fundamental solution and long time behavior of the porous medium equation in hyperbolic space}, J. Math. Pures Appl. (9) \textbf{104} \rm (2015), 454--484.

\bibitem{WNC-NonAnal} Y. Wang, P. Niu, X. Cui, \textit{Harnack estimates for a quasi-linear parabolic equation with a singular weight.}
Nonlinear Anal. \textbf{74} \rm (2011),6265--6286.



\end{thebibliography}
\end{document}